\documentclass[amssymb,amscd,amsfonts,refcheck,11pt,graphicx,verbatim,righttag]{amsart}
\usepackage{amsfonts}
\usepackage{amssymb}

 \numberwithin{equation}{section}
 \allowdisplaybreaks[4]
\theoremstyle{plain}

\newcommand{\N}{{\mathbb N}}

\newcommand{\Z}{{\mathbb Z}}
\newcommand{\R}{{\mathbb R}}

\newcommand{\Q}{{\mathbb Q}}
\newcommand{\C}{{\mathcal C}}
\newcommand{\D}{{\mathcal D}}

\newcommand{\M}{{\mathcal M}}

\def\eps{\varepsilon}
\newtheorem{Thm}{Theorem}[section]
\newtheorem{thm}[Thm]{Theorem}
\newtheorem{cor}{Corollary}[section]
\newtheorem{coro}[cor]{Corollary}
\newtheorem{pro}{Proposition}[section]
\newtheorem{prop}[pro]{Proposition}
\newtheorem{Rem}{Remark}[section]
\newtheorem{rem}[Rem]{Remark}
\newtheorem{Lem}{Lemma}[section]
\newtheorem{lem}[Lem]{Lemma}

\def \proof{\bigbreak\noindent{\it Proof.~~}}

 \numberwithin{equation}{section}
\begin{document}
%Topmatter
\title{Multifractal analysis and localized asymptotic behavior for almost additive potentials}

\author{Julien Barral}
\address{LAGA (UMR 7539), D\'epartement de Math\'ematiques, Institut Galil\'ee, Universit\'e
 Paris 13, 99 avenue Jean-Baptiste Cl\'ement , 93430  Villetaneuse, France}
\email{barral@math.univ-paris13.fr}
\author{Yan-Hui Qu}
\address{Department of Mathematics, Tsinghua University, Beijing 100084, China
}
\email{jyh02@mails.tsinghua.edu.cn}

\thanks{The authors thank De-Jun Feng for valuable discussions}
\begin{abstract}
We conduct the multifractal analysis of the level sets of the asymptotic behavior of
almost-additive continuous potentials $(\phi_n)_{n=1}^\infty$ on a
topologically mixing subshift of finite type $X$ endowed itself with
a metric associated with such a potential. We work without bounded
distorsion property assumption. We express the whole Hausdorff
spectrum in terms of  a conditional variational principle, as well as a new large deviations principle. Our approach provides
 a new description of the structure of the spectrum in terms of {\it weak} concavity. Another new point is that we consider sets of points at which the asymptotic behavior of $\phi_n(x)$ is localized, i.e.
  depends on the point $x$ rather than being equal to a constant. Specifically, we compute the Hausdorff dimension of sets of the form
   $\{x\in X: \lim_{n\to\infty} \phi_n(x)/n=\xi(x)\}$, where $\xi$ is a given continuous function.
   This is naturally related to Birkhoff's ergodic theorem and  has interesting geometric applications to fixed points in the asymptotic average for dynamical systems in $\R^d$,
    as well as the fine local behavior of the harmonic measure on conformal planar Cantor sets.
\end{abstract}

\maketitle

\section{Introduction} We say that $(X,T)$ is a {\it topological dynamical system}
(TDS) if $X$ is a compact metric space and $T$ is a continuous
mapping from $X$ to itself. We denote by $\M(X,T)$ the set of
invariant probability measures on~$(X,T)$.

We say that $\Phi=(\phi_n)_{n=1}^\infty$ is {\it almost additive} if
$\phi_n$ is continuous on $X$   and there is a positive constant
$C(\Phi)>0$ such that
$$
-C(\Phi)+ \phi_n+\phi_p\circ T^n\leq \phi_{n+p}\leq C(\Phi)+
\phi_n+\phi_p\circ T^n,\quad \forall\,n,p\in \N.
$$
By subadditivity, for every $\mu\in \M(X,T)$, $\displaystyle
\Phi_*(\mu):=\lim_{n\to\infty} \int_X \frac{\phi_n}{n}\,
\text{d}\mu$ exists, and we define the compact  convex set $
L_\Phi=\{ \Phi_*(\mu):\mu\in \M(X,T)\}. $ We denote by
$\C_{aa}(X,T)$ the collection of almost-additive potentials on~$X$.

The ergodic theorem naturally raises the following question. Given
$\Phi$ an almost additive potential taking values in $\R^d$ (this
means that $\Phi=(\Phi^1,\cdots,\Phi^d)$ with each $\Phi^i\in
\C_{aa}(X,T)$) and $\xi:X\to\R^d$ a continuous function, what is the
Hausdorff dimension of the set
$$
E_\Phi(\xi):=\Big \{x\in X:
\lim_{n\to\infty}\frac{\phi_n(x)}{n}=\xi(x)\Big \}?
$$
When $\xi(x)\equiv\alpha$ is constant, this question has been solved
for some $C^{1+\varepsilon}$ conformal dynamical systems, sometimes assuming restrictions on the regularity of $\Phi$, and this problem is known as
the multifractal analysis of Birkhoff averages, and more generally
almost additive potentials \cite{Collet,Rand,PeW,PW,P, O,FF,BS,FFW,BSS,FLW,FeLa02,Fen03,Olsen,BD}. Moreover,
 the optimal results are expressed in terms of a variational principle of the following form:
 $E_\Phi(\alpha)\neq \emptyset$ if and only if $\alpha\in L_\Phi$ and in this case
\begin{equation}\label{varprinc}
\dim_H E_\Phi(\alpha)=\max\left \{\frac{h_\mu(T)}{\int_X \log \|DT\|
\, \text{d}\mu}: \mu\in \M(X,T),\  \Phi_*(\mu)=\alpha\right \},
\end{equation}
the supremum being attained by a unique Gibbs measure if $\Phi$ is
the sequence of Birkhoff sums of a H\"older potential, and $\alpha$
is in the  interior of $L_\Phi$. To our best knowledge no result is
known for $\dim_H  E_\Phi(\xi)$ for non constant $\xi$. We are going
to give an answer to this question when $(X,T)$ is a topologically
mixing subshift of finite type endowed with a metric associated with
a negative almost additive potential, and consider geometric
realizations on Moran sets like those studied in \cite{Bar96}, the
main examples being  $C^1$ conformal repellers and $C^1$ conformal
iterated function systems (see section~\ref{examples} for precise
definitions and statements). In the setting outlined above, if $d=1$
and $\xi$ takes its values in $L_\Phi$, we find the natural
variational formula
$$
\dim_H  E_\Phi(\xi)=\max\left \{\frac{h_\mu(T)}{\int_X \log \|DT\|
\, \text{d}\mu}: \mu\in \M(X,T),\  \Phi_*(\mu)\in \xi(X)\right\}.
$$
As an application of this kind of results, we obtain unexpected results like the following one:
 Let $d\in\N_+$ and $(m_1,\dots,m_d)$ be $d$ integers $\ge 2$. Let $T:[0,1]^d\to[0,1]^d$ be the mapping
  $(x_1,\dots,x_d)\mapsto (m_1x_1\pmod 1,\dots,m_d x_d\pmod 1)$. Consider
$$
\mathcal F=\Big\{x\in [0,1]^d: \lim_{n\to\infty}
\frac{1}{n}\sum_{k=0}^{n-1} T^k x=x\Big \},
$$
the set of those points $x$ which are fixed by $T$ in the asymptotic average.
Then $\mathcal F$ is dense and of full Hausdorff dimension in $[0,1]^d$.

Another application concerns harmonic measure. Let us consider here the special
case of the set $J=C^2\subset \R^2$, where $C$ is the middle third Cantor set.
The harmonic measure on $J$ is the probability measure $\omega$ such that for each
 $x\in J$ and $r>0$, $\omega (B(x,r))$ is the probability that a planar Brownian
  motion started at $\infty$ attains $J$ for the first time at a point of $B(x,r)$
  (see Section \ref{appli2} for more general examples and a reference). For $x\in J$,
   one defines the local dimension of $\omega$ at $x$ as
$\displaystyle d_\omega (x)=\lim_{r\to 0^+} {\log \omega
(B(x,r))}/{\log (r)} $ whenever this limit exists. Let $I$ stand for
the set of all possible local dimensions
 for $\omega$. By using the fact that $\omega$ is a Gibbs measure,  we prove that
  if $\xi:J\to\R_+$ is continuous and $\xi(J)\subset I$, then the set
  $E_\omega(\xi)=\{x\in J: d_\omega(x)=\xi(x)\}$ is dense in $J$ and the following variational formula holds:
$$
\dim_HE_\omega(\xi)=\sup\{\dim_H
E_\omega(\alpha):\alpha\in\xi(J)\},
\text{ where $ E_\omega(\alpha)=\{x\in J: d_\omega(x)=\alpha\}.$}
$$

Our approach necessitates to revisit the case where $\xi$ is
constant. At this occasion,  we complete the work achieved in
\cite{FF,FFW,FLW} by identifying, in our general framework, the Hausdorff dimensions of the sets
$E_\Phi(\alpha)$ with a large deviation spectrum which is equal to  the Legendre transform of  a kind of "metric" pressure; this is a new kind of large deviation principle in this context.
%
% large principle in our general
%context and on the whole domain $L_\Phi$, while in  \cite{BSS,BD}
%such a result is obtained under strong assumptions on the potentials
%and only in the  interior of $L_\Phi$. Also, we adopt an alternative
%approach to
% that of \cite{FLW} (which only considers additive potentials) regarding the Hausdorff dimensions estimations; 
 Moreover, our  approach brings out an interesting new property for the structure of the Hausdorff spectrum
 $\alpha\mapsto \dim_H E_\Phi(\alpha)$. We call this property {\it weak} concavity;
 it is between concavity and quasi-concavity. This structure turns out to be crucial
 both in establishing the large deviation principle and our results on fixed points in the asymptotic average.

The paper is organized as follows. In Section~\ref{statements} we
give basic definitions and state our main results  on subshift of
finite type.  In Section~\ref{examples} we give the  geometric
realizations. The other sections provide the proofs of our
results.

 %%%%%%%%%%%%%%%%%%%%%%%%%%%%%%%%%%%%%%%%%%%%%%%%%%%%%%%%%%%%%%%%%%%%%%%%%%%%%%%%%%%%%%%%%%%%%%%%%%%%%%%%%%%%%%%%%%%%%%%%%%%%%%%%%%%%
 %%%%%%%%%%%%%%%%%%%%%%%%%%%%%%%%%%%%%%%%%%%%%%%%%%%%%%%%%%%%%%%%%%%%%%%%%%%%%%%%%%%%%%%%%%%%%%%%%%%%%%%%%%%%%%%%%%%%%%%%%%%%%%%%%%%%
 %%%%%%%%%%%%%%%%%%%%%%%%%%%%%%%%%%%%%%%%%%%%%%%%%%%%%%%%%%%%%%%%%%%%%%%%%%%%%%%%%%%%%%%%%%%%%%%%%%%%%%%%%%%%%%%%%%%%%%%%%%%%%%%%%%%%

\section{Definitions and main results}\label{statements}

\subsection{Definitions. Recalls on thermodynamic formalism}

\subsubsection{Thermodynamic formalism for almost additive potentials}

Given $\Phi\in \C_{aa}(X,T),$ define
$\Phi_{\max}:=\max(\phi_{1})+C(\Phi) $ and
$\Phi_{\min}:=\min(\phi_{1})-C(\Phi)$. Define
$\|\Phi\|:=|\Phi_{\max}|\vee|\Phi_{\min}|.$ By the almost additivity property
we easily get
\begin{equation}\label{max-min}
n\Phi_{\min}\leq \phi_n(x)\leq n\Phi_{\max}, \ \ \forall\  n\in\N.
\end{equation}
Consequently we have $\|\phi_n\|_\infty\leq n\|\Phi\|$.

Define two collections of special almost additive potentials on $X$
as
$$
\C_{aa}^+(X,T):=\{\Phi\in \C_{aa}(X,T): \Phi_{\min}>0\}\ \text{ and
}\ \C_{aa}^-(X,T):=\{\Phi\in \C_{aa}(X,T): \Phi_{\max}<0\}.
$$
For $\Phi\in\C_{aa}^-(X,T)$ we get
$
\phi_{n+1}(x)\leq \phi_n(x)+\phi_1(T^nx)+C(\Phi)\leq
\phi_n(x)+\Phi_{\max}<\phi_n(x),
$
So $\{\phi_n:n\in\N\}$ is a strictly decreasing sequence of
functions.

If $\Phi=(\Phi^1,\cdots,\Phi^d)$ is such that each
$\Phi^j\in\C_{aa}(X,T)$, then we call $\Phi$ a {\it vector-valued
almost additive potential} and write $\Phi\in\C_{aa}(X,T,d)$. In
this case $\Phi=(\phi_n)_{n=1}^\infty$ with
$\phi_n=(\phi_n^1,\cdots,\phi_n^d).$ We set $
\Phi_{\max}:=(\Phi^1_{\max},\cdots,\Phi^d_{\max})$ and $\Phi_{\min}:=(\Phi^1_{\min},\cdots,\Phi^d_{\min}). $ Define
$\displaystyle \|\Phi\|:=\Big(\sum_{j=1}^{d}\|\Phi^j\|^2\Big)^{1/2}$ and $\|\Phi\|_{\mbox{\tiny \rm
lim}}:=\limsup_{n\to\infty}\|\phi_n\|_\infty/n. $ We have
 $\|\phi_n\|_\infty\leq n\|\Phi\|$.

%Write  $\Phi\in \C_{aa}^+(X,T,d)$ if $\Phi=(\Phi^1,\cdots,\Phi^d)$
%such that $\Phi_j\in \C_{aa}^+(X,T)$ for $j=1,\cdots,d.$

Given  $u,v\in \R^d$, we write $[u,v]:=\{tu+(1-t)v: 0\leq t\leq 1\}$
to denote the closed interval connecting $u$ and $v$.  If $u_i\leq
v_i$ for $i=1,\cdots,d$, then we write $u\leq v.$ For
$\Phi\in\C_{aa}(X,T,d)$ define
$C(\Phi):=(C(\Phi^1),\cdots,C(\Phi^d))$, then we also have the
following vector version formula:
$$
-C(\Phi)+ \phi_n+\phi_p\circ T^n\leq \phi_{n+p}\leq C(\Phi)+
\phi_n+\phi_p\circ T^n,\quad \forall\quad  n,p\in \N.
$$
For $\mu\in \M(X,T),$ define
$\Phi_\ast(\mu):=(\Phi^1_\ast(\mu),\cdots,\Phi^d_\ast(\mu))$. Define
$L_\Phi:=\{\Phi_\ast(\mu):\mu\in\M(X,T)\}$. Given $\Phi,\Psi\in
\C_{aa}(X,T,d)$, define $\Phi+\Psi:=(\phi_n+\psi_n)_{n=1}^\infty.$
We have $\Phi+\Psi\in \C_{aa}(X,T,d)$ with
$C(\Phi+\Psi)=C(\Phi)+C(\Psi)$.

The simplest almost additive potentials are the additive ones. Given
$\phi:X\to\R^d$ continuous, define
$\phi_n=S_n\phi:=\sum_{j=0}^{n-1}\phi\circ T^j$ and define
$\Phi=(\phi_n)_{n=1}^{\infty}$. In this case
$\phi_{n+p}=\phi_n+\phi_p\circ T^n$, thus $\Phi\in\C_{aa}(X,T,d)$.
Such a $\Phi$   is called an {\it additive potential}. In fact
$\phi_n$ is the $n$-th Birkhoff sum of $\phi.$ Given an additive
potential $\Phi=(S_n\phi)_{n=1}^{\infty}$, if $\phi$ is H\"{o}lder
continuous, we say that $\Phi$ is {\it H\"{o}lder continuous}.
The simplest H\"{o}lder continuous potentials are the constant potentials
$(n\alpha)_{n=1}^\infty$, $\alpha\in\R^d$, that we also denote as~$\alpha.$

We collect some useful facts here, see \cite{FH} for proofs.

\begin{prop}[\cite{FH}]\label{bsic-aa}
 If $(X,T)$ is a TDS and $\Phi\in \C_{aa}(X,T,d)$,
then the mapping $\Phi_\ast: \M(X,T)\to L_\Phi$ is a continuous
surjection, where $ \M(X,T)$ is endowed with the weak star topology.
The set $L_\Phi$ is a compact set in $\R^d.$ If moreover $(X,T)$
satisfies the specification property (see for instance \cite{GP97} for the definition), then $L_\Phi$ is convex and
$E_\Phi(\alpha)\ne \emptyset$ if and only if $\alpha\in L_\Phi.$

\end{prop}

Like in Remark 3.7 in \cite{FH}, if we define an
equivalence relation on $\C_{aa}(X,T,d)$ by $\Phi\sim \Psi$ if there
exists $\alpha\in\R^d$ such that  $\|\Phi-\Psi-\alpha\|_{\rm
lim}=0,$ then it is not hard to see that the quotient space
$\C_{aa}(X,T,d)/\sim$ with the norm $\|\cdot\|_{\rm lim}$ is a
separated Banach space. The class of $\Phi$ is denoted $\bar \Phi$. We have the following relation between the
dimension of $L_\Phi$ and the dimension of the subspace $\langle
\bar\Phi^1,\cdots,\bar\Phi^d\rangle$ (see Section~\ref{proofs} for the proof):

\begin{prop}\label{dim-L-phi}
 Assume $(X,T)$ is a TDS with specification.  Let $\Phi=(\Phi^1,\cdots,\Phi^d)\in\C_{aa}(X,T,d)$. Then
$L_\Phi$ is of  dimension $d$ if and only if
$\bar\Phi^1,\cdots,\bar\Phi^d\in \C_{aa}(X,T)/\sim$ are linearly
independent.
\end{prop}

The thermodynamic formalism for almost additive potentials has been
studied in several works ~\cite{Fal88,Bar96,FeLa02,Fen04,B,Mum06,BD,CFH}. For our purpose, we only need to consider
 the subshift of finite type case. Let $(\Sigma_A,T)$ be a subshift of
finite type.
 Given
$\Phi\in\C_{aa}(\Sigma_A,T)$, the topological pressure can be
defined as
\begin{equation}\label{pressure-topo}
P(T,\Phi):=\lim_{n\to\infty}\frac{1}{n}\log \sum_{w\in \Sigma_{A,n}}
\exp(\sup_{x\in[w]}\phi_n(x)).
\end{equation}
Usually we write $P(\Phi)$ for $P(T,\Phi)$ when there is no
confusion. The following extension of the classical variational principle (see \cite{Ruelle}) holds:

\begin{thm}{\cite{B,BD,CFH}}\label{varia-principle}
Let $(\Sigma_A,T)$ be a subshift of finite type. For any
$\Phi\in\C_{aa}(\Sigma_A,T)$, we have
$
P(T,\Phi)=\sup\{h_\mu(T)+\Phi_\ast(\mu): \mu\in\M(\Sigma_A,T)\}.
$
\end{thm}

%%%%%%%%%%%%%%%%%%%%%%%%%%%%%%%%%%%%%%%%  weak Gibbs metric   %%%%%%%%%%%%%%%%%%%%%%%%%%%%%%%%%%%%%%%%%%%
\subsubsection{Weak Gibbs metric on subshift of finite type}
Let $(\Sigma_A,T) $ be a topologically mixing subshift of finite
type with alphabet $\{1,\cdots,m\}$, where $A$ is a $m\times m$
matrix with entries $0$ and $1$ such that $A^{p_0}>0$ for some
$p_0\in\N$ and $T$ is the shift map. We endow $\Sigma_A$ with a
metric naturally associated with a potential
 $\Psi\in\mathcal{C}_{aa}^-(\Sigma_A,T)$. This kind of metrics have been considered in
 \cite{GP97} and \cite{KS04} associated with additive potentials.

Note that by endowing $\Sigma_A$ with the  standard
  metric $d_1$  defined as $d_1(x,y)=m^{-|x\wedge y|}$
  (where $|x\wedge y|$ is the length of the common prefix of $x$ and $y$), $(\Sigma_A,
  d_1)$ is a compact metric space and $(\Sigma_A,T)$ is a TDS satisfying the specification property.
 Let $\Sigma_{A,n}$ be the set of the admissible words of length $n$
 and let
$\Sigma_{A,\ast}:=\bigcup_{n\geq 0} \Sigma_{A,n}.$ For
$w\in\Sigma_{A,\ast} $ and $w=w_1\cdots w_n$, we denote the length
of $w$ by $|w|=n$. Given $w\in \Sigma_{A,\ast}\cup \Sigma_{A}$ with
$|w|\geq n$, we denote $w_1\cdots w_n$ by $w|_n$.  Given
$u\in\Sigma_{A,\ast} $ and $v\in \Sigma_{A,\ast} \cup\Sigma_{A}$, if
$u_j=v_j$ for $j=1,\cdots,|u|$, then we say $u$ is a {\it prefix }
of $v$ and write $u\prec v.$ For $u=u_1\cdots u_n\in \Sigma_{A,n}$,
$u^\ast$ stands for $u|_{n-1}$. For $x, y\in \Sigma_{A,\ast}\cup
\Sigma_A$ such that $x\ne y$, $x\wedge y$ stands for the common
prefix of $x$ and $y$ of maximal length.

Recall that  $A^{p_0}(i,j)>0$ for all $1\leq i,j\leq m$,
consequently $A^{p_0+2}(i,j)>0$. For each $i,j$ we fix
 $w(i,j)\in \Sigma_{A,p_0}$ such that $iw(i,j)j$ is admissible.
Define $\Xi:=\{w(i,j): 1\leq i,j\leq m\}.$

Given a continuous function $\phi:\Sigma_A\to \R^d$, we define
\begin{equation}\label{normn}
\|\phi\|_n:=\sup_{x|_n=y|_n} |\phi(x)-\phi(y)|,
\end{equation}
and for $\Phi\in\C_{aa}(\Sigma_A,T,d)$ we write
$\|\Phi\|_n:=\|\phi_n\|_n$. Writing $\Phi=(\Phi^1,\cdots,\Phi^d),$  we
have 
\begin{equation}\label{component}
(\sum_{j=1}^{d}\|\Phi^j\|_n^2)^{1/2}\leq \sqrt{d}\|\Phi\|_n.
\end{equation}

For $\Phi\in \C_{aa}(\Sigma_A,T)$ and $w\in\Sigma_{A,n}$ we define
$$
\Phi[w]:=\sup\{\exp({\phi_n(x)}):x\in [w]\}.
$$

Now we fix a $\Psi\in  \C_{aa}^-(\Sigma_A,T).$ For $x,y\in \Sigma_A$
define
$$
d_\Psi(x,y):=
\begin{cases} \Psi[x\wedge y], & \text{ if }
x\ne y\\
0, & \text{ if } x=y.\end{cases}
$$
\begin{prop}\label{metric-sym}
 $d_\Psi$ is an ultra-metric  on $\Sigma_A.$ If $x\in \Sigma_A$ and $r>0$, the closed ball
$B(x,r)$ is the cylinder $[x|_n]$, where $n$ is the unique integer
such that $\Psi[x|_{n-1}]> r$ and $\Psi[x|_n]\leq r$. Each
cylinder $[w]$ is a ball with $\mathrm{diam}([w])=\Psi[w]$.
\end{prop}

The proof is elementary and we omit it.
%\proof We only
%need to show the triangular  inequality. Let $x\ne
% y\in \Sigma_A$. For any $z\in\Sigma_A$,  either $x\wedge z\prec x\wedge y$, or  $y\wedge z\prec x\wedge
% y$. Thus either $\Psi[x\wedge z]\geq \Psi[x\wedge y]$ or $\Psi[y\wedge z]\geq \Psi[x\wedge
%y]$. Then we have
% $$d_\Psi(x,y)=\Psi[x\wedge y]\leq \max\{\Psi[x\wedge z],\Psi[y\wedge
%z]\}=\max\{d_\Psi(x,z),d_\Psi(y,z)\}.$$
% Thus $d_\Psi$ is ultra   metric.
%Let $n$ be the unique number such that $\Psi[x|_{n-1}]> r$ and
%$\Psi[x|_n]\leq r$. If $y\in [x|_n]$, then $  x|_{n}\prec y$ and
%$d_\Psi(x,y)=\Psi[x\wedge y]\leq\Psi[x|_n]\leq r$. On the other hand
%if $y\not\in [x|_n]$, then $x\wedge y\prec x|_{n-1}$ and
%$d_\Psi(x,y)=\Psi[x\wedge y]\geq\Psi[x|_{n-1}]> r $. Thus we have
%$[x|_n]= B(x,r).$ \hfill$\Box$
%\begin{rem}{\rm  If we take $\psi_n(x)\equiv -n\log m$, then  $d_\Psi=d_1.$
%In section  \ref{examples} we will take  special $\Psi$ such that
%the metric $d_\Psi$ is related to the geometric structures of the
% conformal repeller.
%}\end{rem}
For the metric space $(\Sigma_A,d_\Psi)$ we define
$${\mathcal B}_n(\Psi)=\{w\in\Sigma_{A,\ast}:[w] \text{ is a closed ball of }  \Sigma_A \text{ with radius } e^{-n}\}\quad(n\ge 0).$$
We note that $\Sigma_A=\bigcup_{w\in{\mathcal B}_n(\Psi)}[w]$ for each $n\ge 0$.

\subsubsection{Three dimension functions} We introduce three functions which will turn
out to take the same values on $L_\Phi$ and provide the Hausdorff and packing dimensions of
the sets $E_\Phi(\alpha)$. They correspond to different point of views to estimate these dimensions,
 namely box-counting of balls intersecting $E_\Phi(\alpha)$,  variational principle for entropy like \eqref{varprinc}
   and  Legendre transform of a kind of {\it metric} pressure. The proofs of the propositions stated in this section are given in Section~\ref{proofs}.

\noindent {\bf (1) Box-counting type function, the large deviation
spectrum:} fix $\Psi\in \C_{aa}^-(\Sigma_A,T)$ and $\Phi\in
\C_{aa}(\Sigma_A,T,d)$. Define  $d_\Psi$ and ${\mathcal B}_n(\Psi)$
as above. Given  $\alpha\in L_\Phi$, $n\geq 1$ and $\epsilon>0$,
define
$$F(\alpha,n,\epsilon,\Phi,\Psi):=\Big \{u\in {\mathcal B}_n(\Psi):\ \text{ there exists}\  x\in[u] \ \text{such that
}\  \Big |\frac{\phi_{|u|}(x)}{|u|}-\alpha\Big |<\epsilon \Big \}.$$
Let $f(\alpha,n,\epsilon,\Phi,\Psi)$ be the cardinality of
$F(\alpha,n,\epsilon,\Phi,\Psi)$.

\begin{prop}\label{dim-formular-1} For any
$\Psi\in \C_{aa}^-(\Sigma_A,T)$, the limit
\begin{equation}\label{full-dim}
D(\Psi):=\lim_{n\to\infty}\frac{\log \#{\mathcal B}_n(\Psi)}{n}
\end{equation}
exists.  Moreover there exist constants $C_2(\Psi)>C_1(\Psi)>0$ such
that
\begin{equation}\label{ful-dim-1}
C_1(\Psi)\log m\leq D(\Psi)\leq C_2(\Psi)\log m.
\end{equation}
 For any $\Phi\in
\C_{aa}(\Sigma_A,T,d)$ and any $\alpha\in L_\Phi$, we have
$$\lim_{\epsilon\to 0}\liminf_{n\to\infty}\frac{\log f(\alpha,n,\epsilon,\Phi,\Psi)}{n}=
\lim_{\epsilon\to 0}\limsup_{n\to\infty}\frac{\log
f(\alpha,n,\epsilon,\Phi,\Psi)}{n}=:\Lambda_{\Phi}^{\Psi}(\alpha).$$
The function $\Lambda_{\Phi}^{\Psi}:L_\Phi\to \R$ is  upper
semi-continuous.
\end{prop}

We will prove that $\Lambda_\Phi^\Psi(\alpha)$ is the Hausdorff
dimension of $E_\Phi(\alpha)$ for all $\alpha\in L_\Phi$. The
function $\Lambda_\Phi^\Psi$ has more regularity than upper
semi-continuity. To state it we  need several standard notations
from convex analysis. Recall that  a subset $M$ of $ \R^d$ is an {\it affine
subspace} if $\lambda x+(1-\lambda)y\in M$ for every $x,y\in M$ and
$\lambda\in \R.$ Given $A\subset\R^d$, the {\it affine hull} of $A$
is the smallest affine subspace $M$ of $\R^d$ such that $M\supset A$ and is denoted
by $\mathrm{aff}(A).$ For a convex set $A$, we define
$\mathrm{ri}(A),$ the {\it relative interior} of $A$ as
$\mathrm{ri}(A):=\{x\in \mathrm{aff}(A): \ \exists \epsilon>0,
(x+\epsilon B)\cap \mathrm{aff}(A)\subset A\},$ where
$B=B(0,1)\subset \R^d$ is the unit ball.  Let $A\subset \R^d$ be a
convex set and $h: A\to \R$ be a function. If
 there exists $c\geq 1$ such that for any $\alpha,\beta\in A$, we can
find $\gamma_1=\gamma_1(\alpha,\beta),
\gamma_2=\gamma_2(\alpha,{\beta})\in [c^{-1},c]$ such that for any
$\lambda\in[0,1]$
\begin{equation}\label{lower-semi-conti}
\lambda h(\alpha)+(1-\lambda)h(\beta)\leq
h\Big (\frac{\lambda\gamma_1\alpha+(1-\lambda)\gamma_2\beta}
{\lambda\gamma_1+(1-\lambda)\gamma_2}\Big ),
\end{equation}
then we call $h$  a {\it weakly concave function} on $A.$ Note
that if $c=1$, we go back to the usual concept of concave function. Also,
$h(\gamma)\ge \min (h(\alpha),h(\beta))$ if $\gamma\in [\alpha,\beta]\subset A$, thus $h$ is quasi-concave.

 \begin{prop}\label{regularity}
  The function $\Lambda_\Phi^\Psi: L_\Phi\to \R$ is
 bounded, positive and weakly concave. It is  continuous on any closed interval $I\subset
 L_\Phi$ and  on $\mathrm{ri}(A)$, where $A\subset L_\Phi$ is any convex set.
  Consequently it is continuous on ${\rm ri}(L_\Phi).$ If moreover $L_\Phi$ is a convex polyhedron,
 then $\Lambda_\Phi^\Psi$ is continuous on $L_\Phi.$
  Assume $I=[\alpha_0,\alpha_1]\subset L_\Phi$  and $\alpha_{\max}\in
 I$ such that $\Lambda_\Phi^\Psi(\alpha_{\max})=\max\{\Lambda_\Phi^\Psi(\alpha): \alpha\in
 I\}$, then $\Lambda_\Phi^\Psi$ is decreasing from $\alpha_{\max}$ to
 $\alpha_j$, $j=0,1.$
 \end{prop}

 \begin{rem}{\rm  Large deviations spectra for the Hausdorff dimension estimation of sets like $ E_\Phi(\alpha)$
 have been considered since the first studies of multifractal properties of Gibbs or weak Gibbs measures
  and then extended to the study of Birkhoff averages \cite{Collet,Rand,BMP,PeW,PW,P,BS,O,FF,FFW,BSS,FLW}.
  Until now, in the situations where such a spectrum may be non-concave \cite{BSS,BD,FLW},
  no description of its regularity like that of Proposition~\ref{regularity} had been given.
   Moreover, the methods used in the papers mentioned above seem not adapted to provide this information.
% Specifically, under the assumptions of \cite{BSS}, it is not necessary to consider such a spectrum to estimate Hausdorff dimensions, because the $E_\Phi(\alpha)$ is studied for those $\alpha$ of the relative interior of $L_\Phi$, and there one can show that such a set carries a Gibbs measure; in \cite{FLW}, the large deviation spectrum is defined in a different way, by using cylinders of the same generation in sets like $F(\alpha,n,\epsilon,\Phi,\Psi)$, and .
 }\end{rem}

\noindent {\bf (2) Function associated with a conditional variational principle:}  For $\alpha\in L_\Phi$ let
$$
{\mathcal E}_\Phi^\Psi(\alpha):=\sup\left
\{\frac{h_\mu(T)}{-\Psi_\ast(\mu)}:\mu\in\M(\Sigma_A,T)\ \
\text{such that } \ \ \Phi_\ast(\mu)=\alpha\right\}.
$$
\noindent {\bf (3) Pressure type function and its Legendre transform-like associated function:} at first we define a kind of
pressure function.

\begin{prop}\label{gene-thermo}
 Fix $\Phi\in
\C_{aa}(\Sigma_A,T,d)$ and $\Psi\in \C_{aa}^-(\Sigma_A,T)$. Let
$z,\alpha\in\R^d.$ Then the equation
\begin{equation}\label{solution}
P(\langle z,\Phi-\alpha\rangle+\tau_{\Phi}^{\Psi}(z,\alpha)\Psi)=0
\end{equation}
 has a unique solution
$\tau_{\Phi}^\Psi(z,\alpha).$ Moreover  the following variational principle holds
\begin{equation}\label{geometricpressure}
\tau_{\Phi}^\Psi(z,\alpha)=\sup\left \{\frac{h_\mu(T)+ \langle
z,\Phi_\ast(\mu)-\alpha\rangle }{-\Psi_\ast(\mu)}
 : \mu\in {\mathcal M}(\Sigma_A,T)\right \},
 \end{equation}
and one also has %When $\Phi$ and $\Psi$ are H\"{o}lder potentials, we have
\begin{equation}\label{pressure-metric}
\tau_{\Phi}^\Psi(z,\alpha)=\lim_{n\to\infty}\frac{1}{n}\log
\sum_{w\in {\mathcal B}_n(\Psi)} \exp(\sup_{x\in[w]}\ \langle
z,(\phi_{|w|}(x)-|w|\alpha)\rangle).
\end{equation}
\end{prop}

We will identify the function
\begin{equation}\label{gene-Legendre-trans}
\tau_{\Phi}^{\Psi\star}(\alpha):=\inf\{\tau_{\Phi}^\Psi(z,\alpha):z\in\R^d\}
\end{equation}
with $\Lambda_\Phi^\Psi$ on $L_\Phi$. This is the large deviations principle announced in the introduction.

\begin{rem}\label{remlambda}{\rm  (1) The function $\tau_{\Phi}^\Psi(z,\alpha)$ defined as in $(\ref{solution})$
 was first introduced in~\cite{BSS} when $\Phi$ and $\Psi$ are H\"older potentials,
  and also in ~\cite{BD} where H\"older condition on $\Phi$ is replaced by the bounded distorsion property for almost additive potentials.

%(2) By comparing $(\ref{pressure-topo}) $ with \eqref{geometricpressure} and
%$(\ref{pressure-metric})$, we see that
%$\tau_{\Phi}^\Psi(z,\alpha)$ is a kind of pressure function, which
%is related to the metric $d_\Psi$, or to the geometry for the
%geometric realizations of this model.

\noindent (2)  $\tau_{\Phi}^{\Psi\star}$ is a  generalization of Legendre
transform as noted in~\cite{BSS}. For the special constant potential
 $\Psi=(-n)_{n=1}^{\infty}$, $\tau_{\Phi}^\Psi(z,\alpha)=P(\langle z,\Phi\rangle)-\langle
z,\alpha\rangle$, thus $\tau_{\Phi}^{\Psi\star}$ is the classical
Legendre transform of the pressure function $P(\langle
z,\Phi\rangle)$. }\end{rem}

\subsection{Main results on topologically mixing subshift of finite type}

Throughout this subsection we fix $\Phi\in
\C_{aa}(\Sigma_A,T,d)$ and $\Psi\in \C_{aa}^-(\Sigma_A,T)$. We work
on the metric space $(\Sigma_A,d_\Psi)$. We write $\dim_H^\Psi E,
\dim_P^\Psi E, \dim_B^\Psi E $ for the Hausdorff, packing and box
dimensions of $E\subset \Sigma_A$. For convenience we write
$\D_\Phi^\Psi(\alpha):=\dim_H^\Psi E_\Phi(\alpha)$.

%\newpage

 \begin{thm}[{\bf Multifractal analysis of the level sets $E_\Phi(\alpha)$}]\label{main-one-sided}$\ $

\begin{enumerate}
\item $E_\Phi(\alpha)\ne\emptyset$ if and only if $\alpha\in L_\Phi$.
For  $\alpha\in L_\Phi$ we have
$$
\D_\Phi^\Psi(\alpha)=
\Lambda_\Phi^\Psi(\alpha)={\mathcal E}_\Phi^\Psi(\alpha)=\tau_{\Phi}^{\Psi\star}(\alpha),$$
and the function $\D_\Phi^\Psi$ is weakly concave.

\item $\dim_H^\Psi\Sigma_A=\dim_B^\Psi\Sigma_A=D(\Psi)=\max\{\Lambda_\Phi^{\Psi}(\alpha): \alpha\in L_\Phi\}$.
\end{enumerate}
\end{thm}

%%%%%%%%%%%%%%%%%%%%%%%%%%%%%%%%%%%%%%%%%%%%%%%%%%%%%%%%%%%%%%%%%%%%%%%%%%%%%%%%%%%%%%%%%%%%%%%%%%%%%%%

%Now we come to our second main result.

%%%%%%%%%%%%%%%%%%%%%%%%%%%%%%%%%%%%%%%%%%%%%%%%%%%%%%%%%%%%%%%%%%%%%%%%%%%%%%%%%%%%%%%%%%%%%%%%%%%%%%%%%%%
 \begin{thm}[{\bf Localized asymptotic behavior}]\label{main-fun-level-one-sided}
 Assume $\xi:\Sigma_A\to\R^d$ is continuous
 and
$\xi(\Sigma_A)\subset \mathrm{aff}(L_\Phi)$. 

\begin{enumerate}
\item $\dim_H^\Psi E_\Phi(\xi)\geq
\sup\{\mathcal{D}_\Phi^\Psi(\alpha):\alpha\in \xi(\Sigma_A)\cap
\mathrm{ri}(L_\Phi)\}.$ 

\item If $\xi(\Sigma_A)\subset L_\Phi$
then $E_\Phi(\xi)$ is dense in $\Sigma_A$.

\item If
$\sup\{\mathcal{D}_\Phi^\Psi(\alpha):\alpha\in\xi(\Sigma_A)\cap
\mathrm{ri}(L_\Phi)\}= \sup\{\mathcal{D}_\Phi^\Psi(\alpha):\alpha\in
\xi(\Sigma_A)\cap L_\Phi\}$, then ${\dim}_H^\Psi
E_\Phi(\xi)=\dim_P^\Psi
E_\Phi(\xi)=\sup\{\mathcal{D}_\Phi^\Psi(\alpha):\alpha\in
\xi(\Sigma_A)\cap L_\Phi\}$.

\item If $d=1$ and $\xi(\Sigma_A)\subset L_\Phi$, then $E_\Phi(\xi)$
is dense and $ \dim_H^\Psi E_\Phi(\xi)=\dim_P^\Psi
E_\Phi(\xi)=\sup\{\mathcal{D}_\Phi^\Psi(\alpha):\alpha\in
\xi(\Sigma_A)\}. $
\end{enumerate}\end{thm}

\begin{rem}\label{connexion}{\rm   (1) Recall Remark~\ref{remlambda}(2). In \cite{FF,FFW,Fen03},  where the metric is the standard one, the  equality $\Lambda_\Phi^\Psi(\alpha)=\inf_{z\in\R^d}P(\langle z,\Phi\rangle)-\langle
z,\alpha\rangle$ is established, and both functions are concave. In our work, the weak concavity of $\Lambda_\Phi^\Psi$ turns out to be crucial in proving the equality $\Lambda_\Phi^\Psi(\alpha)=\tau_{\Phi}^{\Psi\star}(\alpha)$ in full generality. This equality can be read as a new large deviation principle thanks to the expression \eqref{pressure-metric} giving $\tau_{\Phi}^{\Psi}$ as a kind of metric pressure. 
 
\noindent (2) In ~\cite{BSS,BD}, assuming more regularity than
continuity for $\Phi$ and $\Psi$, namely bounded distorsion property
or H\"older continuity in the case where the potential is additive,
the equality $\D_\Phi^\Psi(\alpha)= {\mathcal
E}_\Phi^\Psi(\alpha)=\tau_{\Phi}^{\Psi\star}(\alpha)$ is shown only
for $\alpha\in \mathrm{int}(L_\Phi)$, where ${\rm int}(A)$ denotes
the interior of $A\subset\R^d$. The argument is strongly based on
the differentiability of the pressure function in these cases. 

\noindent
(3) In \cite{FLW}, the authors consider the case of additive potentials $\Phi$ and $\Psi$,
and work under the assumption that $\Psi$ corresponds to a H\"older potential. They show $\D_\Phi^\Psi(\alpha)=
{\mathcal E}_\Phi^\Psi(\alpha)$ for all $\alpha\in L_\Phi$. Here we work under weaker regularity assumptions on
 $\Psi$, and both $\Phi$ and $\Psi$ are almost additive. Also, we use a different method
 to compute the function $\D_\Phi^\Psi(\alpha)$, namely concatenation of Gibbs measures.
 Such a method has been used successfully in \cite{KS04} to deal with the special sets
 $E_\Psi^\Psi(\alpha)$ when $\Psi$ is additive as well as in \cite{BaFe09} to deal with
 the asymptotic behavior of almost additive potentials in the different context of full-shifts
 endowed with self-affine metrics (the spectrum is always concave in this case). Here, we need to
 refine such approach in order to remove some delicate points in our geometric
  application to attractors of $C^1$ conformal iterated function systems.
% We notice that an alternative approach to compute $\D_\Phi^\Psi(\alpha)$ consists in using the large deviation spectrum whose computation is based on
%   box-counting of balls of the same diameter than cylinders of the same generation.
%    This simplifies a lot the Hausdorff dimensions calculations, and this is also more closely related to
%     a kind of pressure like the one introduced in  \eqref{geometricpressure} and
%$(\ref{pressure-metric})$. This is also crucial to get the weak
%concavity of $\Lambda_\Phi^\Psi$. 
We also mention that
in the case where the metric is the standard one on a full-shift, the equalities $\D_\Phi^\Psi(\alpha)=
{\mathcal E}_\Phi^\Psi(\alpha)=\tau_{\Phi}^{\Psi\star}(\alpha)$ (there the spectrum is concave)
have been obtained in \cite{Fen03} when $\Phi$ is built from 
 Birkhoff products of continuous positive matrices. There, the computation of Hausdorff
 dimension  uses concatenation of words, like in \cite{FF,FFW,FLW}.
 }\end{rem}

\begin{rem}\label{remloc}
{\rm \noindent (1) The proof Theorem~\ref{main-fun-level-one-sided} uses the weak concavity of the spectrum
$\D_\Phi^\Psi$. It also requires to
concatenate
 Gibbs measures in a more elaborated way than to determine~$\D_\Phi^\Psi$.

 \noindent(2) In fact we shall prove a slightly more general result than Theorem~\ref{main-fun-level-one-sided}(1): (1') Suppose that $\xi$ is bounded and continuous
outside a subset $E$ of $\Sigma_A$, and $\xi(\Sigma_A)\subset \mathrm{aff}(L_\Phi)$.
 If $\dim_H^\Psi E<\sup\{\D_\Phi^\Psi(\alpha):\alpha\in \xi (\Sigma_A\setminus E)\cap \mathrm{ri}(L_\Phi)\}$,
 then $\dim_H^\Psi E_\Phi(\xi)\geq \sup\{\D_\Phi^\Psi(\alpha):\alpha\in \xi (\Sigma_A\setminus E)\cap \mathrm{ri}(L_\Phi)\}.$

 \noindent(3) An extension of Theorem~\ref{main-fun-level-one-sided}(4)  is given in the final remark of Section~\ref{appli2}.}
\end{rem}

\section{Geometric results}\label{examples}

In this section we show how the main results of the previous section can be applied to
 multifractal analysis on conformal repellers and on attractors of conformal IFS
 satisfying the strong open set condition. Such sets  fall in
 the Moran-like geometric
 constructions considered in \cite{Bar96,P}. At first
 we describe this kind of construction (Section~\ref{gengeom}). Then we state
  the geometric results deduced from Theorems~\ref{main-one-sided} and \ref{main-fun-level-one-sided}
  (Section~\ref{maingengeom}). We give our application to {\it fixed points
    in the asymptotic average for dynamical systems} in $\R^d$ in Section~\ref{asympave}. Finally, we give an application
to the local scaling properties of weak Gibbs measures in
Section~\ref{appli2}, special example of which is the harmonic measure
on planar conformal Cantor sets.

\subsection{General setting of geometric realization}\label{gengeom}

  Let $(\Sigma_A,T) $ be a topologically mixing subshift of finite
type with alphabet $\{1,\cdots,m\}$ and $\Psi\in
\C_{aa}^-(\Sigma_A,T)$. Let $X$ be $\R^{d^\prime}$ or be a
connected, $d^\prime$-dimensional $C^1$ Riemannian manifold.
Consider a family of sets $\{R_w: w\in \Sigma_{A,\ast}\}$, where
each $R_w\subset X$ is a compact set with nonempty interior. We
assume that this family of compact sets satisfies the following
conditions:

(1) $R_w\subset R_{w^\prime}$ whenever $w^\prime\prec w$.

(2)\ For any integer $n>0$, the interiors of distinct $R_w, w\in
  \Sigma_{A,n}$ are disjoint.

 (3)\ Each $R_w$ contains a ball of radius $\underline{r}_w$ and is
  contained in a ball of radius $\overline{r}_w$.

  (4)\ There exists a  constant $K>1$ and a negative sequence $\eta_n=o(n)$ such that
  for every $w\in\Sigma_{A,\ast}$,
  \begin{equation}\label{conformal-basic}
  K^{-1}\exp(\eta_{|w|})\Psi[w]\leq\underline{r}_w\leq \overline{r}_w\leq K\Psi[w].
   \end{equation}

Let $\displaystyle J=\bigcap_{n\geq 0}\bigcup_{w\in
\Sigma_{A,n}}R_w.$ We call $J$ the {\it limit set of the family
$\{R_w:w\in\Sigma_{A,\ast} \}$.} We can define the coding map $\chi:
\Sigma_A\to J$ as
$\displaystyle \chi(x)=\bigcap_{n\geq 1}R_{x|_n}, \ \ \forall x\in
\Sigma_A.$
It is clear that $\chi$ is continuous and surjective.

We say that $J$ is a {\it Moran type geometric realization of
$\Sigma_A$ with potential $\Psi$.}

For this kind of construction we have the following useful
observation:
%(see Lemma~\ref{dim-compare-lem}):

\begin{prop}\label{dim-compare-prop}
 Let $J$ be a Moran type geometric realization of
$\Sigma_A$ with almost additive potential $\Psi$, then for any
$E\subset J$ we have $\dim_H E=\dim_H^\Psi(\chi^{-1}(E))$.
\end{prop}

In this paper we consider two classes of Moran type geometric
realizations of $\Sigma_A.$

{\bf (1)}\ Topologically mixing $C^1$ conformal repeller $(J,g)$. We refer
the book \cite{P} for the definitions and the basic properties
related to conformal repellers. It is well known that in this case
$(J,g)$ has a Markov partition $\{R_1,\cdots,R_m\}$. For each
$w=w_1\cdots w_n,$ define $R_w:=R_{w_1}\cap
g^{-1}(R_{w_2})\cap\cdots\cap g^{-n+1}(R_{w_n})$. Define
$\psi(x)=-\log |g^\prime(\chi(x))|$ and
$\Psi=(S_n\psi)_{n=1}^{\infty}$. By the definition of $R_w$ and the
property of Markov partition, the condition (1) and (2) are checked
directly.  (3) and (4) are stated in \cite{P} (Proposition 20.2),
except that for (4) we have an additional term
$\exp(\eta_{|w|})=\exp(-\|\Psi\|_{|w|})$. This is because we only
assume $\psi$ to be continuous rather than H\"older continuous.
Thus $J$ is a Moran type geometric realization of $\Sigma_A$ for
some primitive matrix $A$ and the potential $\Psi.$ Moreover in this
case we have $\chi\circ T=g\circ \chi.$

{\bf (2)} Attractors of $C^1$ conformal IFS satisfying the strong open set
condition. For completeness we recall the related definitions. Let
$U\subset\R^d$ be a non-empty open set. A map $f: U\to U$ is {\it
contracting} if there exists $0<\gamma<1$ such that $|f(x)-f(y)|\leq
\gamma|x-y|$ for all $x,y\in U.$ Let $\{f_1,\cdots,f_m\}$ be a
collection of contracting maps from $U$ to $U$ and suppose that for
some closed set $X\subset U$ we have $f_j(X)\subset X$ for each $j$.
Then, it is well known that there is a unique non-empty compact set
$J\subset X$ such that $J=\bigcup_{j=1}^{m}f_j(J)$. Such a family is
called an {\it Iterated Function System (IFS)}, of which $J$ is the
{\it attractor }. This IFS is said to satisfy the {\it open set
condition (OSC)} if there is a non-empty open set $V\subset U$ such
that $f_j(V)\subset V$ for each $j$ and $f_i(V)\cap
f_j(V)=\emptyset$ for $i\ne j.$ The {\it strong open set condition
(SOSC)} holds if moreover this open set $V$ can be chosen with
$V\cap K\ne \emptyset.$

 A $C^1$-map $f:U\to\R^d$ is {\it conformal}
if the differential $f^\prime(x):\R^d\to\R^d$ satisfies
$|f^\prime(x)y|=|f^\prime(x)||y|\ne 0$ for all $x\in U$ and
$y\in\R^d$,  $y\ne 0$. We say that an IFS $\{f_1,\cdots,f_m\}$ is a
{\it $C^1$ conformal IFS} if each $f_j$ is an injective conformal
map.  We refer to \cite{PRSS} for more details.

 Assume $\{f_1,\cdots,
f_m\}$ is a $C^1$ conformal IFS satisfying the SOSC. Let $J$ be its
attractor.  Define $ \psi(x)=\log |f_{x_1}^\prime(\chi(Tx))|\ \
\text{and}\ \ \Psi=(S_n\psi)_{n=1}^{\infty}. $ Let $V$ be an open
set such that the SOSC holds. For $w=w_1\cdots w_n,$ define
$R_w=f_{w}(\overline{V})$, where $f_{w}:=f_{w_1}\circ \cdots\circ
f_{w_n}$. Due to the SOSC,  (1) and (2) hold for $\{R_w: w\in
\Sigma_{A,\ast}\}$. Moreover, arguments similar to those used to
prove Proposition 20.2 in \cite{P}  show that (3) and (4) also hold.
Thus, $\{R_w: w\in \Sigma_{A,\ast}\}$ is a Moran type geometric
realization of $\Sigma_A$ with potential $\Psi$. Notice that here
$\Sigma_A$ is the full shift $\Sigma_m.$ By the uniqueness of the
attractor it is easy to verify that the attractor $J$ is the limit
set of the family $\{R_w: w\in \Sigma_{A,\ast}\}$.

\subsection{Multifractal analysis on Moran type geometric realizations}\label{maingengeom}
We are going to conduct multifractal analysis on Moran type
geometric realizations, thus we need a dynamics $g$ on $J$ so that
$(J,g)$ be a factor of some $(\Sigma_A,T)$. For $C^1$ conformal
repellers, there is such a natural dynamic. For the attractor of a
$C^1$ conformal IFS, there is no such one in general, the difficulty
coming from those points having several codings. However, under the
SOSC, we can naturally define such a $g$ by removing a "negligible"
part of $J$:

Let $\{f_1,\cdots,f_m\}$ be a $C^1$ conformal IFS satisfying the
SOSC. Let $V$ be an open set such that the SOSC holds. By \cite{PRSS}, such an open set always exists as soon as the mappings $f_i$ are $C^{1+\epsilon}$ and the OSC holds. Define
$\widetilde Z_\infty:=\bigcup_{w\in \Sigma_{A,\ast}} f_w(
\partial V)$ and $
Z_\infty:=\chi^{-1}(\widetilde Z_\infty)$. We have the following
lemma (proved in Section~\ref{sec7}):

\begin{lem}\label{boundary}
 The set $\Sigma_{A}\setminus  Z_\infty$ is not empty and
 $\chi: \Sigma_{A}\setminus  Z_\infty\to J\setminus
 \widetilde Z_\infty$ is a bijection. Moreover $T(\Sigma_A\setminus Z_\infty)\subset \Sigma_A\setminus
 Z_\infty$, $T(Z_\infty)\subset Z_\infty$ and for any Gibbs measure $\mu$ on $\Sigma_A$ we have $\mu(Z_\infty)=0.$
\end{lem}

By the previous lemma we can define the mapping $\widetilde g: J\setminus
\widetilde Z_\infty\to J\setminus \widetilde Z_\infty$  as
 $\widetilde g(x)=\chi \circ T\circ \chi^{-1}$.
  By construction we have $\chi\circ T= \widetilde g\circ\chi$ over $\Sigma_{A}\setminus  Z_\infty$.

Let $J$ be a Moran type geometric realization of $(\Sigma_A,T)$. We
set $\widetilde J=J$ when $J$ is a $C^1$ conformal repeller and
$\widetilde J=J\setminus \widetilde Z_\infty$ when $J$ is the
attractor of a $C^1$ conformal IFS satisfying the SOSC.

Given a sequence of functions  $\Phi=(\phi_n)_{n=1}^\infty$ from
$\widetilde J$ to $\mathbb{R}^d$ and $\alpha \in \mathbb{R}^d$, we
set $\displaystyle E_\Phi(\alpha)=\Big \{x\in \widetilde J:
\lim_{n\to\infty}{\phi_n(x)}/{n}=\alpha\Big \}. $ We also use the
notation $\D_\Phi(\alpha)=\dim_H E_\Phi(\alpha)$ and we define
$
L_\Phi=\{\alpha\in\R^d: E_\Phi(\alpha)\neq\emptyset\}.
$

When $J$ is a conformal repeller the system $(J,g)$ is naturally a TDS. For
$\Phi\in \C_{aa}(J,g,d)$, if we define $\widetilde
\Phi:=(\phi_n\circ\chi)_{n=1}^{\infty}$, since $g\circ
\chi=\chi\circ T$, we have $\widetilde \Phi\in
\C_{aa}(\Sigma_A,T,d)$ with $C(\widetilde \Phi)=C(\Phi)$. And for
$\alpha\in\R^d$ we have $E_\Phi(\alpha)=\chi(E_{\widetilde
\Phi}(\alpha)).$

When $J$ is the attractor of a $C^1$ conformal IFS satisfying the
SOSC,  if $\varphi$ is a continuous function from $J$ to
$\mathbb{R}^d$, it generates the additive potential
$\widetilde\Phi=(S_n\phi)_{n=1}^\infty$ on $(\Sigma_A,T)$, where
$\phi=\varphi\circ\chi$, and  it also defines
 $\Phi=(S_n\varphi)_{n=1}^\infty$ on $(\widetilde J,\widetilde g)$.  Then for
 $\alpha\in\R^d$ we have
$ E_\Phi(\alpha)=\chi (E_{\widetilde \Phi}(\alpha)\setminus
Z_\infty). $

\begin{thm}\label{appli-one-sided}
Let $J$ be a Moran type geometric realization of $(\Sigma_A,T)$. If
$J$ is a $C^1$ conformal repeller, let $\Phi\in \C_{aa}(J,g,d)$ and
define $\widetilde{\Phi}$ as above. If $J$ is the attractor of a
$C^1$ conformal IFS satistying the SOSC, let $\varphi$ be a
continuous map from $J$ to $\R^d$, and define the additive potential
$\widetilde\Phi=(S_n\phi)_{n=1}^\infty$ on $(\Sigma_A,T)$ with
$\phi=\varphi\circ\chi$ and $\Phi=(S_n\varphi)_{n=1}^\infty$ on
$(\widetilde J,\widetilde g)$. Then
\begin{enumerate}
\item $L_\Phi=L_{\widetilde\Phi}$;  for  $\alpha\in L_{\Phi}$ we have  $\dim_H
E_\Phi(\alpha)=\dim_PE_\Phi(\alpha)$ and
$$
\D_\Phi(\alpha)=\D_{\widetilde{\Phi}}^\Psi(\alpha)=
\Lambda_{\widetilde{\Phi}}^\Psi(\alpha)=\mathcal{E}_{\widetilde{\Phi}}^\Psi(\alpha)=\tau_{\widetilde{\Phi}}^{\Psi\star}(\alpha).$$

\item  $\dim_H J=\dim_B J=D(\Psi)=\max\{\D_\Phi(\alpha): \alpha\in L_{{\Phi}}\}$.
\end{enumerate}
\end{thm}

\begin{rem}
For the case of conformal repellers, the connection between
Theorem~\ref{appli-one-sided} and the other works \cite{BSS,FLW,BD}
is similar to that done in Remark~\ref{connexion}(2) and (3).
\end{rem}

For the set $E_\Phi(\xi)$ we have the following
result:

\begin{thm}\label{appli-fun-level-one-sided}
Let $J$ be a Moran type geometric realization of $(\Sigma_A,T)$,
which is either a $C^1$ conformal repeller or the attractor of a
$C^1$ conformal IFS satisfying the SOSC. Let $\Phi$ and $\widetilde
\Phi$ be the same as in Theorem \ref{appli-one-sided}. Let
$\xi:J\to\R^d$ be continuous and $\displaystyle
E_\Phi(\xi)=\Big\{x\in\widetilde J:
\lim_{n\to\infty}{\phi_n(x)}/{n}=\xi(x)\Big \}. $ If $\xi(J)\subset
\mathrm{aff}(L_{\Phi})$, then

\begin{enumerate}
\item  $\dim_H E_\Phi(\xi)\geq \sup\{\D_\Phi(\alpha):\alpha\in \xi(
J)\cap \mathrm{ri}(L_{{\Phi}})\},$ and $E_\Phi(\xi)$ is dense if
$\xi(J)\subset L_\Phi$.

\item  If  $\sup\{\D_\Phi(\alpha):\alpha\in \xi(J)\cap
\mathrm{ri}(L_{{\Phi}})\}= \sup\{\D_\Phi(\alpha):\alpha\in
\xi(J)\cap L_{{\Phi}}\}$, then ${\dim}_H E_\Phi(\xi)=\dim_P
E_\Phi(\xi)=\sup\{\D_\Phi(\alpha):\alpha\in \xi(J)\cap
L_{{\Phi}}\}$.

\item If $d=1$ and $\xi(J)\subset L_\Phi$, then  ${\dim}_H
E_\Phi(\xi)=\dim_P E_\Phi(\xi)=\sup\{\D_\Phi(\alpha):\alpha\in
\xi(J)\}$ and $E_\Phi(\xi)$ is dense.
\end{enumerate}
\end{thm}

 %%%%%%%%%%%%%%%%%%%%%%%%%%%%%%%%%%%%%%%%%%%%%%%%%%%%%%%%%%%%%%%%%%%%%%%%%%%%%%%%%%%%%%%%%%%%%%%%%%%%%%%%%%%%

\subsection{Application to fixed points  in the asymptotic average for dynamical systems in $\R^d$} \label{asympave}

Suppose that $(J,g)$ is a dynamical system with $J\subset \R^d$.  We say that $x\in J$ is a {\it fixed point of $g$
in the asymptotic average} if $\displaystyle \lim_{n\to\infty}\frac{1}{n}\sum_{k=0}^{n-1} g^k x=x$.
 We are interested in the Hausdorff dimension of the set of all  such  points:
 $$
\mathcal{F}(J,g)=\Big\{x\in J: \lim_{n\to\infty}\frac{1}{n}\sum_{k=0}^{n-1} g^k x=x\Big \}.
$$
If $\xi$ stands for the identity map on $J$ and $\Phi$ stands for the additive potential
 associated with the potential $\xi$, in our setting we have $\mathcal{F}(J,g)=E_\Phi(\xi)$.

The set $L_\Phi$ is contained in the convex hull of $J$, and it
contains the set  of the fixed points of $g$.  An example of trivial
situation is provided by the unit circle endowed with dynamic
$g(z)=z^2$ in $\mathbb{C}$. There, $\mathcal{F}(J,g)=\{1\}$. How
about general conformal repellers and attractors of conformal IFS?
This question is non trivial in general.  We are  going to describe
a class of conformal IFS, namely self-similar generalized Sierpinski
carpets, for which the situation is non trivial and we have a
complete answer.

% an interesting situations. We start
%with a result about the Hausdorff dimension of
%$\mathcal{F}(\widetilde J,\widetilde g)$
% associated with self-similar generalized Sierpinski carpets.
%Then we apply this result to two particular cases.

%\subsubsection{Fixed points in the asymptotic average for self-similar generalized Sierpinski carpets}
We consider a special self-similar IFS $\{f_1,\cdots,f_m\}$ on
$\R^d$:
 $f_j(x)=\rho_jx+c_j,  \ 0<\rho_j<1,  \ (1\leq j\leq m)$. We assume further the SOSC fulfills. Let $x_j$ stand for the unique fixed point of
 $f_j$ and let $J$ be the attractor  of this IFS.
 Notice that the
 mappings $f_j$ have no rotation part, thus the convex hull of $J$ satisfies ${\rm Co}(J)={\rm Co}\{x_1,\cdots,x_m\}=:\Delta$, and is a polyhedron. We further
 assume that ${\rm Co}(J)$ has dimension $d$ (otherwise we can define this IFS in
 a smaller
 affine subspace).

 Let $W$ stand for the open set such that the SOSC
 holds. It is ready to see that $V:=W\cap \Delta$ is also an open
 set such that SOSC holds. We can define the
 dynamics $\widetilde{g}$ on $\widetilde{J}=J\setminus \widetilde Z_\infty$, where $\widetilde Z_\infty$ is defined as in the previous subsection.

Now we have the following result whose proof is given in Section~\ref{sec7}.

 \begin{thm}\label{carpet}
  Let $\Phi={\rm id}_J.$ Then ${\mathcal F}(\widetilde{J},\widetilde{g})$ is dense and
  $\dim_{H}{\mathcal F}(\widetilde{J},\widetilde{g})=\sup\{\D_\Phi(\alpha): \alpha\in J\}. $
  Moreover if the point at which $\D_\Phi$ attains its maximum belongs to $J$,
  then ${\mathcal F}(\widetilde{J},\widetilde{g})$ is of full Hausdorff dimension.
 \end{thm}
We have the following corollary, in which the lower bound for the Hausdorff dimension follows directly from Theorem~\ref{carpet} and the upper bound follows from standard estimates based on the bounds provided in Section~\ref{upperbd}.  
\begin{coro}
Let $N\in \mathbb N_+$ and  let $d_1,\dots d_N$ be $N$ positive integers. Consider $N$
 self-similar IFS without rotations components $\{f_1^{(j)},\cdots,f_{m_j}^{(j)}\}_{1\le j\le N}$
 living respectively in $\R^{d_j}$. Denote by $J_j$, $1\le j\le N$, their respective attractors as
  well as the corresponding dynamical systems $(\widetilde J_j,\widetilde g_j)$.
   Let $\widetilde J=\prod_{j=1}^N \widetilde J_j\subset \R^{\sum_{j=1}^Nd_j}$
   be endowed with the dynamics $\widetilde g=(\widetilde g_1,\dots,\widetilde g_N)$.
   We have $\dim_H\mathcal{F} (\widetilde J,\widetilde g)=\sum_{j=1}^N \dim_H\mathcal{F}
   (\widetilde J_j,\widetilde g_j)=\sum_{j=1}^N\sup\{\D_{\Phi_j}(\alpha): \alpha\in
   J_j\}$, where $\Phi_j={\rm Id}_{\R^{d_j}}$.
\end{coro}

Both the previous results yield the result presented in the introduction of the paper:
\begin{thm}
Let $d\in\N_+$ and $(m_1,\dots,m_d)$ be $d$ integers $\ge 2$. Set
$J=[0,1]^d$ and let $g:J\to J$ be the mapping
$(x_1,\dots,x_d)\mapsto (m_1x_1\pmod 1,\dots,m_d x_d\pmod 1)$. Then
$\mathcal{F}(J,g)$ is dense and of full Hausdorff dimension in $[0,1]^d$.
\end{thm}
To see this, for a fixed integer $m\ge 2$ let $g_m:[0,1]\to[0,1]$ be the mapping $x\mapsto m x\pmod 1.$
%\begin{rem} Each point in $E$ is called an asymptotic fixed point of
%$\widetilde{\sigma}.$ At the first glance this result is a bit
%surprising, since
% by the ergodic theorem, for Lebesgue a.e.
%$x\in [0,1]^d$, we have
%$$
%\frac{\sum_{j=0}^{n-1}\widetilde{\sigma}^jx}{n}\to \int_{[0,1]^d} x\
%dx=(\frac{1}{2},\cdots,\frac{1}{2}).
%$$
%This theorem shows that although the set of asymptotic points has
%Lebesgue measure $0$,  it is of full dimension.
%\end{rem}
Let $(\Sigma_m,T)$ be the full shift over alphabet
$\{0,\cdots,m-1\}$, where   $\Sigma_m$ is endowed with the usual
metric  $ d_1(x,y)=m^{-|x\wedge y|}. $ If we take $\psi\equiv -\log
m$ and $\Psi=(n\psi)_{n=1}^{\infty}$, then $d_1=d_\Psi.$ Define a
map $\chi:\Sigma_m\to [0,1]$ as $ \chi(x)=\sum_{n=1}^\infty
{x_n}/{m^n}. $ Then $\chi$ is continuous and surjective. Consider
the IFS $\{f_j:j=0,\cdots,m-1\}$ defined as $f_j(x)=(x+j)/m$. It is
seen that the SOSC holds with $V=(0,1)$.
 Let
$ \widetilde Z_\infty:=\Big \{\sum_{j=1}^n x_jm^{-j}: n\in\N;
x_j=0,\cdots,m-1\Big \}\cup \{1\}. $ Define the dynamics
$\widetilde{g}$ on $\widetilde J_m=[0,1]\setminus\widetilde
Z_\infty$ as in the previous section. Then it is easy to check that
$\widetilde{g}={g_m}|_{\widetilde J_m}.$ Let $\Phi={\rm
id}_{[0,1]}$.  By theorem \ref{carpet} we get $ \dim_H
\mathcal{F}(\widetilde
J_m,g_m)=\sup\{\D_{{\Phi}}(\alpha):\alpha\in[0,1]\}. $ By the law of
large number applied to the measure of maximal entropy we get
$\D_\Phi(1/2)=1.$ We conclude by noticing that
$\mathcal{F}(J,g)=\prod_{i=1}^d\mathcal{F}([0,1],g_{m_i})$.

Next we consider concrete examples of carpets in the unit square.

\subsubsection*{{\bf  Heterogeneous carpets in the unit square}}
In order to fully illustrate our purpose, we consider an IFS
$S_0=\{f_1,\cdots,f_N\}$ in $\R^2$ made of contractive similitudes
without rotations such that the squares $f_i([0,1]^2)$ form a tiling
of $[0,1]^2$. All these situations have been determined in
\cite{Brooks}.
 In this way, $]0,1[^2$ can be chosen as the open set such that the SOSC holds,
  and the boundaries of the sets $f_i(]0,1[^2)$ have big intersections.
 The  picture in the left of Figure 1 give an example of this kind of
IFS. This IFS contains $15$ dilation maps, and the dynamics on this
attractor is highly non trivial.

 Let $\Phi$ denote ${\rm Id}_{\R^2}$. For each $\emptyset\neq S\subset S_0$, we denote by $J_S$ the attractor of the IFS $S$.
 The dynamics $\widetilde g_S$ defined on $\widetilde J_S$ is the restriction  of  $\widetilde g_{S_0}$
 to $\widetilde J_{S}$. The set $\mathcal{F}(\widetilde J_{S_0},\widetilde g_{S_0})$ is of full Hausdorff dimension,
  since $J_{S_0}=[0,1]^2$. If $S\neq S_0$, we have the variational formula $\dim_H \mathcal{F}
  (\widetilde J_{S},\widetilde g_{S})=\sup\{D_{\Phi}(\alpha):\alpha\in J_S\}$, and in general it
  is hard to know whether $ \mathcal{F}(\widetilde J_{S},\widetilde g_{S})$ is of full dimension or not
  in $J_S$. However, here are two  simple examples illustrating both possibilities.

 We consider the case of the regular tiling associated with the IFS $ S_0=\Big \{f_{i,j}:x\mapsto
 \frac{x}{3}+\frac{(i,j)}{3}:0\le i,j\le 2\Big \}$. Then, let $S_1=\{f_{0,0},f_{0,2},f_{2,0},f_{2,2}\}$ and $S_2=S_1\cup\{f_{1,1}\}$.
We claim that $\mathcal{F}(\widetilde J_{S_1},\widetilde g_{S_1})$
is not of full Hausdorff dimension,
 while $\mathcal{F}(\widetilde J_{S_2},\widetilde g_{S_2})$ is.

The simpler situation is that of $S_2$. In this case, $G=(1/2,1/2)$,
the center of symmetry of $J_{S_2}$
 is the fixed point of $f_{1,1}$ and it belongs to $L_{\Phi}$. Moreover, it is obvious that the uniform measure
 (or Parry measure) on $J_{S_2}$ is carried by the set $E_{\Phi}(G)$. This yields the result by Theorem~\ref{carpet},
  and $\dim_H \mathcal{F}(\widetilde J_{S_2},\widetilde g_{S_2})=\log (5)/\log (3)$.

In the case of $S_1$, the point $G$ is still the center of symmetry
of $J_{S_1}$, so  $\mathcal D_{\Phi}$ reaches its maximum at $G$.
However, $G$ does not belong to $J_{S_1}$.  Since $\Phi$ is H\"older
continuous and the tiling is regular, we know that $\mathcal
D_{\Phi}$ is strictly concave. By using the symmetry, one deduces
that the restriction of $\mathcal D_{\Phi}$ to $J_{S_1}$ reaches its
maximum at any of the four points $(1/3,1/3)$, $(1/3,2/3)$,
$(2/3,1/3)$ and $(2/3,2/3)$. This yields $\displaystyle \dim_H
\mathcal{F}(\widetilde J_{S_1},\widetilde
g_{S_1})=\D_\Phi((1/3,1/3))<\log (4)/\log (3)=\dim_H J_{S_1}$.

%%%%%%%%%%%%%%%%%%%%%%%%%%%%%%%%%%%55
\begin{figure}[h]
\setlength{\unitlength}{0.35cm}
\begin{picture}(36,12)(0,0)

\put(16,2){
\begin{picture}(9,9)\setlength{\unitlength}{0.5cm}
\multiput(0,0)(0,4){2}{\line(0,1){2}}
\multiput(2,0)(0,4){2}{\line(0,1){2}}
\multiput(4,0)(0,4){2}{\line(0,1){2}}
\multiput(6,0)(0,4){2}{\line(0,1){2}}
\multiput(0,0)(4,0){2}{\line(1,0){2}}
\multiput(0,2)(4,0){2}{\line(1,0){2}}
\multiput(0,4)(4,0){2}{\line(1,0){2}}
\multiput(0,6)(4,0){2}{\line(1,0){2}}
\put(2.75,2.75){$\bullet$} \put(2.75,2.1){G} \put(2.25,-1.5){IFS
$S_1$}
\end{picture}
}

\put(28,2){
\begin{picture}(9,9)\setlength{\unitlength}{0.5cm}
\multiput(0,0)(0,4){2}{\line(0,1){2}}
\multiput(2,0)(0,4){2}{\line(0,1){2}}
\multiput(4,0)(0,4){2}{\line(0,1){2}}
\multiput(6,0)(0,4){2}{\line(0,1){2}}
\multiput(0,0)(4,0){2}{\line(1,0){2}}
\multiput(0,2)(4,0){2}{\line(1,0){2}}
\multiput(0,4)(4,0){2}{\line(1,0){2}}
\multiput(0,6)(4,0){2}{\line(1,0){2}}
\multiput(2,2)(0,2){2}{\line(1,0){2}}
\multiput(2,2)(2,0){2}{\line(0,1){2}}
\put(2.75,2.75){$\bullet$} \put(2.75,2.1){G} \put(2.25,-1.5){IFS
$S_2$} 
\end{picture}
}

\put(-2,-2){
\begin{picture}(9,9)\setlength{\unitlength}{0.13cm}
\multiput(0,0)(0,37){2}{\line(1,0){37}}\multiput(0,0)(37,0){2}{\line(0,1){37}}%37
\multiput(0,0)(0,19){2}{\line(1,0){19}}\multiput(0,0)(19,0){2}{\line(0,1){19}}\put(8,9){19}% 19
\multiput(19,0)(0,18){2}{\line(1,0){18}}\multiput(19,0)(18,0){2}{\line(0,1){18}}\put(26,8){18}% 18
\multiput(19,18)(0,1){2}{\line(1,0){1}}\multiput(19,18)(1,0){2}{\line(0,1){1}}% 1
\multiput(20,18)(0,3){2}{\line(1,0){3}}\multiput(20,18)(3,0){2}{\line(0,1){3}}\put(21,18.5){3}% 3
\multiput(23,18)(0,8){2}{\line(1,0){8}}\multiput(23,18)(8,0){2}{\line(0,1){8}}\put(26.25,21){8}% 8
\multiput(31,18)(0,6){2}{\line(1,0){6}}\multiput(31,18)(6,0){2}{\line(0,1){6}}\put(33.5,20){6}% 6
\multiput(0,19)(0,18){2}{\line(1,0){18}}\multiput(0,19)(18,0){2}{\line(0,1){18}}\put(7.5,27){18}% 18
\multiput(18,19)(0,2){2}{\line(1,0){2}}\multiput(18,19)(2,0){2}{\line(0,1){2}}% 2
\multiput(18,21)(0,5){2}{\line(1,0){5}}\multiput(18,21)(5,0){2}{\line(0,1){5}}\put(20,23){5}% 5
\multiput(18,26)(0,11){2}{\line(1,0){11}}\multiput(18,26)(11,0){2}{\line(0,1){11}}\put(22.5,30){11}% 11
\multiput(31,24)(0,1){2}{\line(1,0){1}}\multiput(31,24)(1,0){2}{\line(0,1){1}}% 1
\multiput(31,25)(0,1){2}{\line(1,0){1}}\multiput(31,25)(1,0){2}{\line(0,1){1}}% 1
\multiput(32,24)(0,5){2}{\line(1,0){5}}\multiput(32,24)(5,0){2}{\line(0,1){5}}\put(34,25.5){5}% 5
\multiput(29,29)(0,8){2}{\line(1,0){8}}\multiput(29,29)(8,0){2}{\line(0,1){8}}\put(32.5,32){8}% 8
\put(30,26.5){3}
\end{picture}
}
\end{picture}
\caption{}
\end{figure}
%%%%%%%%%%%%%%%%%%%%%%%%%%%%%%

\subsection{Localized results for Gibbs measures}\label{appli2}
Let $\{f_1,\cdots,f_m\}$ be a homogenous self-similar IFS in
$\mathbb C$ satisfying the {\it  strong separation condition}, that is, each
function $f_j$ has the form $f_j(z)=a_jz+b_j$ where
$0<\rho=|a_j|<1,$ and there exists a topological closed  disk $D$ such that
$f_j(D)\subset D$  and  the $f_j(D)$ are pairwise disjoint. There is a natural
coding map $\chi:\Sigma_m\to J$. Moreover if we define
$\psi(x)\equiv\log \rho$ for $x\in\Sigma_m$, and
$\Psi=(S_n\psi)_{n=1}^{\infty}$, then
$\chi:(\Sigma_m,d_\Psi)\to(J,|\cdot|)$ is a bi-Lipschitz
homeomorphism.

 Let $\phi: J\to \R$ be  continuous and define $\tilde \phi=\phi\circ
 \chi.$ By subtracting a constant potential if necessary,  we can assume $P(T,\tilde\phi)=0.$ There
exists a weak Gibbs measure $\tilde\mu$ on $\Sigma_m$ (see \cite{K})
 such that
 $$d_{\tilde\mu}(x):=\lim_{r\to 0}\frac{\log \mu(B(x,r))}{\log
 r}=\lim_{n\to\infty}\frac{S_n\tilde\phi(x)}{n\log\rho}$$
 in the sense that either both the limits do not exist, either they
 exist and are equal. Define $\mu:=\chi_\ast(\tilde\mu)$. By the bi-Lipschitz property of $\chi$ and the strong separate condition, we can
 easily conclude that
 $d_\mu(y)=\lim_{n\to\infty}S_n\phi(y)/n\log\rho$ for any $y\in
 J.$ Let $\Phi=(S_n\phi)_{n=1}^\infty.$ If we define $E_\mu(\alpha)=\{y\in J: d_\mu(y)=\alpha\}$, then
 we get $E_\Phi(\alpha)=E_\mu(\alpha/\log \rho)$ for any $\alpha\in L_\Phi.$

 By  applying Theorem
\ref{appli-fun-level-one-sided} for $d=1$, we have the following
property regarding the local property of
 weak Gibbs measure:

 \begin{coro}\label{weakgibbs}
Let $\mu$ be the weak Gibbs measure associated with $\phi.$ Then the
set of all possible local dimension for $\mu$ is the interval $
L_\Phi/\log\rho.$ Assume $\xi:J\to\R$ is
 continuous and $\xi(J)\subset  L_\Phi/\log\rho$,
 then
 $$
\dim_H\{x\in J: d_\mu(x)=\xi(x)\}=\sup\{\dim_H
E_\mu(\alpha):\alpha\in\xi(J)\}.
 $$
 \end{coro}

 Now let $\omega$ stand for the harmonic measure on $J$. It is well
 known that (see for example the survey paper \cite{Ma}) there
 exists a H\"older continuous function $\phi:J\to\R$ such that $w\asymp
 \mu$, where $\mu$ is the equilibrium state of $\phi.$
By a direct application of the above corollary we have the following
property:
 \begin{coro}\label{harmonic}
Assume $\omega$ is the harmonic measure on $J$ and $I$ is the set of
all possible local dimension for $\omega$.  Assume $\xi:J\to\R$ is
 continuous and $\xi(J)\subset I$.
 Then
 $$
\dim_H\{x\in J: d_\omega(x)=\xi(x)\}=\sup\{\dim_H
E_\omega(\alpha):\alpha\in\xi(J)\}.
 $$
 \end{coro}

\noindent
{\bf Final remark.} At least when $d=1$, it is not difficult to extend the results obtained in
 this paper by considering $\Upsilon=(\gamma_n)_{n\ge 1}\in \C_{aa}^+(\Sigma_A,T)$
 and the more general level sets $E^\Psi_{\Phi/\Upsilon}(\xi)=\{x\in \Sigma_A: \lim_{n\to\infty}\phi_n(x)/\gamma_n(x)=\xi(x)\}$;
  when $\xi$ is constant, such sets have been considered in the contexts examined in \cite{BSS,BD}.
   The formula is that if the continuous function $\xi$ takes values in the set
   $L_{\Phi/\Upsilon}=\{\Phi_*(\nu)/\Upsilon_*(\nu):\nu\in\mathcal{M}(\Sigma_A,T)\}$,
   then $\dim_H^\Psi (E^\Psi_{\Phi/\Upsilon}(\xi))=\sup\{-h_\nu(T)/\Psi_*(\nu):
   \nu\in \mathcal{M}(\Sigma_A,T),\ \Phi_*(\nu)/\Upsilon_*(\nu)\in \xi(\Sigma_A)\}$.
   When $\Upsilon =-\Psi$, this can be applied to the local dimension of Gibbs measures
    associated with H\"older potentials $\varphi$ on any $C^1$ conformal repeller of a map $f$,
    since in this case we know from \cite{P} that such a measure is doubling so that the
     local dimension is directly related to the asymptotic behavior of  $S_n\varphi/S_n(-\log \|Df\|)$.
      Consequently, Corollary~\ref{harmonic} can be extended to harmonic measure on more general conformal repellers (see \cite{Ma}).

\section{Proofs of Propositions~\ref{dim-L-phi}, ~\ref{dim-formular-1}, \ref{regularity} and \ref{gene-thermo}}\label{proofs}

\subsection{Proof of proposition \ref{dim-L-phi}}
Suppose that $\Psi,\Psi^\prime\in \C_{aa}(X,T)$ and
$\|\Psi-\Psi^\prime\|_{\rm lim}=0$. Then  $\mu_\ast(\Psi)=\mu_\ast(\Psi^\prime)$ for every $\mu\in\M(X,T)$. From this and
the definition  of $L_\Phi$,
 we can easily show that  if $L_\Phi$
is of dimension $d$ then $\bar\Phi^1,\cdots,\bar\Phi^d$
are linearly independent.

To show the other direction, at first we assume $d=1$.  By \cite{FH}
lemma 3.5, we have $L_\Phi=[\beta_1,\beta_2]$, where
$\displaystyle\beta_1=\lim_{n\to\infty}\inf_{x\in X}\phi_n(x)/n$ and
$\displaystyle\beta_2=\lim_{n\to\infty}\sup_{x\in X}\phi_n(x)/n$. We
need to show that if $\Phi\not\sim 0$, then $L_\Phi$ is a non
degenerate interval in $\R.$ Otherwise, we get $\beta_1=\beta_2$,
then we get $\|\Phi-\beta_1\|_{\rm lim}=0,$ thus  $\Phi\sim 0$,
which is a contradiction.

Now assume  $d>1$ and $\bar\Phi^1,\cdots,\bar\Phi^d$ are
linearly independent.  If $L_\Phi$ has dimension strictly less than
$d$, then there exists a non-zero vector $p\in \R^d$ such that
$p\cdot L_\Phi$ is a singleton, that is $L_{p\cdot\Phi}=p\cdot
L_\Phi$ is a singleton. By what has been shown, we conclude that
$p\cdot\Phi\sim 0$, that is $p_1\bar\Phi^1+\cdots+
p_d\bar\Phi^d=0,$ which is a contradiction. \hfill $\Box$

\subsection{Proof of Proposition ~\ref{dim-formular-1}} We need some preliminary results.

\begin{lem}\label{multiplictive}
Let  $\Phi\in \C_{aa}(\Sigma_A,T)$. Let $C=C(\Phi)$.
\begin{enumerate}
\item $\lim_{n\to\infty}\|\Phi\|_n/n=0.$

\item  For any $u,v\in \Sigma_{A,\ast} $ such that $uv\in  \Sigma_{A,\ast}$ we have
$$ \exp(-C-\|\Phi\|_{|u|})\Phi[u]\Phi[v]\leq \Phi[uv]\leq \exp(C)\Phi[u]\Phi[v].$$

\item  For $w=u_1w_1\cdots u_nw_nu_{n+1}\in\Sigma_{A,\ast}$, let
$k=\sum_{j=1}^{n+1}|u_j|$. We have
\begin{equation}\label{222}
\exp(-2nC+k\Phi_{\min})\prod_{j=1}^{n}\Phi[w_j]\exp(-\|\Phi\|_{|w_j|})
\leq \Phi[w]\leq \exp(2nC+k\Phi_{\max})\prod_{j=1}^{n}\Phi[w_j].
\end{equation}

\item If  $\Phi\in \C_{aa}^-(\Sigma_A,T)$, then $\Phi[v]\leq\Phi[u]$ for
$u\prec v$.
\end{enumerate}
\end{lem}

%%%%%%%%%%%%%%%%%%%%%%%%%%%%%%%%%%%%%%%%%%%%%%%%%%%%%%%%%%%%%%%%%%%%%%%%%%%%%%%%%%%%%%%%%%%%%%%%%%%%%%%%%%

\proof\ (1) \  Fix $k\in\N$, let $n=kp+l$ with $0\leq l< k$. Then by
almost additivity we get
$$
\sum_{j=0}^{p-1}\phi_k(T^{kj}x)+\phi_l(T^{pk}x)-pC\leq\phi_n(x)\leq
\sum_{j=0}^{p-1}\phi_k(T^{kj}x)+\phi_l(T^{pk}x)+pC.
$$
This yields $ \|\Phi\|_n\leq
\sum_{j=1}^p\|\phi_k\|_{kj}+2k\|\Phi\|+2pC$ and $
\frac{\|\Phi\|_n}{n}\leq
\frac{\sum_{j=1}^p\|\phi_k\|_{kj}}{kp}+\frac{2\|\Phi\|}{p}+\frac{2C}{k}.
$ When $k$ is fixed, since $\phi_k$ is continuous, we know that
$\sum_{j=1}^p\|\phi_k\|_{kj}/p\to 0$ as $p\to\infty.$ Then the
result follows easily.

(2)  Let $|u|=n$ and $|v|=k$. Given $x\in[uv]$ we have $x\in [u]$ and
$T^nx\in [v]$ and $\phi_{n+k}(x)\leq C+\phi_{n}(x)+\phi_{k}(T^nx).$
Thus $\sup_{x\in [uv]}\phi_{n+k}(x)\leq C+\sup_{x\in [u]}\phi_{n}(x)+\sup_{y\in [v]}\phi_{k}(y)$.
Consequently we have $\Phi[uv]\leq \exp(C)\Phi[u]\Phi[v]$.

On the other hand take  $x_0\in [v]$ and $y_0\in [u]$ such that
$\phi_{k}(x_0)=\sup_{x\in [v]}\phi_{k}(x) $ and
$\phi_{n}(y_0)=\sup_{y\in [u]}\phi_{n}(y) $. Let $\widetilde{x}=ux_0,$
then we have
$$
\sup_{x\in[uv]}\phi_{n+k}(x)\geq \phi_{n+k}(\widetilde{x})
\geq\phi_{n}(\widetilde{x})+\phi_{k}(x_0)-C \geq
\phi_{n}(y_0)+\phi_{k}(x_0)-C-\|\Phi\|_n.
$$
Consequently,  $\exp(-C-\|\Phi\|_n)\Phi[u]\Phi[v]\leq
\Phi[uv].$

(3) It is similar to the proof of (2).

%(3)\ Let $k_0=0$, $k_1=|u_1|, k_2=|u_1w_1|, \cdots,
%k_{2n}=|u_1w_1\cdots u_nw_n|$. Fix $x\in[w]=[u_1w_1\cdots
%u_nw_nu_{n+1}]$. We have $T^{k_{2(j-1)}}(x)\in [u_j]$ and
%$T^{k_{2j-1}}(x)\in[w_j]$. Thus
%\begin{eqnarray*}
%\phi_{|w|}(x)&\leq&
%\sum_{j=1}^{n+1}\phi_{|u_j|}(T^{k_{2(j-1)}}(x))+\sum_{j=1}^{n}\phi_{|w_j|}(T^{k_{2j-1}}(x))+2nC\\
%&\leq&
%\sum_{j=1}^{n+1}|u_j|\Phi_{\max}+\sum_{j=1}^{n}\phi_{|w_j|}(T^{k_{2j-1}}(x))+2nC\\
%&=& \sum_{j=1}^{n}\phi_{|w_j|}(T^{k_{2j-1}}(x))+2nC+k\Phi_{\max}.
%\end{eqnarray*}
%From this we easily get $\Phi[w]\leq
%\exp(2nC+k\Phi_{\max})\prod_{j=1}^{n}\Phi[w_j].$ The other inequality
%follows similarly.

(4)\ If $\Phi\in \C_{aa}^-(\Sigma_A,T)$, then $\phi_n(x)\leq
\phi_m(x)$ for any $m\leq n$. Since $[v]\subset[u]$ for $u\prec v$,
by definition we get $\Phi[v]\leq\Phi[u]$. \hfill$\Box$

\begin{rem}{\rm  \label{vector-potential}
By repeating the proof of (1), one can show that we still have $\|\Phi\|_n/n\to 0$
as $n\to\infty$ for
$\Phi\in\C_{aa}(\Sigma_A,T,d)$.
}\end{rem}

%%%%%%%%%%%%%%%%%%%%%%%%%%%%%%%%%%%%%%%%%%%%%%%%%%%%%%%%%%%%%%%%%%%%%%%%%%%%%%%%%%%%%%%%%%%%%%%%%%%%%%%%%%%%%%
%%%%%%%%%%%%%%%%%%%%%%%%%%%%%%%%%%%%%%%%%%%%%%%%%%%%%%%%%%%%%%%%%%%%%%%%%%%%%%%%%%%%%%%%%%%%%%%%%%%%%%%%%%%%%
%%%%%%%%%%%%%%%%%%%%%%%%%%%%%%%%%%%%%%%%%%%%%%%%%%%%%%%%%%%%%%%%%%%%%%%%%%%%%%%%%%%%%%%%%%%%%%%%%%%%%%%%%%%%%555555555555
%%%%%%%%%%%%%%%%%%%%%%%%%%%%%%%%%%%%%%%%%%%%%%%%%%%%%%%%%%%%%%%%%%%%%%%%%%%%%%%%%%%%%%%%%%%%%%%%%%%%%%%%%%%%%555555555555

\begin{lem} \label{Moran-covering}
Let $\Psi\in  \C_{aa}^-(\Sigma_A,T).$
\begin{enumerate}
\item Let $C_1(\Psi)=1/|\Psi_{\min}|$ and
$C_2(\Psi)=1+1/|\Psi_{\max}|.$ For any $w\in{\mathcal B}_n(\Psi)$ we
have
\begin{equation}\label{constant}
C_1(\Psi)n\leq|w|\leq C_2(\Psi)n.
\end{equation}

\item For any $w\in{\mathcal B}_n(\Psi)$ we have
\begin{equation}\label{444}
\exp(-C(\Psi)-\|\Psi\|_{|w|}+\Psi_{\min})e^{-n}\leq\Psi[w]\leq
e^{-n}.
\end{equation}

\item The balls in  $\{[w]: w\in{\mathcal B}_n(\Psi)\}$ are pairwise
disjoint.

\item  If $u\prec v$ are such that $u\in {\mathcal B}_{n_1}(\Psi)$ and
$v\in {\mathcal B}_{n_2}(\Psi)$, then
$$|v|-|u|\leq \frac{\Psi_{\min}-\|\Psi\|_{|v|}-(n_2-n_1)-2C}{\Psi_{\max}}.$$
\end{enumerate}
\end{lem}

\proof\  (1) By $(\ref{max-min})$ we have
$e^{|w|\Psi_{\min}}\leq\Psi[w]\leq e^{|w|\Psi_{\max}}$ for any $w\in
\Sigma_{A,\ast}$. If $w\in {\mathcal B}_n(\Psi)$, we have
$e^{\Psi_{\min}|w|}\leq \Psi[w]\leq e^{-n}<\Psi[w^\ast]\leq e^{\Psi_{\max}(|w|-1)}.$
Thus $C_1(\Psi)n=-n/\Psi_{\min}\leq|w|\leq
n(1-1/\Psi_{\max})=C_2(\Psi)n.$

(2) By definition $\Psi[w]\leq e^{-n}$.  Assume $|w|=k$.  By Lemma
\ref{multiplictive}(2)
$$\Psi[w]=\Psi[w^\ast w_k]\geq \exp(-C-\|\Psi\|_{k-1})\Psi[w^\ast]\Psi[w_k]\geq \exp(-C-\|\Psi\|_k-n+\Psi_{\min}).$$

(3)  $d_\Psi$ is ultra-metric.

%Assume $w,\widetilde{w}\in {\mathcal B}_n(\Psi)$ and $[w]\cap
%[\widetilde{w}]\ne \emptyset$. Take $x\in [w]\cap [\widetilde{w}]$. By
%proposition \ref{metric-sym}, $B(x,e^{-n})=[x|_k]$ where $k$ is the
%unique number such that $\Psi[x|_{k-1}]>e^{-n}$ and $\Psi[x|_k]\leq
%e^{-n}$. Since $x\in [w]$, we know that $\Psi[x|_j]=\Psi[w|_j]$ for
%$j=1,\cdots, |w|$. Thus $B(x,e^{-n})=[w]$. The same argument shows
%that $B(x,e^{-n})=[\widetilde{w}]$. So we have $[w]=[\widetilde{w}]$.

(4) Write $v=uw$. Then $|w|=|v|-|u|$ and
$$e^{-C-\|\Psi\|_{|v|}-n_2+\Psi_{\min}}\leq \Psi[v]=\Psi[uw]\leq e^{C}\Psi[u]\Psi[w]\leq e^{C-n_1+|w|\Psi_{\max}}.$$
 \hfill$\Box$

Given $\Phi\in \C_{aa}(\Sigma_A,T,d)$ and two constants $C_2\geq C_1>0$.
For each $n\in \N$ we define
\begin{equation}\label{variation}
\|\Phi\|_{n}^\star:=\max\{\|\Phi\|_l:C_1n\leq l\leq C_2n\}.
\end{equation}
By remark \ref{vector-potential},  we have
$\|\Phi\|_{n}^\star/n\to 0$ when $n\to\infty.$

%%%%%%%%%%%%%%%%%%%%%%%%%%%%%%%%%%%%%%%%%%%%%%%%%%%%%%%%%%%%%%%%%%%%%%%%%%%%%%%%%%%%%%%%%%%%%%%%%%%%%%%%%%%%%

%%%%%%%%%%%%%%%%%%%%%%%%%%%%%%%%%%%%%%%%%%%%%%%%%%%%%%%%%%%%%%%%%%%%%%%%%%%%%%%%%%%%%%%%%%%%%%%%%%%%%%5555

\noindent {\bf Proof of Proposition \ref{dim-formular-1}}.\ We fix
$\Phi$ and $\Psi$ and write $F(\alpha,n,\epsilon),
f(\alpha,n,\epsilon)$ and $\Lambda$ for
$F(\alpha,n,\epsilon,\Phi,\Psi), f(\alpha,n,\epsilon,\Phi,\Psi)$ and
$\Lambda_\Phi^\Psi$ to simplify the notation.

At first we show that $\log f(\alpha,n,\epsilon)$, as a sequence of
$n$, has a kind of subadditivity property. More precisely, for any
$\epsilon>0$, there exist an $N\in \N$ and $\beta_n>0$ such that
$\log\beta_n=o(n)$ and
$f(\alpha,n,\epsilon)^p\leq
\beta_n^pf(\alpha,(n+\widetilde{c})p,2\epsilon)$ for any $ n\geq N,$
and any $ p\geq 1,$
where $\widetilde{c}=\lfloor-p_0\Psi_{\max}-2C(\Psi)\rfloor$. Recall that $p_0$ is a fixed positive integer such that $A^{p_0}>0.$

 In fact for
$w_1,\cdots,w_p\in F(\alpha,n,\epsilon),$ let
$w=\overline{w}_1\cdots \overline{w}_p$, where
$\overline{w}_j=w_ju_j$ with $u_j\in \Xi$ such that $w_ju_jw_{j+1}$
is admissible. Recall (see \eqref{normn}) that for  any cylinder
$[u]$ and any $x,\widetilde{x}\in [u]$, we have
$|\psi_{|u|}(x)-\psi_{|u|}(\widetilde{x})|\leq \|\Psi\|_{|u|}.$
Now for any $x\in [w]$,  let $ s_0=0, s_k=\sum_{j=1}^{k}(|w_j|+p_0)\
(1\leq k\leq p)$ and define $x^{k}=T^{s_{k-1}}x$. We have
$|w|=s_{p}$ and $x^k\in[w_k]$ for $k=1,\cdots,p.$
 Take
$y^k\in[w_k]$ such that $\Psi[w_k]=\exp ({\psi_{|u|}(y^k)})$. Then
\begin{eqnarray*}
\psi_{|w|}(x)
&\geq& \sum_{k=1}^p \psi_{|w_k|}(x^k)+p_0p\Psi_{\min}-(2p-1)C(\Psi)\\
&\geq& \sum_{k=1}^p \psi_{|w_k|}(y^k)-\sum_{k=1}^p
\|\Psi\|_{|w_k|}+p(p_0\Psi_{\min}-2C(\Psi)).
\end{eqnarray*}
 Thus
by Lemma \ref{Moran-covering}(2)
\begin{eqnarray*}
\Psi[w]&\geq& \exp(\psi_{|w|}(x))\\
&\geq& \exp\Big (\sum_{k=1}^p \psi_{|w_k|}(y^k)-\sum_{k=1}^p \|\Psi\|_{|w_k|}+p(p_0\Psi_{\min}-2C(\Psi))\Big )\\
&=&(\prod_{k=1}^{p}\Psi[w_k])\exp\Big (-\sum_{k=1}^p \|\Psi\|_{|w_k|}+p(p_0\Psi_{\min}-2C(\Psi))\Big )\\
&\geq&\exp(-pn)\exp(-2\sum_{k=1}^p \|\Psi\|_{|w_k|}+p((p_0+1)\Psi_{\min}-3C(\Psi)))\\
&\geq&\exp(-pn)\exp\Big (p((p_0+1)\Psi_{\min}-3C(\Psi)-2p\|\Psi\|_{n}^\star)\Big)\geq \exp(-p(n+c_1(n)),
\end{eqnarray*}
where  $c_1(n)=-(p_0+1)\Psi_{\min}+3C(\Psi)+2\|\Psi\|_{n}^\star>0$
and $\|\Psi\|_{n}^\star$ is defined as in $(\ref{variation})$ with
constants $C_2(\Psi)\geq C_1(\Psi)>0$. Lemma \ref{multiplictive}(1) yields $c_1(n)/n\to 0.$

 By Lemma \ref{multiplictive}(3) we also
have
$$\Psi[w]\leq \exp(2pC(\Psi)+p_0p\Psi_{\max})(\prod_{k=1}^{p}\Psi[w_k])\leq \exp(-p(n-p_0\Psi_{\max}-2C(\Psi)).$$
By definition of $\widetilde c$, there exists $u\in {\mathcal
B}_{p(n+\widetilde{c})}(\Psi)$ such that $u$ is the prefix of $w$.
Write ${w}=uw^\prime$.

\noindent {\bf Claim:}\ \ $|w^\prime|\leq p(ac_1(n)+b)$ for some
constant $a,b>0.$

Indeed, we have
$e^{-p(n+c_1(n))}\leq\Psi[{w}]\leq e^C\Psi[u]\Psi[w^\prime]\leq e^Ce^{-p(n-p_0\Psi_{\max}-2C)} e^{|w^\prime|\Psi_{\max}}.$
Thus $|w^\prime|\leq p (c_1(n)+p_0\Psi_{\max}+3C)/(-\Psi_{\max}).$

Now since $w_k\in F(\alpha,n,\epsilon)$ we can find
$x_{k}\in[w_{k}]$ such that
$|\frac{\phi_{|w_k|}(x_k)}{|w_k|}-\alpha|<\epsilon.$ Take $x\in
[w]$; in particular, $x\in[u]$. Define $s_k$ and $ x^{k}$ as above.
We have $|{w}|=s_{p}$ and $x^k\in[w_k]$ for $k=1,\cdots,p.$ By
almost additivity, we get
$$
\phi_{|u|}(x)+\phi_{|w^\prime|}(T^{|u|}x)-C(\Phi)\leq\phi_{|w|}(x)\leq
\phi_{|u|}(x)+\phi_{|w^\prime|}(T^{|u|}x)+C(\Phi)
$$
(this is a vector inequality). Notice
that if $\beta_1,\beta_2\in \R^d$ are such that $\beta_1>0$ and
$-\beta_1\leq \beta_2\leq \beta_1$, then $|\beta_2|\leq |\beta_1|.$
Thus we have
\begin{eqnarray*}
\phi_{|u|}(x)&=&\phi_{|w|}(x)+\eta_0=\sum_{k=1}^{p}\phi_{|w_ku_k|}(x^k)+\eta_1+\eta_0=\sum_{k=1}^{p}\phi_{|w_k|}(x^k)+\eta_2+\eta_1+\eta_0\\
&=&\sum_{k=1}^p\phi_{|w_k|}(x_k)+\eta_3+\eta_2+\eta_1+\eta_0=(\sum_{k=1}^p|w_k|)\alpha+\eta_4+\eta_3+\eta_2+\eta_1+\eta_0,
\end{eqnarray*}
where $|\eta_0|\leq |w^\prime|\|\Phi\|+|C(\Phi)|\leq
p(ac_1(n)+b)\|\Phi\|+|C(\Phi)|$; $\eta_1\leq (p-1)|C(\Phi)|$;
$|\eta_2|\leq p(p_0\|\Phi\|+|C(\Phi)|)$; $|\eta_3|\leq
\sum_{k=1}^{p}\|\Phi\|_{|w_k|}\leq p\|\Phi\|_{n}^\star$;
$|\eta_4|\leq (\sum_{k=1}^p|w_k|)\epsilon.$ Since
$s_{p}=\sum_{k=1}^p|w_k|+p_0p$ and $|w_k|\geq C_1 n$, we have
$s_p\geq C_1np$, and
\begin{eqnarray*}
\Big |\frac{\phi_{|u|}(x)}{|u|}-\alpha\Big |\leq\frac{|((\sum_{k=1}^p|w_k|)-|u|)\alpha|+|\eta_4|+|\eta_3|+|\eta_2|+|\eta_1|+|\eta_0|}{|u|}.
\end{eqnarray*}
Moreover, $|u|=s_p-|w^\prime|\geq pC_1n-p(ac_1(n)+b)$  and $c_1(n)/n,
\|\Phi\|_{n}^\star/n \to 0$,  so that we can choose $N(\epsilon)$ big enough
such that
$|\phi_{|u|}(x)/|u|-\alpha|\leq 2\epsilon$ when $n\geq N(\epsilon)$. Consequently $u\in
F(\alpha,p(n+\widetilde{c}),2\epsilon).$ Thus, $f(\alpha,p(n+\widetilde{c}),2\epsilon)\geq f(\alpha,n,\epsilon)^p/m^{p(ac_1(n)+b)}.$
If we take $\beta_n=m^{ac_1(n)+b}$, then we get the desired
subadditivity.

Next we show that $$\lim_{\epsilon\to 0}\liminf_{n\to\infty}\frac{\log
f(\alpha,n,\epsilon)}{n}= \lim_{\epsilon\to
0}\limsup_{n\to\infty}\frac{\log f(\alpha,n,\epsilon)}{n}.$$

Note that both limit exist since $f(\alpha,n,\epsilon)$ is an
increasing function in  the variable $\epsilon.$ Denote by $\beta$
the left-hand side limit. Then for any $\delta>0$, there exists
$\epsilon_0>0$ such that $\displaystyle \liminf_{n\to\infty}\log
f(\alpha,n,\epsilon_0)/n<\beta+\delta. $
 Fix
$\delta>0$ and $\epsilon_0>0$ as above. To show the equality we only need to show
that
$
\limsup_{n\to\infty}{\log (f(\alpha,n,\epsilon_0/4))}/{n}\leq
\beta+\delta.
$
Fix $n\in \N$. Take a sequence of integers $n_k\to\infty$ such that
$ f(\alpha,n_k,\epsilon_0)<e^{n_k(\beta+\delta)}$ for any $ k\in \N.
$ For each $k,$ write $n_k=(n+\widetilde{c})p_k-l_k$ with $0\leq
l_k<n+\widetilde{c}.$ By the subadditivity property, we have
$f(\alpha,n,\epsilon_0/4)^{p_k}\leq \beta_n^{p_k}f(\alpha,(n+\widetilde{c})p_k,\epsilon_0/2).$
If $w=w_1w_2$ is such that $w_1\in {\mathcal B}_l(\Psi)$ and $w\in {\mathcal
B}_{l+s}(\Psi)$ with $1\leq s\leq (n+\widetilde{c})$, then by Proposition
\ref{Moran-covering}(4) we have
\begin{equation}\label{7}
|w_2|\leq
\frac{\Psi_{\min}-\|\Psi\|_{|w|}-(n+\widetilde{c})-2C(\Psi)}{\Psi_{\max}}.
\end{equation}
Thus $|w_2|/|w|\to 0$ when $|w|\to\infty.$  Choose $l_0$ large
enough so that when $l\geq l_0$ we have
\begin{equation}\label{6}
\frac{|w|}{|w_1|}\leq \frac{3}{2},\ \ \frac{|C(\Phi)|}{|w_1|}\leq
\frac{\epsilon_0}{8}\ \ \text{and }\ \ \frac{|w_2|}{|w_1|}\leq
\frac{\epsilon_0}{8(\|\Phi\|+|\alpha|)}.
\end{equation}

Let $k_0$  such that $(n+\widetilde{c})p_{k_0}\geq l_0+(n+\widetilde{c})$.
Let $k\geq k_0.$  Fix $w\in
F(\alpha,(n+\widetilde{c})p_k,\epsilon_0/2)$. There exists $x\in
[w]$ such that $|\phi_{|w|}(x)-|w|\alpha|\leq |w|\epsilon_0/2.$
Let $w_1\prec w$ such that $\Psi[w_1]\leq e^{-n_k} $ and
$\Psi[w_1^\ast]>e^{-n_k}$. Thus $[w_1]\in {\mathcal B}_{n_k}(\Psi)$. Write
$w=w_1w_2.$ By $(\ref{6})$
 we have
\begin{eqnarray*}
|\phi_{|w_1|}(x)-|w_1|\alpha|&\leq& |\phi_{|w_1|}(x)-\phi_{|w|}(x)|+|\phi_{|w|}(x)-|w|\alpha|+|w_2||\alpha|\\
&\leq&|\phi_{|w|}(x)-|w|\alpha|+|w_2|(\|\Phi\|+|\alpha|)+|C(\Phi)|\\
&\leq &\frac{|w|\epsilon_0}{2}+|w_2|(\|\Phi\|+|\alpha|)+|C(\Phi)|\\
&\leq&\frac{3|w_1|\epsilon_0}{4}+\frac{|w_1|\epsilon_0}{8}+\frac{|w_1|\epsilon_0}{8}=|w_1|\epsilon_0,
\end{eqnarray*}
which means that $w_1\in F(\alpha,n_k,\epsilon_0)$. Moreover by
$(\ref{7})$, we have $|w_2|/|w|\to 0$ when $|w|\to\infty.$ Thus we
can find a sequence $\gamma_k$ such that $|w_2|\leq
\gamma_k=o(|w|)=o(n_k)=o(p_k)$. We can conclude that
$f(\alpha,(n+\widetilde{c})p_k,\epsilon_0/2)\leq m^{\gamma_k}f(\alpha,n_k,\epsilon_0).$
This yields
$$f(\alpha,n,\epsilon_0/4)\leq \beta_n m^{\gamma_k/p_k}f(\alpha,n_k,\epsilon_0)^{1/p_k}\leq \beta_n m^{\gamma_k/p_k}e^{n_k(\beta+\delta)/p_k}.$$
Letting $k\to\infty$ we get $f(\alpha,n,\epsilon_0/4)\leq \beta_n e^{(n+\widetilde{c})(\beta+\delta)}$. Then, letting  $n\to\infty$ we have
$\limsup_{n\to\infty}{\log (f(\alpha,n,\epsilon_0/4))}/{n}\leq \beta+\delta$.

Next we show the upper semi-continuity of $\Lambda(\alpha)$. Let
$\alpha\in L_\Phi$. For any $\eta>0$ there is $\epsilon >0$ such
that
$\liminf_{n\to\infty}\frac{f(\alpha,n,\epsilon)}{n}< \Lambda(\alpha)+\eta.$
Let $\beta\in L_\Phi$ with $|\beta-\alpha|<\epsilon/3.$ Given $w\in
F(\beta,n,\epsilon/3)$, there exists $x\in[w]$ such that
$|\phi_{|w|}(x)/|w|-\beta|\leq \epsilon/3.$ Hence
$|\phi_{|w|}(x)/|w|-\alpha|\leq |\phi_{|w|}(x)/|w|-\beta|+|\beta-\alpha|<\epsilon,$
which means $w\in  F(\alpha,n,\epsilon)$. This proves that $
F(\beta,n,\epsilon/3)\subset F(\alpha,n,\epsilon)$. It follows that
$f(\beta,n,\epsilon/3)\leq f(\alpha,n,\epsilon)$, therefore
$$
\Lambda(\beta)\leq
\liminf_{n\to\infty}\frac{f(\beta,n,\epsilon/3)}{n}\leq
\liminf_{n\to\infty}\frac{f(\alpha,n,\epsilon)}{n}<
\Lambda(\alpha)+\eta.
$$
This establishes the upper semi-continuity of $\Lambda$ at $\alpha.$

By essentially repeating the  proof above (in fact it is much
easier), we can show
$$\liminf_{n\to\infty}\frac{\log \#{\mathcal B}_n(\Psi)}{n}=
 \limsup_{n\to\infty}\frac{\log \#{\mathcal B}_n(\Psi)}{n}.$$ We denote the  limit by $D(\Psi)$. By
 $(\ref{constant})$ there exist constants $0<C_1=C_1(\Psi)<C_2=C_2(\Psi)$ such
 that for any $w\in{\mathcal B}_n(\Psi)$, $C_1n\leq |w|\leq C_2n. $ This yields
 $\#\Sigma_{A,[C_1n]}\leq \#{\mathcal B}_n(\Psi)\leq \#\Sigma_{A,[C_2n]}$ and
 $C'_1\log m\leq D(\Psi)\leq C'_2\log m.$ \hfill$\Box$

Now we come to the weak concavity of the function $\Lambda_\Phi^\Psi.$

\subsection{Proof of Proposition \ref{regularity}}

Let $A\subset\R^d$. We say that $x\in A$ is a {\it local cone point}, or an $\epsilon$-cone point,
if there exists $\epsilon>0$ such that for any $y\in A\cap
B(x,\epsilon)$, the interval $[x,y_\epsilon]\subset A,$ where
$y_\epsilon:=x+\epsilon(y-x)/|y-x|$.
\begin{lem} \label{distorsion-concave}
Let $A\subset \R^d$ be a  convex set and  $h: A\to \R$ be a bounded
weakly concave function. Then $h$ is lower semi-continuous at each
local cone point of $A$. Especially  $h$ is lower semi-continuous on
$\mathrm{ri}(A)$ and on any closed interval $I\subset A$. It is
lower semi-continuous on $A$ if $A\subset \R^d$ is a convex closed
polyhedron.
\end{lem}

\proof \  Let $\beta\in A$ be a $\epsilon$-cone point of $A$ for
some $\epsilon>0$. Suppose that $h$ is not lower
semi-continuous at $\beta$. Thus we can find $\eta>0$ and
$\alpha_n\in A\cap B(\beta,\epsilon)$ such that $\alpha_n\to\beta$
and $h(\alpha_n)\leq h(\beta)-\eta.$ Define
$\alpha_n^\prime=\beta+\epsilon(\alpha_n-\beta)/|\alpha_n-\beta|$,
then $\alpha_n^\prime\in A$ since $\beta$ is a $\epsilon$-cone
point. Let $\lambda_n\in [0,1]$ such that
$$
\alpha_n=\frac{\lambda_n\gamma_{1}(\alpha_n^\prime,\beta)\alpha_n^\prime+(1-\lambda_n)\gamma_{2}(\alpha_n^\prime,\beta)\beta}
{\lambda_n\gamma_1(\alpha_n^\prime,\beta)+(1-\lambda_n)\gamma_2(\alpha_n^\prime,\beta)}.
$$
Since $\gamma_1,\gamma_2\in[c^{-1},c]$ and $\alpha_n\to \beta$ we
conclude that $\lambda_n\to 0.$ Since $h$ is bounded,  by
$(\ref{lower-semi-conti})$ we get
$
h(\alpha_n)\geq \lambda_nh(\alpha_n^\prime)+(1-\lambda_n)h(\beta)\to h(\beta) \quad ({\rm as}\ n\to\infty),
$
which is in contradiction with the choice of $\alpha_n.$ So $h$ is
lower semi-continuous at $\beta.$

Since each $x\in E$ is a local cone point of $E$ when  $E$ is ${\rm
ri}(A)$, or $E$ is a closed interval in $A$, or $E$ is A itself  and
$A$ is a convex closed polyhedron, the other results follow. \hfill
$\Box$

%%%%%%%%%%%%%%%%%%%%%%%%%%%%%%%%%%%%%%%%%%%%%%%%%%%%%%%%%%%%%%%%%%%%%%%%%%%%%%%%%%%%%%%%%%%%%%%%%%%%%%%%%%

%%%%%%%%%%%%%%%%%%%%%%%%%%%%%%%%%%%%%%%%%%%%%%%%%%%%%%%%%%%%%%%%%%%%%%%%%%%%%%%%%%%%%%%%%%%%%%%%%%%%%%%%%%%%%%%%%%%%%%%%%

\noindent{\bf Proof of Proposition \ref{regularity}.}\ At first we
show that $\Lambda_\Phi^\Psi$ is bounded and
 positive. Fix $\alpha\in L_\Phi.$
  By definition  $\Lambda_\Phi^\Psi(\alpha)\leq
  D(\Psi)$. On the other hand since $\alpha\in L_\Phi$, for any
  $\epsilon>0$, when $n$ large enough, $F(\alpha,n,\epsilon)\ne
  \emptyset.$ Consequently $\Lambda_\Phi^\Psi(\alpha)\geq 0.$
Thus  $\Lambda_\Phi^\Psi(L_\Phi)\subset [0,D(\Psi)]$.

Next we show that $\Lambda_\Phi^\Psi$ is weakly concave. Let
$\alpha,\beta\in L_\Phi$. For any $w_1,\cdots, w_p\in
F(\alpha,n,\epsilon)$ and any $w_{p+1}, \cdots, w_{p+q}\in
F(\beta,n,\epsilon)$, let $w=\overline{w}_1\cdots\overline{w}_{p+q}$
where $\overline{w}_j=w_ju_j$ with $u_j\in \Xi$ such that $w$ is
admissible. By the same argument as for Proposition
\ref{dim-formular-1}, we can show that
$
\exp(-(p+q)(n+c_1(n)))\leq \Psi[w]\leq \exp(-(p+q)(n+\widetilde{c}))
$
with the same $c_1(n)$ and $\widetilde{c}$ as in Proposition
\ref{dim-formular-1}, which means that there exists $u\prec w$ such
that $u\in {\mathcal B}_{(p+q)(n+\widetilde{c})}(\Psi)$. Write
${w}=uw^\prime$. We also have $|w^\prime|\leq (p+q)(ac_1(n)+b)$ with
the same $(a,b)$ as in that proposition.

Now define $F_k(\alpha,n,\epsilon)=\Sigma_{A,k}\cap
F(\alpha,n,\epsilon)$. We have $\displaystyle F(\alpha,n,\epsilon)=\bigcup_{C_1n\leq k\leq
C_2n}F_k(\alpha,n,\epsilon),$
where $C_i=C_i(\Psi)$ for $i=1,2.$ Define
$f_k(\alpha,n,\epsilon)=\#F_k(\alpha,n,\epsilon).$ Choose $k_0$ such
that
$f_{k_0}(\alpha,n,\epsilon)=\max_{C_1n\leq k\leq
C_2n}f_k(\alpha,n,\epsilon).$
Then $f_{k_0}(\alpha,n,\epsilon)\geq
f(\alpha,n,\epsilon)/(C_2-C_1)n.$ Write $k_0=\gamma_n(\alpha)n,$
thus $\gamma_n(\alpha)\in [C_1,C_2]$. Likewise we can find
$\gamma_n(\beta)\in[C_1,C_2]$ such that
$f_{\gamma_n(\beta)n}(\beta,n,\epsilon)\geq
f(\beta,n,\epsilon)/(C_2-C_1)n.$

Fix a subsequence $n_k\uparrow \infty$ such that
$\gamma_{n_k}(\alpha)\to \gamma(\alpha)$ and $\gamma_{n_k}(\beta)\to
\gamma(\beta)$ as $k\to\infty$.  Take $w_1,\cdots,w_p\in
F_{\gamma_{n_k}(\alpha)n_k}(\alpha,n_k,\epsilon)$ and
$w_{p+1},\cdots,w_{p+q}\in
F_{\gamma_{n_k}(\beta)n_k}(\beta,n_k,\epsilon)$. Choose $x_j\in
[w_j]$
 such that
$$\begin{cases} |\phi_{|w_j|}(x_j)-|w_j|\alpha|\leq |w_j|\epsilon, & \text{ if }
1\leq j\leq p\\
|\phi_{|w_j|}(x_j)-|w_j|\beta|\leq |w_j|\epsilon, & \text{ if }
p+1\leq j\leq p+q .\end{cases}$$

  Let
$w=\overline{w}_1\cdots\overline{w}_{p+q}$ and write $w=uw^\prime$
such that $u\in {\mathcal B}_{(p+q)(n_k+\widetilde{c})}(\Psi)$. Then we know
that
$|w|=p(\gamma_{n_k}(\alpha)n_k+p_0)+q(\gamma_{n_k}(\beta)n_k+p_0)$
and $|u|=|w|-|w^\prime|.$ Now for any $x\in [w]$, define $x^1=x$ and
$x^j=T^{\sum_{l=1}^{j-1}|w_l|+p_0}x$ for $j\geq 2.$ Then we have
\begin{eqnarray*}
\phi_{|u|}(x) &=&\phi_{|w|}(x)+\eta_0
=\sum_{j=1}^{p+q}\phi_{|w_j|}(x^j)+\eta_1+\eta_0
=\sum_{j=1}^{p+q}\phi_{|w_j|}(x_j)+\eta_2+\eta_1+\eta_0\\
&=&p\gamma_{n_k}(\alpha)n_k\alpha+q\gamma_{n_k}(\beta)n_k\beta+\eta_3+\eta_2+\eta_1+\eta_0\\
&=&p\gamma(\alpha)n_k\alpha+q\gamma(\beta)n_k\beta+\eta_4+\eta_3+\eta_2+\eta_1+\eta_0,
\end{eqnarray*}
where $|\eta_0|\leq |w^\prime|\|\Phi\|+|C(\Phi)|\leq
(p+q)(ac_1(n_k)+b)\|\Phi\|+|C(\Phi)|,$ $|\eta_1|\leq
(p+q)(p_0\|\Phi\|+2C(\Phi))$; $|\eta_2|\leq
(p+q)\|\Phi\|_{n_k}^\star$; $|\eta_3|\leq
n_k(p\gamma_{n_k}(\alpha)+q\gamma_{n_k}(\beta))\epsilon$, and
$|\eta_4|\leq
pn_k|\alpha||\gamma_{n_k}(\alpha)-\gamma(\alpha)|+qn_k|\beta||\gamma_{n_k}(\beta)-\gamma(\beta)|$.

This yields that for $k$ large enough, $u\in
{F}((p\gamma(\alpha)\alpha+q\gamma(\beta)\beta)/(p\gamma(\alpha)+q\gamma(\beta)),(n_k+\widetilde{c})(p+q),2\epsilon).
$ Thus we conclude that
\begin{eqnarray*}
&&f\Big (\frac{p\gamma(\alpha)\alpha+q\gamma(\beta)\beta}{p\gamma(\alpha)+q\gamma(\beta)},(n_k+\widetilde{c})(p+q),2\epsilon\Big )\\
&\geq&
f_{\gamma_{n_k}(\alpha)n_k}(\alpha,n_k,\epsilon)^pf_{\gamma_{n_k}(\beta)n_k}(\beta,n_k,\epsilon)^qm^{-(p+q)(ac_1(n_k)+b)}\\
&\geq&
f(\alpha,n_k,\epsilon)^pf(\beta,n_k,\epsilon)^q[(C_2-C_1)n_k]^{-p-q}m^{-(p+q)(ac_1(n_k)+b)}.
\end{eqnarray*}
Combining this with Proposition \ref{dim-formular-1} we get
$$\lambda\Lambda_\Phi^\Psi(\alpha)+(1-\lambda)\Lambda_\Phi^\Psi(\beta)\leq
\Lambda_\Phi^\Psi\left (\frac{\lambda\gamma(\alpha)\alpha+(1-\lambda)\gamma(\beta)\beta}{\lambda\gamma(\alpha)+(1-\lambda)\gamma(\beta)}\right)$$
%$$
%\frac{p}{p+q}\Lambda_\Phi^\Psi(\alpha)+\frac{q}{p+q}\Lambda_\Phi^\Psi(\beta)\leq
%\Lambda_\Phi^\Psi\left (\frac{p\gamma(\alpha)\alpha+q\gamma(\beta)\beta}{p\gamma(\alpha)+q\gamma(\beta)}\right ).
%$$
for any $\lambda=\frac{p}{p+q}\in[0,1]\cap\Q$. Since
$\Lambda_\Phi^\Psi$ is upper semi-continuous, we conclude that this
formula holds for any $\lambda\in [0,1]$. Thus $\Lambda_\Phi^\Psi$
is weakly concave.

Assume $A\subset L_\Phi$ is a convex set, and $I\subset L_\Phi$ is a
closed interval. By Lemma \ref{distorsion-concave},
$\Lambda_\Phi^\Psi$ is lower semi-continuous on $\mathrm{ri}(A)$ and
  $I$.
 Combining this  with the upper semi-continuity yields the continuity on $\mathrm{ri}(A)$  and
 $I$. Taking $A=L_\Phi$ we get the continuity on ${\rm ri}(L_\Phi)$.

 Now assume $L_\Phi$ is a polyhedron. By Lemma
 \ref{distorsion-concave}, $\Lambda_\Phi^\Psi$ is lower
 semi-continuous on $L_\Phi$. This,  together with the upper
 semi-continuity yields the continuity on $L_\Phi.$

 Let $I=[\alpha_1,\alpha_2]\subset L_\Phi$ and $\alpha_{\max}\in
 I$ as defined in the proposition. Assume  $\Lambda_\Phi^\Psi$ is not decreasing from $\alpha_{\max}$ to
 $\alpha_1$. Since $\Lambda_\Phi^\Psi$ is continuous on $I$,
 we can find $\beta_1,\beta_2, \beta_3\in [\alpha_1,\alpha_{\max}]$ such that
 $\beta_2\in [\beta_1,\beta_3]$ and $
 \Lambda_\Phi^\Psi(\beta_1)=\Lambda_\Phi^\Psi(\beta_3)>\Lambda_\Phi^\Psi(\beta_2),
 $
which is in contradiction with the fact that $\Lambda_\Phi^\Psi$ is
quasi-concave, since it is weakly concave. Thus
 $\Lambda_\Phi^\Psi$ is decreasing from $\alpha_{\max}$ to
 $\alpha_1$. The same proof  shows that $\Lambda_\Phi^\Psi$ is decreasing from $\alpha_{\max}$ to
 $\alpha_2$.
\hfill$\Box$

\subsection{Proof of Proposition \ref{gene-thermo}} We begin with a lemma
about the Lipschitz continuity of the pressure functions.

\begin{lem}\label{pre-sloope}
Let $\Phi, \Psi\in \C_{aa}(\Sigma_A,T)$ and  define
$f(\lambda)=P(\Phi+\lambda\Psi)$. If $\lambda_1<\lambda_2$ we have
$\Psi_{\min}(\lambda_2-\lambda_1)\leq f(\lambda_2)-f(\lambda_1)\leq \Psi_{\max}(\lambda_2-\lambda_1).$
\end{lem}

\proof\ Let $w\in \Sigma_{A,n}$, and pick up $x_1,x_2\in [w]$ such
that $(\phi_n+\lambda_j\psi_n)(x_j)=\sup_{x\in [w]}(\phi_n+\lambda_j\psi_n)(x).$
Then we have
\begin{eqnarray*}
(\phi_n+\lambda_2\psi_n)(x_2)&\leq&(\phi_n+\lambda_2\psi_n)(x_1)+\|\Phi+\lambda_2\Psi\|_n\\
&=&(\phi_n+\lambda_1\psi_n)(x_1)+(\lambda_2-\lambda_1)\psi_n(x_1)+\|\Phi+\lambda_2\Psi\|_n\\
&\leq&(\phi_n+\lambda_1\psi_n)(x_1)+(\lambda_2-\lambda_1)n\Psi_{\max}+\|\Phi+\lambda_2\Psi\|_n;
\end{eqnarray*}
since $\|\Phi+\lambda_2\Psi\|_n/n\to 0$, this yields $f(\lambda_2)-f(\lambda_1)\leq \Psi_{\max}(\lambda_2-\lambda_1)$.
The other inequality can be proved similarly.
 \hfill $\Box$

\noindent {\bf Proof of Proposition \ref{gene-thermo}.}\ Define
$f(\lambda)=P(\langle z,\Phi-\alpha\rangle+\lambda\Psi). $
 By Lemma \ref{pre-sloope} for $\lambda_1<\lambda_2$ we have
 $
\Psi_{\min}(\lambda_2-\lambda_1)\leq f(\lambda_2)-f(\lambda_1)\leq
\Psi_{\max}(\lambda_2-\lambda_1),
 $
 where $\Psi_{\min}\leq \Psi_{\max}<0.$ Thus $f(\lambda)=0$ has a
 unique solution, which is $\tau_{\Phi}^\Psi(z,\alpha)$.  By Theorem \ref{varia-principle} we have
$$0=P(\langle z,\Phi-\alpha\rangle+\tau_{\Phi}^\Psi(z,\alpha)\Psi)=
\sup\left\{h_\mu+(\langle
z,\Phi-\alpha\rangle+\tau_{\Phi}^\Psi(z,\alpha)\Psi)_\ast(\mu):\mu\in
\M(\Sigma_A,T)\right\}.$$  From this  $(\ref{geometricpressure})$
follows easily.

Now we show $(\ref{pressure-metric})$. By considering the potential
$\Phi^\prime=\Phi-\alpha,$ we can restrict ourselves to the case
$\alpha=0.$ At first we assume $\Phi=(S_n\phi)_{n=1}^{\infty}$ and
$\Psi=(S_n\psi)_{n=1}^{\infty}$ with $\phi$ and $\psi$ H\"{o}lder
continuous. Since $\psi$ is H\"{o}lder continuous, for any $w\in
{\mathcal B}_n(\Psi)$ and any $x\in[w]$ we have
$S_{|w|}\psi(x)\thickapprox -n.$
 Let $\mu$ be the unique equilibrium state of $\langle z,\phi\rangle+\tau(z,0)\psi$. For any $w\in
{\mathcal B}_n(\Psi)$, by the Gibbs property we have
\begin{eqnarray*}
\mu([w])&\thickapprox& \exp\Big(\langle
z,S_{|w|}\phi(x)\rangle+\tau(z,0)S_{|w|}\psi(x)\Big)\
\ (\forall x\in[w])\\
&\thickapprox& \exp\Big(\langle
z,S_{|w|}\phi(x)\rangle-\tau(z,0)n\Big)\
\ (\forall x\in[w])\\
&\thickapprox& \exp\Big(\sup_{x\in[w]}\langle
z,S_{|w|}\phi(x)\rangle-\tau(z,0)n\Big).
\end{eqnarray*}
From this and $\sum_{w\in{\mathcal B}_n(\Psi)}\mu([w])=1$,
$(\ref{pressure-metric})$ follows. The general case requires an approximation argument. We postpone
its proof to the end of Section \ref{appendix}. \hfill$\Box$

\section{Proof of  Theorem~ \ref{main-one-sided}}\label{proof of theo}
Our plan is the following: at first we show that
$\D_\Phi^\Psi(\alpha)\leq \Lambda_\Phi^\Psi(\alpha)\leq {\mathcal
E}_\Phi^\Psi(\alpha)\leq \D_\Phi^\Psi(\alpha),$ then we show
$\Lambda_\Phi^\Psi(\alpha)=\tau_{\Phi}^{\Psi\ast}(\alpha).$ We
divide this into four steps:

\subsection{
$\D_\Phi^\Psi(\alpha)\leq \Lambda_\Phi^\Psi(\alpha)$}\label{upperbd} We prove a
slightly  more general result for the upper bound.  Given
$\Phi\in\C_{aa}(\Sigma_A,T,d)$ and $\Omega\subset L_\Phi$, define
$\displaystyle E_\Phi(\Omega):=\bigcup_{\alpha\in
\Omega}E_\Phi(\alpha). $

\begin{prop}\label{upper-bound} For any compact set $\Omega\subset L_\Phi$
 we have $\dim_P^\Psi E_\Phi(\Omega)\leq \sup\{\Lambda_\Phi^\Psi(\alpha):
\alpha\in\Omega\}$.
In particular, if $\alpha\in L_\Phi$ we have
$\D_\Phi^\Psi(\alpha)\leq \dim_P^\Psi E_\Phi(\alpha)\leq
\Lambda_\Phi^\Psi(\alpha)$.
\end{prop}

%%%%%%%%%%%%%%%%%%%%%%%%%%%%%%%%%%%%%%%%%%%%%%%%%%%%%%%%%%  upper bound  %%%%%%%%%%%%%%%%%%%%%%%%%%%%%%%%%%%%%%%%%%%%%%%%%%%

\proof\
 Let
$\displaystyle
\Lambda_\Phi^\Psi(\alpha,\epsilon):=\limsup_{n\to\infty}\frac{\log
f(\alpha,n,\epsilon,\Phi,\Psi)}{n}, $ then
$\Lambda_{\Phi}^\Psi(\alpha,\epsilon)\searrow
\Lambda_\Phi^\Psi(\alpha) $ when $\epsilon\searrow 0.$ Fix $\eta>0$,
for each $\alpha\in \Omega,$ there exists $\epsilon_\alpha>0$ such
that for any $0<\epsilon\leq \epsilon_\alpha$ we have
$
\Lambda_\Phi^\Psi(\alpha,\epsilon)<\Lambda_\Phi^\Psi(\alpha)+\eta.
$

Since $\{B(\alpha,\epsilon_\alpha): \alpha\in \Omega\}$ is an open
covering of $\Omega$, we can find a finite covering $\{B(\alpha_1,\epsilon_1),\cdots, B(\alpha_s,\epsilon_s)\},$ where
$\epsilon_j=\epsilon_{{\alpha_j}}$. For each $n\in \N$ define
$$H(n,\eta):=
\bigcup_{j=1}^{s}\bigcup_{w\in F(\alpha_j,n,\epsilon)}[w]\ \ \text{
and }\ \ G(k,\eta):=\bigcap_{n\geq k}H(n,\eta)$$

\noindent{\bf Claim:} \ \ $E_\Phi(\Omega)\subset
\bigcup_{k\in\N}G(k,\eta).$

Indeed,  for any $x\in E_\Phi(\Omega)$, there exists
$\alpha\in \Omega$ such that $\phi_n(x)/n\to \alpha$. There exists
$j\in\{1,\cdots,s\}$ such that $\alpha\in B(\alpha_j,\epsilon_j)$.
Take $N$ large enough so that $|\phi_n(x)/n-\alpha|<
\epsilon_j-|\alpha-\alpha_j| $ for any $n\geq N$. For such an $n$ we have
$|\phi_n(x)/n-\alpha_j|\leq |\phi_n(x)/n-\alpha|+|\alpha-\alpha_j|<\epsilon_j,$
hence $x\in
 H(n,\eta)$ for
all $n$ large enough.

 By the previous claim  we have
\begin{equation}\label{sup}
\dim_P^\Psi E_\Phi(\Omega)\leq \sup_{k\in\N} \dim_P^\Psi G(k,\eta).
\end{equation}

Now we find the desired upper bound for the packing dimension of $ G(k,\eta).$ By definition
it is covered by $\{[w]: w\in F(\alpha_j,n,\epsilon); j=1,\cdots,
s\}$ for any $n\geq k$. Since each element in $\{[w]: w\in
F(\alpha_j,n,\epsilon)\}$ is a ball with radius $e^{-n}$, we
conclude that
\begin{eqnarray*}
\dim_{P}^\Psi G(k,\eta) &\leq& \overline{\dim}_B^\Psi G(k,\eta) \leq
\limsup_{n\to\infty}\frac{\log\sum_{j=1}^{s}
f(\alpha_j,n,\epsilon_j)}{n}\\
& \leq& \sup_{j=1,\cdots,s}\limsup_{n\to\infty}\frac{\log
f(\alpha_j,n,\epsilon_j)}{n}
 =
\sup_{j=1,\cdots,s}\Lambda_\Phi^\Psi(\alpha_j,\epsilon_j)\\
& \leq& \sup_{j=1,\cdots,s}\Lambda_\Phi^\Psi(\alpha_j)+\eta
\leq\sup\{\Lambda_\Phi^\Psi(\alpha): \alpha\in \Omega\}+\eta.
\end{eqnarray*}

Combining this with
$(\ref{sup})$ we get $\dim_P^\Psi E_\Phi(\Omega)\leq \sup\{\Lambda_\Phi^\Psi(\alpha):
\alpha\in\Omega\}+\eta.$ Since $\eta$ is arbitrary, we get the
result.
 \hfill $\Box$

%%%%%%%%%%%%%%%%%%%%%%%%%%%%%%%%%%%%%%%%%%%%%%%%%%%%%%%%%%%%%%%%%%%%%%%%%%%%%%%

\subsection{
$\Lambda_\Phi^\Psi(\alpha)\leq {\mathcal E}_\Phi^\Psi(\alpha)$}  Our
approach is inspired by that of \cite{FF}, which deals with the
case that $\Psi$ is additive and built from a constant negative
potential.

 To show this
inequality we need to approximate the almost additive potentials
$\Phi$ and $\Psi$ by two sequences of H\"{o}lder potentials. We
describe this procedure as follows.

 Given $\Phi\in\C_{aa}(\Sigma_A,T,d)$, for each $k\in\N$ we define
${\Phi}^k\in \C_{aa}(\Sigma_A,T,d)$  as follows. For each
$w\in\Sigma_{A,k}$ choose $x_w\in [w]$. For any $x\in[w]$ define $
\widetilde{\phi}_k(x):={\phi_k(x_w)}/{k}. $ Define
\begin{equation}\label{app-potential}
{\phi}_n^k:=S_n\widetilde{\phi}_k\ \ \text{ and }\ \
{\Phi}^k=({\phi}_n^k)_{n=1}^{\infty}.
\end{equation}
  Thus ${\Phi}^k$ is additive and
$\widetilde{\phi}_k$ depends only on the first $k$ coordinates of
$x\in\Sigma_A.$ Consequently $\Phi^k$ is H\"{o}lder continuous.

\begin{lem}\label{control}
  We have $\Phi_{\min}\leq
\Phi_{\min}^k\leq\Phi_{\max}^k\leq \Phi_{\max}$. Moreover
$
\|\phi_n-{\phi}_n^k\|\leq\frac{n}{k}|C(\Phi)|+4k\|\Phi\|+\frac{\sqrt{d}\|\Phi\|_k}{k}n.
$
Consequently $\|\Phi-\Phi^k\|_{\mbox{\tiny \rm lim}}\to 0$ when $k\to\infty.$
\end{lem}
This lemma will be proved at the end of this  subsection.

\noindent{\bf Proof of $\Lambda_\Phi^\Psi(\alpha)\leq {\mathcal
E}_\Phi^\Psi(\alpha)$.}\ Now, for ${\Phi}\in \C_{aa}(\Sigma_A,T,d)$
and $\Psi\in \C_{aa}^-(\Sigma_A,T)$ define $\Phi^k $ and $\Psi^k$
according to $(\ref{app-potential})$. Fix $\epsilon>0$. By Lemma
\ref{control}, we can find $K(\epsilon)$ such that   for each
 $k\geq K(\epsilon)$ and sufficiently large $n$(related to $k$) we
 have $\|\phi_n-{\phi}_n^k\|_\infty\leq n\epsilon/2$ and $\|\psi_n-{\psi}_n^k\|_\infty\leq n\epsilon/2$. Then
$F(\alpha,n,\epsilon/2,\Phi,\Psi)\subset F(\alpha,n,\epsilon,\Phi^k,\Psi),$ and
consequently $f(\alpha,n,\epsilon/2,\Phi,\Psi)\leq
f(\alpha,n,\epsilon,\Phi^k,\Psi)$.

For any word $w$ such that $|w|\geq k$, we define the integer valued
function $\theta_w:\Sigma_{A,k}\to \N$ as $\theta_w(u)=\#\{j: w_j\cdots w_{j+k-1}=u\}.$
It is clear that
\begin{equation}\label{ltheta}
\sum_u\theta_w(u)=|w|-k+1.
\end{equation}
Let
${\mathcal P}_k^{(n)}=\{\theta_w: w=w_1w_2 \ \text{admissible}, \ w_1\in F(\alpha,n,\epsilon,\Phi^k,\Psi), |w_2|=k-1\}.$
Since $w_1\in \mathcal B_n(\Psi)$ we have  $|w_1|\leq C_2(\Psi)n,$ thus
$\#{\mathcal P}_k^{(n)}\leq (C_2n)^{m^k}$. For each $\theta\in
{\mathcal P}_k^{(n)}$, let ${\mathcal T}(\theta)$ be the collection
of all $w_1w_2$ such that $w_1\in F(\alpha,n,\epsilon,\Phi^k,\Psi),
|w_2|=k-1$ and $\theta_{w_1w_2}=\theta.$ Then we have
$$f(\alpha,n,\epsilon,\Phi^k,\Psi)\leq \sum_{\theta\in{\mathcal P}_k^{(n)}}\#{\mathcal T}(\theta)
\leq (C_2n)^{m^k}\max_{\theta\in{\mathcal P}_k^{(n)}}\#{\mathcal
T}(\theta).$$
Thus $\displaystyle \frac{\log f(\alpha,n,\epsilon/2,\Phi,\Psi)}{n}\leq \frac{\log f(\alpha,n,\epsilon,\Phi^k,\Psi)}{n}\leq
\max_{\theta\in{\mathcal P}_k^{(n)}}\frac{\log \#{\mathcal
T}(\theta)}{n}+m^kO(\frac{\log n}{n}).$

Following ~\cite{FF} we define $\triangle_k^+$, the set of all
positive functions $p$ on $\Sigma_{A,k}$ satisfying the following
two relations:
$$
\sum_{w\in\Sigma_{A,k}}p(w)=1;\ \  \sum_{w}p(ww_1w_2\cdots
w_{k-1})=\sum_{w}p(w_1w_2\cdots w_{k-1}w).
$$

It is  known (see \cite{FF}) that for any $\eta>0$, there is a positive integer
 $N=N(\eta)$ such that for any $w\in \Sigma_{A,l+k-1}$ with
$l>N$, there exists a probability vector $p\in \triangle_k^+$ such
that
$$\Big|\frac{\theta_w(u)}{l}-p(u)\Big |<\eta, \ \ p(u)>\frac{\eta}{m^{k+1}}.$$
We discard the trivial case where $\Phi\equiv 0$ and fix  $\eta>0$ such that
$\eta<\epsilon/(m^k\|\Phi\|).$

Now we fix $\theta\in {\mathcal P}_k^{(n)}$ with $n$ large enough,
and write $\theta=\theta_{ww^\prime}$ with $w\in
F(\alpha,n,\epsilon,\Phi^k,\Psi)$ and $|w^\prime|=k-1$. Notice that
 any word $v\in {\mathcal T}(\theta)$ can be written as $v_1v_2$ with $|v_2|=k-1$ and
$|v_1|$ also equal to a constant (this is due to \eqref{ltheta}) that we denote by $l_\theta$.  Fix a
$p\in\triangle_k^+$ as described above. Consider the Markov measure
$\nu_p$ corresponding to $p$ (see ~\cite{FF} for the definition and
related properties). For any word $v=v_1v_2\in{\mathcal T}(\theta)$
with $|v_1|=l_\theta$ (and $v_1\in F(\alpha,n,\epsilon,\Phi^k,\Psi)$) and $|v_2|=k-1$, we
have
\begin{eqnarray*}
\nu_p([v_1v_2])=p(v_1v_2|_k)\prod_{|u|=k} t(u)^{\theta(u)} \geq
\frac{\eta}{m^{k+1}}\prod_{|u|=k} t(u)^{\theta(u)}:=\rho,
\end{eqnarray*}
where $\displaystyle t(a_1\cdots a_k)=\frac{p(a_1\cdots a_k)}{\sum_{\epsilon}p(a_1\cdots a_{k-1}\epsilon)}$.
 Also $\rho\#{\mathcal T}(\theta)\leq \nu_p(\bigcup_{v\in{\mathcal T}(\theta)}[v])\leq 1$. Thus,
$$\#{\mathcal T}(\theta)\leq \frac{1}{\rho}=\frac{m^{k+1}}{\eta}\prod_{|u|=k}t(u)^{-\theta(u)}.$$
Since $C_1(\Psi)n\leq l_\theta\leq C_2(\Psi)n$ and
$\eta/m^{k+1}\leq t(u)\leq 1$, we have
\begin{eqnarray*}
\frac{\log\#{\mathcal T}(\theta)}{l_\theta}&\leq& O(\frac{k}{n})+
O(\frac{|\log\eta|}{n})-\sum_{|u|=k}\frac{\theta(u)}{l_\theta}\log
t(u)\\
&\leq&O(\frac{k}{n})+ O(\frac{|\log\eta|}{n})-\sum_{|u|=k}p(u)\log
t(u)+m^k\eta(|\log \eta|+(k+1)\log m)\\
&=& h(\nu_p)+ O(\frac{k}{n})+
O(\frac{|\log\eta|}{n})+m^k\eta(|\log\eta|+(k+1)\log m).
\end{eqnarray*}

Next we estimate $n/l_\theta.$ Let
$x_0\in[ww^\prime]$. By $(\ref{444})$ we have
$$
 -n-C(\Psi)-2\|\Psi\|_{l_\theta}+\Psi_{\min}\leq
\psi_{l_\theta}(x_0)\leq\sup_{x\in[w]}\psi_{l_\theta}(x)\leq -n,
$$
and by Lemma \ref{control}, when $k$ and $n$ are large enough we
have $\|\psi_{l_\theta}-{\psi}^k_{l_\theta}\|_\infty\leq
l_\theta\epsilon$, thus
\begin{equation}\label{comparible}
 -n-C(\Psi)-2\|\Psi\|_{n}^\star+\Psi_{\min}-C_2n\epsilon\leq
{\psi}_{l_\theta}^k(x_0)\leq -n+C_2n\epsilon.
\end{equation}
Also,
\begin{eqnarray*}
\frac{{\psi}_{l_\theta}^k(x_0)}{l_\theta}&=&\sum_{|u|=k}\frac{\theta(u)}{l_\theta}\widetilde{\psi}_k(x_u)
=\sum_{|u|=k}p(u)\widetilde{\psi}_k(x_u)+m^kO(\eta)\\
&=&\int\widetilde{\psi}_k d\nu_p+m^kO(\eta)
={\Psi}^k_\ast(\nu_p)+m^kO(\eta)
=\Psi_\ast(\nu_p)+O(\epsilon)+m^kO(\eta).
\end{eqnarray*}
Combining this with $(\ref{comparible})$ and the fact that
$\|\Psi\|_n^\star/n=o(1)$ we get
$$\frac{n}{l_\theta}=-\Psi_\ast(\nu_p)+O(\epsilon)+m^kO(\eta)+o(1).$$
As a result we get
$$
\frac{\log\#{\mathcal T}(\theta)}{n}=\frac{\log\#{\mathcal
T}(\theta)}{l_\theta}\cdot\frac{l_\theta}{n} \leq\frac{h(\nu_p)+
O(\frac{k}{n})+ O(\frac{|\log\eta|}{n})+m^k\eta(|\log\eta|+(k+1)\log
m)}{-\Psi_\ast(\nu_p)+O(\epsilon)+m^kO(\eta)+o(1)}.
$$

 Since $w\in F(\alpha,n,\epsilon,\Phi^k,\Psi)$, there exists $y_0\in[w]$ such that
 $|{\phi}_{l_\theta}^k(y_0)/l_\theta-\alpha|\leq \epsilon.$ We  have
 \begin{eqnarray*}
 &&|\Phi_\ast(\nu_p)-\alpha|\\
 &\leq&  |{\Phi}^k_\ast(\nu_p)-\alpha|+|{\Phi}_\ast(\nu_p)-{\Phi}^k_\ast(\nu_p)|
 \leq  |\int\widetilde{\phi}_k
 d\nu_p-\alpha|+\epsilon
 =  |\sum_{|u|=k}p(u)\widetilde{\phi}_k(x_u)-\alpha|+\epsilon\\
 &\leq&  |\sum_{|u|=k}\frac{\theta_{ww^\prime}(u)}{l_\theta}\widetilde{\phi}_k(x_u)-\alpha|+m^k\eta\|\widetilde{\phi}_k\|+\epsilon \ \
 (\theta_{ww^\prime}=\theta,w\in F(\alpha,n,\epsilon,\Phi^k,\Psi))\\
 &\leq&  |{\phi}_{l_\theta}^k(x)/l_\theta-\alpha|+m^k\eta\|\Phi\|+\epsilon \ \ (\text{ for any } x\in[ww^\prime])\\
 &\leq&|{\phi}_{l_\theta}^k(x)/l_\theta-{\phi}_{l_\theta}^k(y_0)/l_\theta|+
 |{\phi}_{l_\theta}^k(y_0)/l_\theta-\alpha|+m^k\eta\|\Phi\|+\epsilon \leq \frac{k\|\Phi\|}{C_1n} +m^k\eta\|\Phi\|+2\epsilon.
\end{eqnarray*}

By our choice of $\eta$ we have $m^k\eta\|\Phi\|<\epsilon$. Moreover,
when $n\geq k\|\Phi\|/(C_1\epsilon)$ we have $k\|\Phi\|/(C_1n)\leq
\epsilon.$ Letting $n\to\infty$ and then $\eta\to 0$ we conclude that
$$
\limsup_{n\to\infty}\frac{\log
f(\alpha,n,\epsilon/2)}{n}\leq\sup_{|\Phi_\ast(\nu)-\alpha|\leq
4\epsilon}\frac{h(\nu)}{-\Psi_\ast(\nu)+O(\epsilon)}.
$$
Notice that the set of invariant measures $\nu$ such that
$|\Phi_\ast(\nu)-\alpha|\leq 4\epsilon$ is compact, so by using the
upper semi-continuity of $h(\nu)$  and letting $\epsilon$ tend to 0 we can
find an invariant measure $\nu_0$ such that $\Phi_\ast(\nu_0)=\alpha$ and
$$\Lambda_\Phi^\Psi(\alpha)=\lim_{\epsilon\to 0} \limsup_{n\to\infty}\frac{\log f(\alpha,n,\epsilon/2)}{n}
\leq \frac{h(\nu_0)}{-\Psi_\ast(\nu_0)}\leq {\mathcal
E}_\Phi^\Psi(\alpha).$$

\noindent
{\bf Proof of Lemma~\ref{control}.} At first we assume $\Phi\in\C_{aa}(\Sigma_A,T)$. By
$(\ref{max-min})$ we get $\Phi_{\min}\leq\widetilde{\phi}_k\leq
\Phi_{\max}.$ Since $\Phi^k$ is additive, we have
$
\Phi_{\min}\leq\widetilde{\phi}_{k\min}=
\Phi_{\min}^k\leq\Phi_{\max}^k=\widetilde{\phi}_{k\max}\leq \Phi_{\max}.
$

For $n\in \N$, write $n=pk+l$ with $0\leq l<k.$ For $1\leq j\leq k$
we have
\begin{eqnarray*}
\phi_n(x)&\leq&
\phi_j(x)+\sum_{l=0}^{p-2}\phi_k(T^{j+lk}x)+pC(\Phi)+\phi_{k+l-j}(T^{j+(p-1)k})\\
&\leq& \sum_{l=0}^{p-2}\phi_k(T^{j+lk}x)+pC(\Phi)+2k\|\Phi\|.
\end{eqnarray*}
So we get
\begin{eqnarray*}
\phi_n(x) &\leq&
\sum_{j=1}^k\sum_{l=0}^{p-2}\phi_k(T^{j+lk}x)/k+pC(\Phi)+2k\|\Phi\|\\
&=&\sum_{j=1}^{(p-1)k}\phi_k(T^jx)/k+pC(\Phi)+2k\|\Phi\|
\leq\sum_{j=0}^{n-1}\phi_k(T^jx)/k+pC(\Phi)+4k\|\Phi\|\\
&\leq&\sum_{j=0}^{n-1}\widetilde{\phi}_k(T^jx)+pC(\Phi)+4k\|\Phi\|+\frac{\|\Phi\|_k}{k}n
={\phi}_n^k(x)+pC(\Phi)+4k\|\Phi\|+\frac{\|\Phi\|_k}{k}n.
\end{eqnarray*}
Similarly we have $ \phi_n(x) \geq
{\phi}_n^k(x)-pC(\Phi)-4k\|\Phi\|-\frac{\|\Phi\|_k}{k}n,$ hence $
\|\phi_n-\widetilde{\phi}_n^k\|_\infty\leq
pC(\Phi)+4k\|\Phi\|+\frac{\|\Phi\|_k}{k}n,
$ and $
\|\Phi-\Phi^k\|_{\mbox{\tiny \rm
lim}}\leq
\frac{C(\Phi)}{k}+\frac{\|\Phi\|_k}{k}\to 0\quad  (\text{as $k\to\infty$}) .
$

If $\Phi=(\Phi^{(1)},\cdots,\Phi^{(d)})\in\C_{aa}(\Sigma_A,T,d)$, applying the result just proven to each component of $\Phi$ and using $(\ref{component})$ we get the result. \hfill $\Box$

 \subsection{ ${\mathcal E}_\Phi^\Psi(\alpha)\leq
 \D_\Phi^\Psi(\alpha)$}

 It is the content of Proposition \ref{lower-unify}.
  
  %Nevertheless the reader can find a specific proof in \cite{BQ}.

Until now we have shown that $\D_\Phi^\Psi(\alpha)={\mathcal
E}_\Phi^\Psi(\alpha)=\Lambda_\Phi^\Psi(\alpha).$

%%%%%%%%%%%%%%%%%%%%%%%%%%%%%%%%%%%%%%%%%%%%%%%%%%%%%%%%%%%%%%%%%%%%%%%%%%%%%%%%%%%%%%%%%%%%%%%%%%%%%%%%%%%%%%%%%%%%%%%%%%%%%%%%

\subsection{
The large deviation principle $\Lambda_\Phi^\Psi(\alpha)=\tau_{\Phi}^{\Psi\star}(\alpha)$,  and
$\dim_H^\Psi(\Sigma_A)=\dim_B^\Psi(\Sigma_A)=D(\Psi)=\max\{\Lambda_\Phi^{\Psi}(\alpha):
\alpha\in L_\Phi\}$} $\ $

At first we show the following simple fact:

\begin{lem} \label{easy-fact}
For any $\alpha\in L_\Phi$ we have
$\Lambda_\Phi^\Psi(\alpha)\leq\tau_{\Phi}^{\Psi\star}(\alpha)$.
\end{lem}

\proof\  Fix $\alpha\in L_\Phi$. By the variational principle, for
any $\mu\in\M(\Sigma_A,T)$
$$
0= P(\langle
z,\Phi-\alpha\rangle+\tau_{\Phi}^{\Psi}(z,\alpha)\Psi)\geq h_\mu(T)+
\langle
z,\Phi_\ast(\mu)-\alpha\rangle+\tau_{\Phi}^{\Psi}(z,\alpha)\Psi_\ast(\mu).
$$
Thus if $\Phi_\ast(\mu)=\alpha$, we get $
\tau_{\Phi}^{\Psi}(z,\alpha)\geq -{h_\mu(T)}/{\Psi_\ast(\mu)}. $
This implies that
$\tau_{\Phi}^{\Psi\star}(\alpha)=\inf\{\tau_{\Phi}^{\Psi}(z,\alpha):z\in\R^d\}\geq
{\mathcal E}_\Phi^\Psi(\alpha).$ Since ${\mathcal
E}_\Phi^\Psi(\alpha)=\Lambda_\Phi^\Psi(\alpha)$,  the result
follows. \hfill$\Box$

\bigskip

Next we show
$\Lambda_\Phi^\Psi(\alpha)=\tau_{\Phi}^{\Psi\star}(\alpha).$ We do
this  at first for H\"{o}lder potentials, then we deal with the
general case by using an approximation procedure.
%%%%%%%
\subsubsection{ {\bf $\Lambda_\Phi^\Psi=\tau_{\Phi}^{\Psi\star}$
 when $\Phi$ and $\Psi$ are H\"older potentials and $L_\Phi$ has dimension~$d$}}
\begin{lem}\label{case-holder}
Assume $\Phi$ and $\Psi$ are H\"{o}lder continuous potentials and
$L_\Phi$ has dimension $d$. Then
\begin{enumerate}
\item $\tau_{\Phi}^{\Psi}(0,\alpha)=D(\Psi)$, consequently
$\tau_{\Phi}^{\Psi\star}(\alpha)\leq D(\Psi)$ for any $\alpha\in
\R^d$. Moreover $D(\Psi)$ is the unique root of $P(\lambda\Psi)=0.$

\item Let $(z,\alpha,\alpha')\in(\R^d)^3$. If $\langle z,\alpha^\prime-\alpha\rangle\geq 0$, then
\begin{equation}\label{sloope-control}
C_1(\Psi)\langle z,\alpha^\prime-\alpha\rangle\leq
\tau_{\Phi}^{\Psi}(z,\alpha)-\tau_{\Phi}^{\Psi}(z,\alpha^\prime)\leq
C_2(\Psi)\langle z,\alpha^\prime-\alpha\rangle.
\end{equation}

\item  Let $\alpha\in\R^d$. Then $\tau_{\Phi}^{\Psi}(\cdot, \alpha)$ is convex.

\item If $\alpha\in \mathrm{int}(L_\Phi)$ and
 $\delta_0>0$ is such that $B(\alpha,\delta_0)\subset L_\Phi$, then for
 any $z\in\R^d$, $
\tau_{\Phi}^{\Psi}(z,\alpha)\geq \delta_0C_1|z|/2
$, where $C_1=C_1(\Psi)=1/|\Psi_{\min}|$. \end{enumerate}
\end{lem}

\proof\ (1) By $(\ref{full-dim})$ and $(\ref{pressure-metric})$ we
get $\tau_{\Phi}^{\Psi}(0,\alpha)=D(\Psi).$ By definition
$D(\Psi)={\tau_\Phi^\Psi}(0,\alpha)$ is the unique root of
$P(\lambda\Psi)=0.$

\noindent
(2) By Proposition \ref{gene-thermo}, we have
$$
{\tau_\Phi^\Psi}(z,\alpha)=\lim_{n\to\infty}\frac{1}{n}\log
\sum_{w\in {\mathcal B}_n(\Psi)} \exp(\sup_{x\in[w]}\ \langle
z,\phi_{|w|}(x)-|w|\alpha\rangle).
$$
For a fixed $w\in\mathcal B_n(\Psi)$ and any $x,y\in[w]$, since
$\Phi$ is H\"{o}lder continuous and $\langle
z,\alpha^\prime-\alpha\rangle\geq 0$, we get
\begin{eqnarray*}
\langle z,\phi_{|w|}(x)-|w|\alpha\rangle&=&\langle
z,\phi_{|w|}(x)-|w|\alpha^\prime\rangle+|w|\langle
z,\alpha^\prime-\alpha\rangle\\
&\leq&\langle
z,\phi_{|w|}(y)-|w|\alpha^\prime\rangle+C_2(\Psi)n\langle
z,\alpha^\prime-\alpha\rangle+C(\Phi,z).
\end{eqnarray*}
Thus $
\max_{x\in[w]}\langle
z,\phi_{|w|}(x)-|w|\alpha\rangle\leq\max_{y\in[w]}\langle
z,\phi_{|w|}(y)-|w|\alpha^\prime\rangle+C_2(\Psi)n\langle
z,\alpha^\prime-\alpha\rangle+C(\Phi,z).
$
From this we get $
{\tau_\Phi^\Psi}(z,\alpha)-{\tau_\Phi^\Psi}(z,\alpha^\prime)\leq
C_2(\Psi)\langle z,\alpha^\prime-\alpha\rangle$. The other
inequality follows similarly.

\noindent
(3) It is a classical result, see for example Lemma 6 in  \cite{PW}.

\noindent (4) Assume $\alpha\in \mathrm{int}(L_\Phi)$ and
 $\delta_0>0$ such that $B(\alpha,\delta_0)\subset L_\Phi$. For any $z\in \R^d$,  let
$\alpha^\prime=\alpha+\delta_0z/|z|$. We  have
\begin{eqnarray*}
&&\sum_{w\in {\mathcal B}_n(\Psi)}\exp(\sup_{x\in[w]}\ \langle
z,(\phi-\alpha)_{|w|}(x)\rangle)
\geq\sum_{w\in
F(\alpha^\prime,n,\delta_0/2,\Phi,\Psi)}\exp(\sup_{x\in[w]}\ \langle
z,(\phi-\alpha)_{|w|}(x)\rangle)\\
&\geq&\sum_{w\in
F(\alpha^\prime,n,\delta_0/2,\Phi,\Psi)}C\exp(|w||z|\delta_0/2) \geq
C\exp(C_1n|z|\delta_0/2)f(\alpha^\prime,n,\delta_0/2,\Phi,\Psi).
\end{eqnarray*}
This yields
$\displaystyle {\tau_\Phi^\Psi}(z,\alpha)\geq \delta_0C_1|z|/2+\limsup_{n\to\infty}
\frac{\log
f(\alpha^\prime,n,\delta_0/2,\Phi,\Psi)}{n}\geq\delta_0C_1|z|/2.$
\hfill $\Box$

\begin{lem}\label{holder-basic}
 If $\Phi$ and $\Psi$ are H\"{o}lder continuous and $L_\Phi$ has dimension $d$, then
$
\Lambda_\Phi^\Psi(\alpha)={\tau_\Phi^\Psi}^\star(\alpha)
$ over $L_\Phi$.
\end{lem}

\proof\ When $\alpha\in \mathrm{int}(L_\Phi)$, the result has been
shown in ~\cite{BSS}. Specifically, there exists $z\in\R^d$ such
that ${\tau_\Phi^\Psi}^\star(\alpha)={\tau_\Phi^\Psi}(z,\alpha)$.

 Now we assume $\alpha$ is in the  boundary of
$L_\Phi$. When $\Phi$ and $\Psi$ are H\"{o}lder continuous, it is
known that $\max\{\Lambda_\Phi^\Psi(\alpha):\alpha\in
L_\Phi\}=\dim_H^\Psi\Sigma_A=\dim_B^\Psi\Sigma_A$ and
$\dim_H^\Psi\Sigma_A$ is the unique root of $P(\lambda\Psi)=0$
\cite{FLW}. Thus by Lemma \ref{case-holder} (1),
$\dim_H^\Psi\Sigma_A=D(\Psi)$. If
$\Lambda_\Phi^\Psi(\alpha)=D(\Psi)$, then by Lemma
\ref{case-holder}(1) we have ${\tau_\Phi^\Psi}^\star(\alpha)\leq
\Lambda_\Phi^\Psi(\alpha).$ So in the following we assume
$\Lambda_\Phi^\Psi(\alpha)<D(\Psi)$. By the regularity of
$\Lambda_\Phi^\Psi$ we can find $\alpha_\ast\in
\mathrm{int}(L_\Phi)$ such that
$\Lambda_\Phi^\Psi(\alpha_\ast)>\Lambda_\Phi^\Psi(\alpha).$ The line
passing through $\alpha$ and $\alpha_\ast$ intersects the boundary
of $L_\Phi$ at another point $\alpha_1.$ Let
$\alpha_{\max}\in[\alpha,\alpha_1]$ such that
$\Lambda_\Phi^\Psi(\alpha_{\max})=\max\{\Lambda_\Phi^\Psi(\beta):\beta\in
[\alpha,\alpha_1]\}$. By Proposition \ref{regularity},
$\Lambda_\Phi^\Psi$ is non-increasing from $\alpha_{\max}$ to
$\alpha$. Let $\beta_0$ be a point in the open interval
$(\alpha,\alpha_{\max})$ such that
$\Lambda_\Phi^\Psi(\beta_0)>\Lambda_\Phi^\Psi(\alpha)$. We have
$\beta_0\in \mathrm{int}(L_\Phi)$. Let
$\alpha_t=t\beta_0+(1-t)\alpha\ (0<t<1)$. Then $\alpha_t\in
\mathrm{int}(L_\Phi)$ and $\alpha_t\to\alpha$ as $t\to 0.$ Let
$z_t\in\R^d$ such that
$\Lambda_\Phi^\Psi(\alpha_t)={\tau_\Phi^\Psi}^\star(\alpha_t)={\tau_\Phi^\Psi}(z_t,\alpha_t)$.

We claim that $\langle z_t,\beta_0-\alpha_{t}\rangle\leq 0$.
Otherwise $\langle z_t,\beta_0-\alpha_{t}\rangle >0$, and by
$(\ref{sloope-control})$ we have
${\tau_\Phi^\Psi}(z_t,\beta_0)\leq {\tau_\Phi^\Psi}(z_t,\alpha_t)-C_1(\Psi)
\langle z_t,\beta_0-\alpha_t\rangle<{\tau_\Phi^\Psi}(z_t,\alpha_t)=\Lambda_\Phi^\Psi(\alpha_t).$
On the other hand we should have ${\tau_\Phi^\Psi}(z_t,\beta_0)\geq
{\tau_\Phi^\Psi}^\star(\beta_0)=\Lambda_\Phi^\Psi(\beta_0)\geq
\Lambda_\Phi^\Psi(\alpha_t),$ which is a contradiction.

Consequently, due to the definition of $\beta_0$  we have $\langle
z_t,\alpha-\alpha_{t}\rangle\geq 0$. Still by
$(\ref{sloope-control})$ we get
${\tau_\Phi^\Psi}(z_t,\alpha)\leq {\tau_\Phi^\Psi}(z_t,\alpha_t)-C_1(\Psi)\langle z_t,\alpha-\alpha_t\rangle
\leq{\tau_\Phi^\Psi}(z_t,\alpha_t)=\Lambda_\Phi^\Psi(\alpha_t).$
Thus
${\tau_\Phi^\Psi}^\star(\alpha)\leq\Lambda_\Phi^\Psi(\alpha_t)$.
Letting $t$ tend to 0, by the continuity of $\Lambda_\Phi^\Psi$ on
the closed interval $[\alpha,\alpha_1]$, we get
${\tau_\Phi^\Psi}^\star(\alpha)\leq \Lambda_\Phi^\Psi(\alpha).$
Combining this with Lemma \ref{easy-fact}, we get the equality.
 \hfill
$\Box$

\subsubsection{ {\bf
$\dim_H^\Psi(\Sigma_A)=\dim_B^\Psi(\Sigma_A)=D(\Psi)=\max\{\Lambda_\Phi^{\Psi}(\alpha): \alpha\in L_\Phi\}$ in the general case}}$\ $

 We need to describe the $\Psi$- and $\Phi$-
dependence of the function $\Lambda_\Phi^\Psi.$  Recall that in
Lemma \ref{upper-bound} we have
 defined
 $\Lambda_{\Phi}^{\Psi}(\alpha,\epsilon)=\limsup_{n\to\infty}{\log( f(\alpha,n,\epsilon,\Phi,\Psi))}/{n}$
 and we know that
 $\Lambda_{\Phi}^{\Psi}(\alpha,\epsilon)\searrow\Lambda_{\Phi}^{\Psi}(\alpha)$ as $\epsilon \searrow 0$.
The following lemma will be proved in Section~\ref{appendix}.
 \begin{lem}\label{Psi-dependence}
  \begin{enumerate}
  \item Assume $\Psi,\Upsilon\in \C_{aa}^{-}(\Sigma_A,T)$, then we have
 \begin{equation}\label{approximation-1}
 |D(\Psi)-D(\Upsilon)|\leq 2\log m\cdot\Big (1+\frac{1}{|\Psi_{\max}|}\Big )\Big (1+\frac{1}{|\Upsilon_{\max}|}\Big )\|\Psi-\Upsilon\|_{\mbox{\tiny \rm lim}}.
 \end{equation}

\item Let  $\Phi,\Theta\in \C_{aa}(\Sigma_A,T,d)$. Let $\beta\in
L_\Theta.$ For any $\alpha\in B(\beta,\eta)\cap L_\Phi$ and
$\delta_0<1/2C_2(\Psi)$ we have
 \begin{equation}\label{approximation-3}
 \Lambda_{\Phi}^{\Psi}(\alpha)\leq \frac{2C_2(\Psi)\log
m}{|\Upsilon_{\max}|}\delta_0+(1-C_2(\Psi)\delta_0)\Lambda_{\Theta}^{\Upsilon}(\beta,
a_0+\kappa\delta_0+2\eta),
 \end{equation}
 where $\delta_0=\|\Psi-\Upsilon\|_{\mbox{\tiny \rm lim}}$,
$a_0=\|\Phi-\Theta\|_{\mbox{\tiny \rm lim}}$, $C_2(\Psi)=1+1/|\Psi_{\max}|$ and
$\kappa=\kappa(\Psi,\Upsilon,\Phi)=14\|\Phi\|C_2(\Psi){|\Upsilon_{\min}|}/
{|\Upsilon_{\max}|}.$
\end{enumerate}
\end{lem}

 %%%%%%%%%%%%%%%%%%%%%%%%%%%%%%%%%%%%%%%%%%%%%%%%%%%%%%%%%%%%%%%%%%%%%%%%%%%%%%%%%%%%%%55

\noindent {\bf Proof of the fact that
$\dim_H^\Psi(\Sigma_A)=\dim_B^\Psi(\Sigma_A)=D(\Psi)=\max\{\Lambda_\Phi^{\Psi}(\alpha):
\alpha\in L_\Phi\}$.} At first assume that $\Psi$ is H\"{o}lder
continuous. By Lemma \ref{case-holder} (1), $P(D(\Psi)\Psi)=0$. Let
$\mu$ be the unique equilibrium state of $D(\Psi)\Psi$. It is well
known (see \cite{Bow75}) that $\mu$ is ergodic, and
$D(\Psi)=\dim_H^\Psi \Sigma_A=\dim_H^\Psi\mu $ (\cite{Bow79}). Let
$\alpha=\Phi_\ast(\mu)$. By the sub-additive ergodic theorem we have
$\mu(E_\Phi(\alpha))=1$, consequently
$\Lambda_{\Phi}^{\Psi}(\alpha)=\D_\Phi^\Psi(\alpha)=D(\Psi).$ Thus,
when $\Psi$ is a H\"{o}lder potential the result holds.

Next we  assume $\Psi\in\C_{aa}^{-}(\Sigma_A,T)$. Define $\Psi^n$
according to $(\ref{app-potential}),$ then
$\lim_{n\to\infty}\|\Psi-\Psi^n\|_{\mbox{\tiny \rm lim}}= 0$ and
$|\Psi_{\max}|\leq|\Psi_{\max}^n|\leq|\Psi_{\min}|.$ By
$(\ref{approximation-1})$ we have $\lim_{n\to\infty} D(\Psi^n)=
D(\Psi)$. Let $\mu_n$ be the unique equilibrium state of
$D(\Psi^n)\cdot\Psi^n$ and define $\alpha_n=\Phi_\ast(\mu_n)$. Then
$\alpha_n\in L_\Phi$ and
$\Lambda_{\Phi}^{\Psi^n}(\alpha_n)=D(\Psi^n).$ Let $\alpha$ be a
limit point of the sequence $\{\alpha_n: n\in \N\}$. Without loss of
generality we assume $\alpha=\lim_{n\to\infty}\alpha_{n}$. By
$(\ref{approximation-3})$ we have
\begin{equation}\label{555}
\Lambda_{\Phi}^{\Psi^n}(\alpha_n)\leq
 \frac{2C_2(\Psi^n)\log
m}{|\Psi_{\max}|}\delta_n+(1-C_2(\Psi^n)\delta_n)
 \Lambda_{\Phi}^{\Psi}(\alpha,
 \kappa_n\delta_n+2\eta_n),
 \end{equation}
where
$
\delta_n:=\|\Psi-\Psi^{n}\|_{\mbox{\tiny \rm lim}},\
 C_2(\Psi^{n})=1+\frac{1}{|\Psi_{\max}^{n}|},\
\kappa_n=14\|\Phi\|C_2(\Psi^{n})\frac{{|\Psi_{\min}|}}{
{|\Psi_{\max}|}},\ \eta_n=|\alpha-\alpha_n|. $ By Lemma
\ref{control} we have  $C_2(\Psi^{n})\leq 1+1/|\Psi_{\max}|, $ thus
we can rewrite $(\ref{555})$ as $ D(\Psi^{n})\leq
d_1\delta_n+\Lambda_{\Phi}^{\Psi}(\alpha, d_2\delta_n+2\eta_n). $
Letting $n$ tend to $\infty$ we get $D(\Psi)\leq
\Lambda_{\Phi}^{\Psi}(\alpha).$ By the definition of box dimension
we have $\dim_B^\Psi\Sigma_A\leq \D(\Psi)$. Thus we have $
D(\Psi)\leq\Lambda_{\Phi}^{\Psi}(\alpha)=\dim_H^\Psi
E_\Phi(\alpha)\leq \dim_H^\Psi \Sigma_A\leq \dim_B^\Psi\Sigma_A\leq
D(\Psi). $\hfill$\Box$

As a consequence of the previous lemma we have the following
corollary:

\begin{coro}\label{appro-po-metr} Given $\Phi\in \C_{aa}(\Sigma_A,T,d)$ and $\Psi\in
\C_{aa}^{-}(\Sigma_A,T)$.
 Define $\Psi^n$ as in $(\ref{app-potential})$. Assume $\Phi^n\in \C_{aa}(\Sigma_A,T,d)$ is
 such that $\|\Phi-\Phi^n\|_{\mbox{\tiny \rm lim}}\to 0$ when $n\to \infty.$
Assume  $\alpha\in L_\Phi, \alpha_n\in L_{\Phi^n}$ and
$\lim_{n\to\infty}\alpha_n=\alpha$. Then $\Lambda_\Phi^\Psi(\alpha)\geq \limsup_{n\to\infty}\Lambda_{\Phi^{n}}^{\Psi^{n}}(\alpha_n).$
\end{coro}
\proof\ By \eqref{approximation-3}, $
 \Lambda_{\Phi^n}^{\Psi^n}(\alpha_n)\leq \frac{2C_2(\Psi^n)\log m}{|\Psi_{\max}|}\delta_n+(1-C_2(\Psi^n)\delta_n)
 \Lambda_{\Phi}^{\Psi}(\alpha, \|\Phi^n-\Phi\|_{\mbox{\tiny \rm lim}}+\kappa_n\delta_n+2\eta_n),
$
where $\delta_n,\kappa_n$ and $\eta_n$ are the same as in the proof
above. Letting $n\to\infty$  we get the result.

\subsubsection{{\bf $\Lambda_\Phi^\Psi(\alpha)={\tau_\Phi^\Psi}^\star(\alpha)$
when $\Phi$ and $\Psi$ are general and $L_\Phi$ has dimension
$d$}}\label{dim-full}
 We need a last intermediate result regarding
the $\Phi$- and $\Psi$- dependence of the function
${\tau_\Phi^\Psi}^\star.$

Given $\Phi\in\C_{aa}(\Sigma_A,T,d)$ and $\Psi\in
\C_{aa}^{-}(\Sigma_A,T)$, let $\Psi^n$ as in
$(\ref{app-potential})$. Assume $\Phi^n\in \C_{aa}(\Sigma_A,T,d)$
are H\"{o}lder and $\lim_{n\to\infty}\|\Phi-\Phi^n\|_{\mbox{\tiny \rm lim}}= 0$.

For $(z,\alpha)\in (\R^d)^2$ let
${\tau_\Phi^\Psi}(z,\alpha,n)$ be the solution of $P(\langle
z,\Phi^n-\alpha\rangle+{\tau_\Phi^\Psi}(z,\alpha,n)\Psi^n)=0$ and
let
${\tau_\Phi^\Psi}^\star(\alpha,n):=\inf\{{\tau_\Phi^\Psi}(z,\alpha,n):z\in\R^d\}$.

\begin{lem}\label{appro-large-de} Assume $L_\Phi$ has dimension $d$.
Then
 \begin{enumerate}
\item Let $(z,\alpha)\in (\R^d)^2$. We have $\lim_{n\to\infty}{\tau_\Phi^\Psi}(z,\alpha,n)={\tau_\Phi^\Psi}(z,\alpha).$
In particular ${\tau_\Phi^\Psi}(0,\alpha)=D(\Psi)$, and consequently
${\tau_\Phi^\Psi}^\star(\alpha)\leq D(\Psi).$ If $\alpha\in
\mathrm{int}(L_\Phi)$, then
$\lim_{n\to\infty}{\tau_\Phi^\Psi}^\star(\alpha,n)=
{\tau_\Phi^\Psi}^\star(\alpha).$

\item  Let $(z,\alpha,\alpha')\in(\R^d)^3$. If $\langle z,\alpha^\prime-\alpha\rangle\geq 0$, then
\begin{equation}\label{sloope-control-general}
C_1\langle z,\alpha^\prime-\alpha\rangle\leq
{\tau_\Phi^\Psi}(z,\alpha)-{\tau_\Phi^\Psi}(z,\alpha^\prime)\leq
C_2\langle z,\alpha^\prime-\alpha\rangle.
\end{equation}
where $C_1=1/|\Psi_{\min}|$ and $C_2=1+1/|\Psi_{\max}|$.
\end{enumerate}
\end{lem}
The proof will be given in Section~\ref{appendix}.

\noindent {\bf Proof of
$\Lambda_\Phi^\Psi(\alpha)={\tau_\Phi^\Psi}^\star(\alpha)$ when
$L_\Phi$ has dimension $d$.} At first assume $\alpha\in
\mathrm{int}(L_\Phi)$. Since $\Lambda_\Phi^\Psi(\alpha)={\mathcal
E}_\Phi^\Psi(\alpha)$,
  there exists an invariant measure $\mu$ such
that $\Phi_\ast(\mu)=\alpha$ and
$\Lambda_\Phi^\Psi(\alpha)=-h_\mu(T)/\Psi_\ast(\mu).$ Choose a
sequence of ergodic Markov measures $\mu_n$ which converges to $\mu$
in the weak-star topology and such that $h_{\mu_n}(T)$ also conveges
to $h_\mu(T)$ as $n$ tends to $\infty$. We claim that we can choose
a sequence $\{\Phi^n: n\ge 1\}$ of
 H\"{o}lder continuous potentials such that
\begin{equation}\label{good-appro}
 \forall\ n\ge 1,\ \Phi^n_\ast(\mu_n)=\alpha\text{ and } \lim_{n\to\infty}\|\Phi-\Phi^n\|_{\mbox{\tiny \rm lim}}= 0.
\end{equation}
Indeed at first, let  $\widetilde{\Phi}^n$ be a sequence associated
with $\Phi$ as in $(\ref{app-potential})$. Then we know that
$\lim_{n\to\infty}\|\Phi-\widetilde{\Phi}^n\|_{\mbox{\tiny \rm
lim}}= 0$ and each $\widetilde{\Phi}^n$ is a H\"{o}lder potential.
Let $\delta_n=\Phi_\ast(\mu)-\widetilde{\Phi}^n_\ast(\mu_n).$
%By
%definition it is seen that
%$|\Phi_\ast(\mu_n)-\Phi^n_\ast(\mu_n)|\leq \|\Phi-\Phi^n\|_{\mbox{\tiny \rm lim}}$.
We have $
|\delta_n|\leq|\Phi_\ast(\mu)-{\Phi}_\ast(\mu_n)|+|\Phi_\ast(\mu_n)-\widetilde{\Phi}^n_\ast(\mu_n)|
\leq
|\Phi_\ast(\mu)-{\Phi}_\ast(\mu_n)|+\|\Phi-\widetilde\Phi^n\|_{\mbox{\tiny
\rm lim}}.
$
Since $\Phi_\ast$ is continuous, we have $\lim_{n\to\infty}
\delta_n=0$. Define $\Phi^n:=\widetilde{\Phi}^n+\delta_n.$ Thus
$\Phi^n$ satisfies $(\ref{good-appro})$.

Since $\|\Phi^n-\Phi\|_{\mbox{\tiny \rm lim}}\to 0$, it is ready to
show that $d_H(L_{\Phi^n},L_\Phi)\to 0.$ Since $L_{\Phi^n}$ and
$L_\Phi$ are all compact convex sets and $L_\Phi$ has dimension $d$,
$L_{\Phi^n}$ has nonempty interior for $n$ large enough.
Consequently $L_{\Phi^n}$ has dimension $d$ for $n$ large enough.

Now, let $\Psi^n$ be associated with $\Psi$ as in
$(\ref{app-potential})$. By Corollary \ref{appro-po-metr}, we have $
\Lambda_\Phi^\Psi(\alpha)\geq
\limsup_{n\to\infty}\Lambda_{\Phi^{n}}^{\Psi^{n}}(\alpha). $ On the
other hand since $\Phi^n_\ast(\mu_n)=\alpha,$  we have
$\Lambda_{\Phi^n}^{\Psi^n}(\alpha)=\mathcal{E}_{\Phi^n}^{\Psi^n}(\alpha)\geq
-h_{\mu_n}(T)/\Psi^n_\ast(\mu_n).$ So we get
$\displaystyle
\liminf_{n\to \infty}\Lambda_{\Phi^n}^{\Psi^n}(\alpha)\geq
\liminf_{n\to\infty}\frac{h_{\mu_n}(T)}{-\Psi^n_\ast(\mu_n)}
=\frac{h_\mu(T)}{-\Psi_\ast(\mu)}=\Lambda_\Phi^\Psi(\alpha).
$
Thus $\Lambda_\Phi^\Psi(\alpha)=\lim_{n\to
\infty}\Lambda_{\Phi^n}^{\Psi^n}(\alpha)$. For large $n$, since
$\Phi^n$ and $\Psi^n$ are H\"{o}lder continuous and $L_{\Phi^n}$ has
dimension $d$, by Lemma \ref{holder-basic} we have
${\tau_\Phi^\Psi}^\star(\alpha,n)=\Lambda_{\Phi^n}^{\Psi^n}(\alpha).$
Now by Lemma \ref{appro-large-de}(1) since $\alpha\in
\mathrm{int}(L_\Phi)$ we have
${\tau_\Phi^\Psi}^\star(\alpha)=\lim_{n\to\infty}{\tau_\Phi^\Psi}^\star(\alpha,n).$
Then the result follows.

Next we assume $\alpha$ is in the  boundary of $L_\Phi.$
 We have
shown that $D(\Psi)=\max\{\Lambda_\Phi^{\Psi}(\alpha): \alpha\in
L_\Phi\}$. Moreover,  by Lemma \ref{appro-large-de},
${\tau_\Phi^\Psi}^\star(\alpha)\leq D(\Psi)$. Thus, since relation
$(\ref{sloope-control-general})$ holds, the same proof as in the
H\"{o}lder case shows that ${\tau_\Phi^\Psi}^\star(\alpha)=
\Lambda_\Phi^\Psi(\alpha)$.
 \hfill $\Box$

\subsubsection{{\bf  ${\tau_\Phi^\Psi}^\star(\alpha)= \Lambda_\Phi^\Psi(\alpha)$ in the general case}}
$\ $

We only need to show the equality when $\dim L_\Phi=s<d.$ Write $\Phi=(\Phi^1,\cdots,\Phi^d),$  $\widehat
\Phi^l:=(\Phi^1,\cdots,\Phi^s,\Phi^l)$ for $l=s+1,\cdots,d$ and
$\widehat\Phi:=(\Phi^1,\cdots,\Phi^s)$. The set $L_{\widehat \Phi}$
is a projection of $L_{\widehat \Phi^{l}}$ and $L_{\widehat \Phi^l}$ is a
projection of $L_\Phi$ for each $l=s+1,\cdots,d.$ Thus ${\rm
dim}L_{\widehat \Phi}\leq {\rm dim}L_{\widehat \Phi^l}\leq s $ for any
$l=s+1,\cdots,d.$ It is clear that $\langle \bar
\Phi^1,\cdots,\bar\Phi^d\rangle$ has dimension at least $s$,
otherwise, the dimension of  $L_\Phi$ will be strictly less than
$s.$ By relabeling we can assume $\bar \Phi^1,\cdots,\bar
\Phi^{s}$ are linearly independent. By Proposition \ref{dim-L-phi},
we have ${\rm dim}L_{\widehat \Phi}=s.$ Thus we have
$
{\rm dim}L_{\widehat \Phi^l}=s,\ \  \text{for any } \ l=s+1,\cdots,d.
$
Again by Proposition \ref{dim-L-phi}, we conclude that for each
$l=s+1,\cdots,d$, $\bar \Phi^l$ is a linear combination of $\bar
\Phi^1,\cdots,\bar \Phi^{s}.$  Thus there exists a unique $d\times
s$-matrix $A$ of rank $s$ and vector $b\in\R^d$ such that $\|A\widehat
\Phi+b-\Phi\|_{\rm lim}=0.$ Consequently $L_\Phi= AL_{\widehat\Phi}+b.$

For $\widehat \Phi\in \C_{aa}(\Sigma_A,T,s)$, since ${\rm dim}L_{\widehat
\Phi}=s$, by the result proven in section \ref{dim-full}, for any
$\beta\in L_{\widehat \Phi}$ we have
$\Lambda_{\widehat\Phi}^\Psi(\beta)=\tau_{\widehat\Phi}^{\Psi\ast}(\beta)$.
Fix $\alpha\in L_\Phi$ and let $\beta \in L_{\widehat\Phi}$ be the unique
vector such that $\alpha=A\beta+b.$ Then
$E_\Phi(\alpha)=E_{\widehat\Phi}(\beta)$, so
$\Lambda_{\Phi}^\Psi(\alpha)=\tau_{\widehat\Phi}^{\Psi\ast}(\beta)$.

On the other hand,
$
\tau_{\Phi}^{\Psi\ast}(\alpha)=\inf_{z\in\R^d}\tau_{\Phi}^\Psi(z,\alpha)\leq\inf_{\widehat
z\in\R^s}\tau_{\Phi}^\Psi(A\widehat z,\alpha),
$
where $\tau_{\Phi}^\Psi(z,\alpha)$ satisfies $P(\langle
z,\Phi-\alpha\rangle+\tau_{\Phi}^\Psi(z,\alpha)\Psi)=0.$ Since
$\|A\widehat \Phi+b-\Phi\|_{\rm lim}=0$, for any $\lambda\in\R$ we have
$P(\langle A\widehat z,\Phi-\alpha\rangle+\lambda\Psi) =P(\langle A\widehat
z,A(\widehat\Phi-\beta)\rangle+\lambda\Psi)=P(\langle A^\ast A\widehat
z,\widehat\Phi-\beta\rangle+\lambda\Psi)$. Thus we get  $\tau_{\Phi}^\Psi(A\widehat
z,\alpha)=\tau_{\widehat\Phi}^\Psi(A^\ast A\widehat z,\beta)$. Since $A$ has
rank $s$, the square matrix $A^\ast A$ also has rank $s$. This yields
$
\tau_{\Phi}^{\Psi\ast}(\alpha)\leq\inf_{\widehat
z\in\R^s}\tau_{\Phi}^\Psi(A\widehat z,\alpha)=\inf_{\widehat
z\in\R^s}\tau_{\widehat\Phi}^\Psi(A^\ast A\widehat z,\beta)=\inf_{\widehat
z\in\R^s}\tau_{\widehat\Phi}^\Psi(\widehat
z,\beta)=\tau_{\widehat\Phi}^{\Psi\ast}(\beta)=\Lambda_\Phi^\Psi(\alpha).
$
Combining this  with  the inverse inequality, we get the result. \hfill
$\Box$

%%%%%%%%%%%%%%%%%%%%%%%%%%%%%%%%%%%%%%%%%%%%%%%%%%%%%%%%%%%%%%%%%%%%%%%%%%%%%%%%%%%%%%%%%%%%%%%%%%%%%%%%%%%%%%%%%%%%%%%%%%%%%%%

\section{Proof of Theorem~\ref{main-fun-level-one-sided}}\label{proofmainth}

We prove the slighly more general result mentioned in Remark~\ref{remloc}(2).
Suppose that $\xi$ is bounded and continuous outside a subset $E$ of $\Sigma_A$,
and $\xi(\Sigma_A)\subset \mathrm{aff}(L_\Phi)$. Also, suppose that
 $\dim_H^\Psi E<\lambda:=\sup\{\D_\Phi^\Psi(\alpha):\alpha\in \xi (\Sigma_A\setminus E)\cap \mathrm{ri}(L_\Phi)\}.$

In order to prepare the proof of our geometric result, we prove a slightly more general result than necessary.

 \begin{prop}\label{lower-unify-func}
 Assume that $Z\subset \Sigma_A$ is a
closed set such that $\mu(Z)=0$ for any Gibbs measure $\mu$
 supported on $\Sigma_A$. For any $\delta>0$ such that $\lambda-\delta> \dim_H^\Psi(E)$,
  we can construct a Moran subset $\Theta\subset\Sigma_A$
  such that $\Theta\setminus E\subset
 E_\Phi(\xi)$, $\dim_H^\Psi \Theta\geq \lambda-\delta$ and  there
 exists an increasing sequence of integers $(g_j)_{j\ge 1}$ such that $T^{g_j}x\not\in Z$ for any $x\in \Theta$ and any
 $j\geq 1.$
 \end{prop}

\proof\ Fix $\delta>0$ such that $\lambda-\delta> \dim_H^\Psi(E)$. Choose $\alpha_0\in
\xi(\Sigma_A\setminus E)\cap \text{ri}(L_\Phi)$ such that
$\Lambda_{{\Phi}}^\Psi(\alpha_0)>\lambda-\delta/2.$ Assume $L_\Phi$
has dimension $d_0\leq d$ and ${\rm
aff}(L_\Phi)=\alpha_0+U(\R^{d_0}\times \{0\}^{d-d_0})$, where $U$ is an orthogonal matrix.

 Since
$\D_{\Phi}^\Psi$ is continuous in $\text{ri}(L_\Phi)$, we can find
$\eta>0$ such that $B(\alpha_0,\eta)\cap \text{aff}(L_\Phi)\subset
\text{ri}(L_\Phi)$ and for any $\alpha\in B(\alpha_0,\eta)\cap
\text{aff}(L_\Phi)$ we have $
|\D_{\Phi}^\Psi(\alpha)-\D_{\Phi}^\Psi(\alpha_0)|<\delta/2. $
Consequently $\D_{\Phi}^\Psi(\alpha)>\lambda-\delta$ for all
$\alpha\in B(\alpha_0,\eta)\cap \text{aff}(L_\Phi)$. Let $n_0\in\N$
such that $2^{-n_0}\sqrt{d_0}<\eta$ and define a sequence of sets as
follows:
$$
\Delta_n:=B(\alpha_0,\eta)\cap \text{aff}(L_\Phi)\cap
(\alpha_0+2^{-n-n_0}U (\Z^{d_0}\times \{0\}^{d-d_0})),\ \ n\geq 0.
$$
Then $\Delta_0\ne\emptyset$ and $\Delta_n\subset\Delta_{n+1}$ for
any $n\geq 0$ and each $\Delta_n$ is a finite set. For each
$\alpha\in \bigcup_{n\geq 0}\Delta_n$, we can find a measure
$\mu_\alpha$ such that ${\Phi}_\ast(\mu_\alpha)=\alpha$ and
$\D_{\Phi}^\Psi(\alpha)={\mathcal
E}_{\Phi}^\Psi(\alpha)=h_{\mu_\alpha}(T)/\gamma_\alpha,$ where
$\gamma_\alpha=-\Psi_\ast(\mu_\alpha)$.

Let $(\varepsilon_j)_{j\ge 1}$ be a positive sequence such that
$\sum_j\eps_j<\infty$. For each $j\ge 1$ and each $\alpha\in
\Delta_j$, we can find a Markov (hence Gibbs) measure
$\mu_{\alpha,j} $ such that
$$
\max
(|h_{\mu_{\alpha,j}}(T)-h_{\mu_\alpha}(T)|,|\beta_{\alpha,j}-\alpha|,|\gamma_{\alpha,j}-\gamma_\alpha|)<
\eps_j<1,
$$
where $\beta_{\alpha,j}={\Phi}_\ast(\mu_{\alpha,j})$ and
$\gamma_{\alpha,j}=-\Psi_\ast(\mu_{\alpha,j})$. Let
$(\varphi^j)_{j\ge 1}$ and $(\psi^j)_{j\ge 1}$ be two sequences of
H\"older potentials defined on $\Sigma_A$ such that $\|\Phi^j-
\Phi\|_{\rm{lim}}< \eps_j$ and $\|\Psi^j- \Psi\|_{\rm{lim}}<
\eps_j$, where $\Phi^j=(S_n\varphi^j)_{n=1}^{\infty}$ and
$\Psi^j=(S_n\psi^j)_{n=1}^{\infty}$. For each $j\ge 1$, $\alpha\in
\Delta_j$ and $s\in \{1,\dots,m\}$ we denote by $\mu_{\alpha,j}^s$
the restriction of $\mu_{\alpha,j}$ to $[s]$ and
 $\nu_{\alpha,j}^{s}$ the probability measure $\mu_{\alpha,j}^{s}/\mu_{\alpha,j}([s])$.  For $N\ge 1$ let
$$
E^{j}_N(\alpha)=\bigcap_{n\ge N}\left\{x\in \Sigma_A: \Big
|\frac{\phi_{n}^j(Tx)}{n}-\alpha\Big |,\, \Big |\frac{\log
\nu^{s}_{\alpha,j}([x|_{n}])}{-n}-h_{\mu_\alpha}(T)\Big |,\, \Big
|\frac{\psi_{n}^j(Tx)}{-n}-\gamma_\alpha\Big |\le 2\eps_j\right\}.
$$
 Notice
that each $\Delta_j$ is a finite set, thus  the ergodicity of each
$\mu_{\alpha,j}$ imply that we can fix an integer $N_j$ such that
$$
\forall \ N\ge N_j, \ \ \forall \alpha\in \Delta_j,\ \forall
s\in\{1,\cdots,m\} \ \ \ \nu_{\alpha,j}^{s}(E^{j}_N(\alpha))\ge
1-\eps_j/2.
$$

Define $ V_N:=\big \{v\in \Sigma_{A,N+1}:  [v]\cap Z=\emptyset\big
\}. $ There exists $\widehat{N}_j\ge 1$ such that for each $N\ge
\widehat{N}_j$,
$$
\nu_{\alpha,j}^{s}\Big( \bigcup_{v\in V_N}[v]\Big )\ge 1-\eps_j/2, \
\ \forall \alpha\in \Delta_j.
$$

Define $
 V_N^j(\alpha)=\{v\in V_N,\ [v]\cap E^{j}_{N_j}(\alpha)\neq\emptyset\}
.$ Thus, if $N\ge \max(N_j,\widehat{N}_j)$ we have
$$
\nu_{\alpha,j}^{s}\Big (\bigcup_{v\in  V_N^j(\alpha)}[v]\Big )\ge
1-\eps_j,\ \ \forall \alpha\in\Delta_j.
$$

Now we can build a measure on $\Sigma_A$ as follows. At first we
define $\vartheta\in\Sigma_{A,*}$ and inductively a sequence of integers $\{g_j:j\geq 0\}$ and a sequence of measures $\{\rho_j:j\geq 0\}$ such that $\rho_j$ is a
measure on $([\vartheta],\sigma([u]:\vartheta\prec u\in \Sigma_{A,g_j}))$ for each $j\ge 0$, and the
measures $\rho_j$ are consistent: for each $j\ge 0$ the
restriction of $\rho_{j+1}$ to $\sigma([u]:\vartheta\prec u\in \Sigma_{A,g_j})$ is
 equal to $\rho_j$.

Fix $x^0\in
\Sigma_A\setminus E$ such that $\xi(x^0)=\alpha_0$ and write
$x^0=x^0_1x^0_2\cdots.$ Choose  $g_0\in\N$ such that
$\text{O}(\xi,[x^0|_{g_0}])\leq 2^{-n_0}$, where $\text{O}(\xi,V)$ stands for the oscillation of $\xi$ over $V$. Write
$\vartheta:=x^0|_{g_0}.$ Define the probability measure $\rho_0$ to be
the trivial probability measure on $([\vartheta],\{\emptyset,[\vartheta]\})$.
 Suppose we
have defined $(g_k,\rho_k)_{0\le k\le j}$ for $j\ge 0$ as desired.
To obtain $(g_{j+1},\rho_{j+1})$ from $(g_j,\rho_j)$, choose any  $
L_{j+1}\ge \max\{N_{j+1},\widehat{N}_{j+1}\}, $ define
$g_{j+1}=g_j+L_{j+1}$.  For every $w\in \Sigma_{A,g_j}$ with
$\vartheta\prec w$, choose $x_w\in[w]$. Since $x_w\in [w]\subset
[\vartheta]$ we have $ |\xi(x_w)-\alpha_0|=|\xi(x_w)-\xi(x^0)|\leq
2^{-n_0}\leq \eta. $ Notice that by our assumption
$\xi(\Sigma_A)\subset \text{aff}(L_\Phi)$, thus $\xi(x_w)\in
B(\alpha_0,\eta)\cap\text{aff}(L_\Phi).$ Take $\alpha_w\in
\Delta_{j+1}$ such that $|\xi(x_w)-\alpha_w|\leq
2^{-j-1-n_0}\sqrt{d_0}$. For each $v\in\Sigma_{A,L_{j+1}}$ such that
$wv$ is admissible,
 let ($\widehat w$ stands for the last letter of $w$)
$$
\rho_{j+1}([wv]):=\rho_j([w])\nu^{\widehat{w}}_{\alpha_w,j+1}([\widehat{w}v]).
$$
By construction the family $\{\rho_j:j\geq 0\}$ is consistent.
Denote by $\rho$ the Kolmogorov extension of the sequence
$(\rho_j)_{j\ge 0}$ to $([\vartheta],\sigma([u]:\vartheta\prec u\in \Sigma_{A,*})$. 

If $\vartheta\prec u$ and $u\in \Sigma_{A,n}$ with $g_j\le n<g_{j+1}$, writing $u=\vartheta
w^1\cdots w^j\cdot v$ with $|w^k|=L_k$ and  $|v|=n-g_i,$ and denoting
$\vartheta w^1\cdots w^k$ by $\widetilde{w}_k$, we have the useful formula
\begin{equation}\label{rho-fun-level}
 \rho([u])=\Big( \prod_{k=1}^j\nu^{\widehat{w^{k-1}}}_{\alpha_{\widetilde{w}_{k-1}},k}([\widehat{w^{k-1}}w^k])\Big )
 \nu^{\widehat{w^{j}}}_{\alpha_{\widetilde{w}_j},j+1}([\widehat{w^{j}}v]).
\end{equation}

Let $ \Theta=\big \{x\in [\vartheta]: \ \forall\ j\ge 1, \
T^{g_j-1}x|_{L_{j+1}+1}\in  V_{L_{j+1}}^{j+1}(\alpha_{x|_{g_j}})\
\big\}. $ By construction, $T^{g_j-1}x\not\in Z$ for any $x\in
\Theta$ and any $j\geq 1.$

Write $\alpha_k:=\alpha_{x|_{g_{k-1}}}$. By construction of $\rho$,
for each $j\ge 1$, by using \eqref{rho-fun-level} we can get
\begin{eqnarray*}
&&\rho(\{x\in [\vartheta]: [x|_{g_j}]\cap \Theta\neq\emptyset\})\\
&=&\sum_{ \substack{ \vartheta w^1\cdots w^j\ \text{admissible},\\
\forall\, 1\le k \le j,\,  \widehat{w^{k-1}}w^k\in V_{L_k}^k(\alpha_k)}} \rho_j(\vartheta w^1\cdots w^j)\\
&=&\sum_{  \substack{ \vartheta w^1\cdots w^{j-1}\ \text{admissible},\\
\forall\, 1\le k \le j-1,\,  \widehat{w^{k-1}}w^k\in
V_{L_k}^k(\alpha_k)}}\rho_{j-1}(\vartheta w^1\cdots w^{j-1})
\sum_{\widehat{w^{j-1}}w^j\in   V_{L_j}^j(\alpha_j)}\nu^{(\widehat{w^{j-1}})}_{\alpha_j,j}([\widehat{w^{j-1}}\cdot w^j])\\
&\ge &\sum_{  \substack{ \vartheta w^1\cdots w^{j-1}\ \text{admissible},\\
\forall\, 1\le k \le j-1,\,  \widehat{w^{k-1}}w^k\in
V_{L_k}^k(\alpha_k)}} \rho_{j-1}(\vartheta w^1\cdots w^{j-1})(1-\eps_j)\ge
\prod_{k=1}^j(1-\eps_k).
\end{eqnarray*}
Consequently, $\rho(\Theta)\ge \prod_{j\ge 1}(1-\eps_j)>0$ since we
assumed that $\eps_j<1$ and $\sum_{j\ge 1}\eps_j<\infty$.

For $\eta\in \{\varphi,\psi\}$ and $j\ge 1$, let
$$
c(\eta^j)=\sup_{n\ge
1}\max_{v\in\Sigma_{A,n}}\max_{x,y\in[v]}|S_n\eta^j(x)-S_n\eta^j(y)|.
$$
This number is finite since each $\eta^j$ is H\"older continuous.
Let $M_j\nearrow \infty$  such that
$$
\forall\ n\ge M_j,\
\max(\|\phi_{n}^j-\phi_n\|_\infty,\|\psi_{n}^j-\psi_n\|_\infty)\le
2\eps_j n.
$$
 The sequence $(L_j)_{j\ge 1}$ can be specified to satisfy the additional properties
$$
 L_j\ge M_{j+1} \ \ \text{and }\ \
\max(K_1(j), K_2(j),K_3(j))\le \eps_jg_j,
$$
(recall that $g_j=g_0+\sum_{k=1}^j L_k$), where
$$
\begin{cases}
K_1(j)=\sum_{k=1}^{j+1}(c(\varphi^k)+c(\psi^k))\\
K_2(j)=\max_{\substack{\alpha\in\Delta_{j+1}\\0\le s\le m-1\\1\le
n\le N_{j+1}}}\max(n|\alpha|,\|\phi_{n}^{j+1}\|_\infty,\|
\log\nu_{\alpha,j+1}^{s}([\cdot|_{n}])\|_\infty,\|\psi_{n}^{j+1}\|_\infty)\\
K_3(j)=(j+1)\max_{1\le n\le M_{j+1}} \max (\|\phi_{n}^{j+1}-
\phi_{n}\|_\infty,\|\psi_{n}^{j+1}-\psi_n\|_\infty)
\end{cases}.
$$

\bigskip

Let us check that $\Theta\setminus E\subset E_{{\Phi}}({\xi}).$ Let
$x\in\Theta\setminus E$, $n\ge g_1$ and $j\ge 1$ such that $g_j\le
n<g_{j+1}$. Since $g_j>L_j\ge M_{j+1}$, we have
$$
|\phi_n (x)-n\xi(x)|\le \|\phi_{n}^{j+1}-\phi_n\|_\infty+
|\phi_{n}^{j+1}(x)-n\xi(x)|\le 2\eps_{j+1}
n+|\phi_{n}^{j+1}(x)-n\xi(x)|.
$$
We have (with $g_{-1}=0$, $\alpha_k=\alpha_{x|_{g_{k-1}}}$ and  $L_0=g_0$)
\begin{eqnarray*}
&&|\phi_{n}^{j+1}(x)-n\xi(x)|\\
&\leq&|\phi_{g_j}^{j+1}(x)-g_j\xi(x)|+|\phi_{n-g_j}^{j+1}(T^{g_j}x)-(n-g_j)\xi(x)|\\
&=&|\phi_{g_j}^{j+1}(x)-\sum_{k=0}^jL_k\alpha_k+\sum_{k=0}^jL_k\alpha_k-
  g_j\xi(x)|+|\phi_{n-g_j}^{j+1}(T^{g_j}x)-(n-g_j)\xi(x)|\\
  &\leq&\sum_{k=0}^j|\phi_{L_k}^{j+1}(T^{g_{k-1}}x)-L_k\alpha_k|+\sum_{k=0}^jL_k|\alpha_k-
  \xi(x)|\ \ \big(=: (I)+(II) \big)\\
  &+&|\phi_{n-g_j}^{j+1}(T^{g_j}x)-(n-g_j)\alpha_{j+1}|+(n-g_j)|\alpha_{j+1}-\xi(x)|\ \ \big(=: (III)+(IV) \big)\\
\end{eqnarray*}
At first we have
%%%%%%%%%%%%%%%%%%%%%%%%%%%%%%%%%%%%%%%%%%%%%%%%%%%%%%%%%%%%%%%%%
\begin{eqnarray*}
(I)+(III)
&\le& \sum_{k=0}^j\|\phi_{L_k}^{j+1}-\phi_{L_k}\|_\infty+\sum_{k=0}^j\|\phi_{L_k}^k-\phi_{L_k}\|_\infty\\
&&+\Big ( \sum_{k=0}^j|\phi_{L_k}^k(T^{g_{k-1}}x)- L_k\alpha_k|\Big
)+ |\phi_{n-g_j}^{j+1}(T^{g_j}x)-(n-g_j)\alpha_{j+1}|.
\end{eqnarray*}
If $L_k\leq M_{j+1}$, then
$\|\phi_{L_k}^{j+1}-\phi_{L_k}\|_\infty\leq K_3(j)/(j+1)$; if $L_k>
M_{j+1}$, then $\|\phi_{L_k}^{j+1}-\phi_{L_k}\|_\infty\leq
2\varepsilon_kL_k.$ Thus we have
$\sum_{k=0}^j\|\phi_{L_k}^{j+1}-\phi_{L_k}\|_\infty\le
K_3(j)+2\sum_{k=0}^j\eps_kL_k$. Since $L_k\geq M_{k+1}\geq M_k$ we
also have $\sum_{k=0}^j\|\phi_{L_k}^k-\phi_{L_k}\|_\infty\le
2\sum_{k=0}^j\eps_kL_k$. Thus both terms are $o(g_j)$ as
$n\to\infty$. Consequently both terms are $o(n).$

For $k=0,\cdots,j$, by the construction of $\Theta$, we have
$T^{g_{k-1}-1}x|_{L_k+1}=x_{g_{k-1}}\cdot (T^{g_{k-1}}x|_{L_k})\in
V_{L_k}^k(\alpha_k)$, so
 $[x_{g_{k-1}}\cdot (T^{g_{k-1}}x|_{L_k})]\cap E^{k}_{N_k}(\alpha_k)\neq\emptyset$. Since $L_k\geq
 N_k,$
there exists $y\in [T^{g_{k-1}}x|_{L_k}]$ such that
$|\phi_{L_k}^k(y)-L_k\alpha_k|\le 2\eps_kL_k$, hence
 $|\phi_{L_k}^k(T^{g_{k-1}}x)-L_k\alpha_k|\le 2\eps_kL_k+c(\varphi^k)$. Similarly, we have
$|\phi_{n-g_j}^{j+1}(T^{g_j}x)-(n-g_j)\alpha_{j+1}|\le
2\eps_{j+1}(n-g_j)+c(\varphi^{j+1})$ if $n-g(j)\ge N_{j+1}$, and we
trivially have
$|\phi_{n-g_j}^{j+1}(T^{g_j}x)-(n-g_j)\alpha_{j+1}|\le 2K_2(j)$
otherwise. This yields
\begin{eqnarray*}
&& \Big ( \sum_{k=0}^j|\phi_{L_k}^k(T^{g_{k-1}}x)- L_k\alpha_k|\Big )+ |\phi_{n-g_j}^{j+1}(T^{g_j}x)-(n-g_j)\alpha_{j+1}|\\
&\le &
 2\sum_{k=0}^j\eps_kL_k+\sum_{k=0}^{j+1}c(\varphi^k)+  2\eps_{j+1}(n-g_j)+2K_2(j)\\
 &\le& K_1(j)+2K_2(j)+  2\eps_{j+1}(n-g_j)+ 2\sum_{k=0}^j\eps_kL_k=o(g_j)=o(n).
 \end{eqnarray*}
Together we get $(I)+(III)=o(n)$.
%%%%%%%%%%%%%%%%%%%%%%%%%%%%%%%%%%%%%%%%%%%%%%%%%%%%%%%%%%%%%%%%%%%%%%%%%%%5
On the other hand, by construction,
$$
|\xi(x)-\alpha_{x|_{g_k}}|\leq
|\xi(x)-\xi(x_{x|_{g_k}})|+|\xi(x_{x|_{g_k}})-\alpha_{x|_{g_k}}|\leq
O(\xi,[x|_{g_k}])+2^{-k-n_0}\sqrt{d_0},
$$
and $\lim_{k\to\infty}O(\xi,[x|_{g_k}])+2^{-k-n_0}\sqrt{d_0}=0$
since $\xi$ is continuous at $x$. Thus we conclude that
$(II)+(IV)=o(n).$ Finally, $| \phi_n (x)-n\xi(x)|=o(n)$, and
$\Theta\setminus E\subset E_{ \Phi}({\xi})$.

Similarly, for any $x\in\Theta$ we can prove that
($\alpha_k:=\alpha_{x|_{g_{k-1}}}$)
\begin{eqnarray*}
&&|-\psi_{n}(x)
-\sum_{k=0}^{j}L_k\gamma_{\alpha_k}-(n-g_j)\gamma_{\alpha_{j+1}}|=o(n)\\
&& |-\log\rho([x|_n])
-\sum_{k=0}^{j}L_kh_{\mu_{\alpha_k}}(T)-(n-g_j)h_{\mu_{\alpha_{j+1}}}(T)|=o(n).
\end{eqnarray*}
By construction $\liminf_{j\to\infty}
h_{\mu_{\alpha_j}}(T)/\gamma_{\alpha_j}\geq \lambda-\delta$. For any
$y\in [x|_n]$ we have $|\psi_n(y)-\psi_n(x)|=o(n)$, thus we get
$\text{diam}([x|_n])=\Psi[x|_n]=\exp(\psi_n(x)+o(n))$. Combining the
above two relations  we conclude that
$\liminf_{n\to\infty}\log\rho([x|_{n}])/\log({\rm
diam}([x|_{n}])\geq \lambda-\delta$. That is the local lower dimension of
$\rho$ at each $x\in\Theta$ is $\geq \lambda-\delta,$ hence
$\dim_H^\Psi(\Theta)\geq \lambda-\delta$ by the mass distribution principle (see \cite{P} for instance). \hfill$\Box$

By essentially repeating the same proof as above (in fact, it is easier), we
can get the following property:

 \begin{prop}\label{lower-unify}
  Assume $Z\subset \Sigma_A$ is a
closed set such that $\mu(Z)=0$ for any Gibbs measure $\mu$
 supported on $\Sigma_A$. For any $\alpha\in L_\Phi,$ we can construct a  subset $\Theta\subset
 E_\Phi(\alpha)$ such that $\dim_H^\Psi\Theta\ge {\mathcal E}_\Phi^\Psi(\alpha) $ and  there exists an integer sequence $g_j\nearrow
 \infty$ such that $T^{g_j}x\not\in Z$ for any $x\in \Theta$ and any
 $j\geq 1.$ In particular, ${\mathcal E}_\Phi^\Psi(\alpha)\leq
 \D_\Phi^\Psi(\alpha)$.
 \end{prop}

 \noindent{\bf Proof of Theorem~\ref{main-fun-level-one-sided}.}
(1')\ Since $\dim_H^\Psi E<\lambda-\delta,$  by the proposition above
we have $\dim_H^\Psi (\Theta\setminus E)=\dim_H^\Psi
\Theta\geq\lambda-\delta.$ Consequently $ \dim_H^\Psi
E_\Phi(\xi)\geq \dim_H^\Psi (\Theta\setminus E)\geq \lambda-\delta.
$ Since $\delta>0$ is arbitrary, we get $\dim_H^\Psi E_\Phi(\xi)\geq
\lambda$.

(2) If $\xi(\Sigma_A)\subset L_\Phi$, the construction of a Moran subset of
 $E_\Phi(\xi)$ can be done around any point of $\Sigma_A$, like in the proof of
 Proposition~\ref{lower-unify-func}. The only difference is that in this case
   the dimension of this set is of no importance. Hence, $E_\Phi(\xi)$  is dense.

(3)\ Now we assume $\xi$ is continuous everywhere. If moreover
$$
\sup\{\D_\Phi^\Psi(\alpha):\alpha\in \xi (\Sigma_A)\cap
\mathrm{ri}(L_\Phi)\}= \sup\{\D_\Phi^\Psi(\alpha):\alpha\in
 {\xi (\Sigma_A)}\cap L_\Phi\}=:\theta,
$$ then at first we have $\dim_H^\Psi E_\Phi(\xi)\geq \theta.$
On the other hand by  definition  we have $E_\Phi(\xi)\subset
E_\Phi( {\xi (\Sigma_A)}\cap L_\Phi)$. Thus by Proposition
\ref{upper-bound}, we have $\dim_P^\Psi E_\Phi(\xi)\leq \theta.$ So
we get $\dim_H^\Psi E_\Phi(\xi)=\dim_P^\Psi E_\Phi(\xi)= \theta$.

(4)\ Assume $d=1$, $\xi$ is continuous everywhere and
$\xi(\Sigma_A)\subset L_\Phi.$ Notice that in this case
$L_\Phi=[\alpha_1,\alpha_2]$ is an interval. Assume
$\alpha_0\in\xi(\Sigma_A)$ such that
$\D_\Phi^\Psi(\alpha_0)=\sup\{\D_\Phi^\Psi(\alpha):\alpha\in
\xi(\Sigma_A)\}$. If $\alpha_0\in (\alpha_1,\alpha_2),$ by (2) we
conclude. Now assume $\alpha_0=\alpha_1$. If $\alpha_1$ is not
isolated in $\xi(\Sigma_A),$ still by (2) and the continuity of
$\D_\Phi^\Psi$, we get the result. If $\alpha_1$ is isolated in
$\xi(\Sigma_A)$, then by the continuity of $\xi$, we can find a
cylinder $[w]\subset \Sigma_A$ such that $\xi([w])=\alpha_1.$ From
this we get $E_\Phi(\xi)\supset E_\Phi(\alpha_1)\cap[w]$. Thus
$\dim_H^\Psi E_\Phi(\xi)\geq \D_\Phi^\Psi(\alpha_1)$ and the result
holds. If $\alpha_0=\alpha_2$, the proof is the same.
 \hfill
$\Box$

\section{Proofs of results in section \ref{examples}}\label{sec7}

We will use the following lemma, which is standard and essentially the same as
Lemma 5.1 in  \cite{GP97} (the proof is elementary).
\begin{lem} \label{dim-compare-lem}
Let $X$ and $Y$ be metric spaces and $\chi: X\to Y$ a surjective
mapping with the following property: there exists a function
$N:(0,\infty)\to \N$ with  $\log N(r)/\log r\to 0$ when $r\to 0$
such that for any $r>0$, the pre-image $\chi^{-1}(B)$ of any
$r$-ball in $Y$ can be covered by at most $N(r)$ sets in $X$ of
diameter less than $r$. Then for any set $E\subset Y$ we have
$\dim_H E\geq \dim_H \chi^{-1}(E).$
\end{lem}

\noindent{\bf Proof of Proposition \ref{dim-compare-prop}.} \
Condition (4) implies that $\chi:(\Sigma_A,d_\Psi)\to (J,d)$ is
Lipschitz continuous, thus we have $\dim_H E\leq
\dim_H^\Psi\chi^{-1}(E).$

For the converse inequality, let us check the condition of the above
lemma. Let $B\subset J$ be a ball of radius $r$, let $n\in \N$ such
that $e^{-n}\leq r< e^{1-n}$. define
 $$
G_B^r=\{w\in {\mathcal B}_n(\Psi): R_w\cap B\ne \emptyset\}.
 $$
 One checks that $\{[w]:w\in G_B^r\}$ is an $r$-covering of
 $\chi^{-1}(B)$. Define $N(r):=\#G_B^r.$
  Let us estimate the number $\#G_B^r.$ Clearly, $\#G_B^r\geq 1.$
  By  condition (4), for each $w\in
 G_B^r$, $R_w$ is contained in a ball of radius $K\Psi[w]\leq
 Ke^{-n}$, thus $\bigcup_{w\in G_B^r}R_w\subset B(y,r+2Ke^{-n})\subset B(y,(e+2K)e^{-n})$, where $y$ is the center of $B$. On
 the other hand, by Lemma \ref{Moran-covering}(1) there exists $C>0$ such that $|w|\leq
 Cn$ for any $w\in{\mathcal B}_n(\Psi)$, thus
 $\eta_{|w|}=o(|w|)=o(n)$ for any $w\in{\mathcal B}_n(\Psi)$.
By  construction, the interiors of the sets $R_w$, $ w\in G_B^r$,
are disjoint and each $R_w$ contains a ball
 of radius $K^{-1}\exp(\eta_{|w|})\Psi[w]=K^{-1}e^{o(n)}\Psi[w]=K^{-1}e^{-n+o(n)}$ by Lemma
 \ref{Moran-covering} (2). Thus $\#G_B^r\leq K^{d^\prime}(e+2K)^{d^\prime}e^{o(n)}$.
 So we conclude that $\log N(r)/\log r=\log \#G_B^r/\log r\to 0$
 as $r\to 0.$ Thus by lemma \ref{dim-compare-lem}, we  can conclude
 that $\dim_H E\geq \dim_H^\Psi\chi^{-1}(E).$
\hfill$\Box$

\noindent{\bf Proof of Lemma \ref{boundary}.}\ At first we show that
$J\cap V= J\setminus \widetilde Z_\infty,$ consequently by the SOSC,
$J\setminus \widetilde Z_\infty\ne\emptyset$ and  we get
$\emptyset\ne \chi^{-1}(J\setminus \widetilde
Z_\infty)=\Sigma_A\setminus Z_\infty.$ In fact
\begin{eqnarray*}
y\in J\setminus \widetilde Z_\infty &\Leftrightarrow& y\in J \text{
and }  \forall\ n\ge 1\ \exists\ x\in\Sigma_A \ \ s.t.\ \ y\in
\text{int} (R_{x|_n})=f_{x|_n}(V) \\
&\Leftrightarrow& y\in J  \text{ and } \forall\ n\ge 1\ \exists!\
x\in\Sigma_A \ \ s.t.\ \ y\in
\text{int} (R_{x|_n})=f_{x|_n}(V)\Leftrightarrow y\in J\cap V.
\end{eqnarray*}

By construction, $\chi: \Sigma_{A}\setminus  Z_\infty\to J\setminus
 \widetilde Z_\infty$ is surjective. Since $ J\setminus \widetilde
Z_\infty=J\cap V,$ it is ready to show that $\chi$ is also
injective.

Next we show that $T(\Sigma_A\setminus Z_\infty)\subset
\Sigma_A\setminus Z_\infty.$ Take $x\in \Sigma_A\setminus Z_\infty$.
If
 $Tx\in Z_\infty,$ then we can find $n_0\in \N$ such that
$\chi(Tx)\in f_{Tx|_{n_0}}(\partial V)$. Consequently $
\chi(x)=f_{x_1}(\chi(Tx))\in f_{x_1}(f_{Tx|_{n_0}}(\partial V))=
f_{x_1}\circ f_{Tx|_{n_0}}(\partial V)= f_{x|_{n_0+1}}(\partial V)
$, which is a contradiction. At last we show that for any Gibbs
measure $\mu$ we have $\mu(Z_\infty)=0.$ Define $ \widetilde
Z_n:=\bigcup_{w\in \Sigma_{A,\ast},\ |w|\leq n}f_w(\partial V) $ and
$Z_n=\chi^{-1}(\widetilde Z_n)$. The sequence $(Z_n)_{n\ge 1}$ is
non decreasing and $Z_\infty=\bigcup_{n\geq 1}Z_n.$ Since the IFS is
conformal we can easily get $T(Z_n)\subset Z_{n-1}$ for $n\geq 1$
and $T(Z_0)\subset Z_0.$ Consequently $T(Z_n)\subset Z_{n}.$ By the
ergodicity we have $\mu(Z_n)=0$ or $1$. By the SOSC,
  $\Sigma_A\setminus Z_n$ is nonempty and open, thus  by
the Gibbs property of $\mu$ we get $\mu(\Sigma_A\setminus Z_n)>0$, hence $\mu(Z_n)=0$. Consequently $\mu(Z_\infty)=0.$

From $T(Z_n)\subset Z_{n}$ we easily get $T(Z_\infty)\subset
Z_\infty.$ \hfill$\Box$

\bigskip

\noindent{\bf Proof of Theorem~\ref{appli-one-sided}.}\ (1) At first
we notice that by the property (4) assumed in the construction of
$J$ the mapping $\chi$ is Lipschitz. This is enough to get the
desired upper bounds from Theorem~\ref{main-one-sided}(1).

Now we deal with the lower bound for dimensions and the equality $L_\Phi=L_{\widetilde\Phi}$.
 We notice that the inclusion $L_\Phi\subset L_{\widetilde\Phi}$ holds by construction.

\noindent{\it Suppose $J$ is a conformal repeller.} Since we have
$\chi\circ T=g\circ \chi$ on $\Sigma_A$ and $\chi$ is surjective, it
is seen that $\chi^{-1}(E_\Phi(\alpha))=E_{\widetilde{\Phi}}(\alpha)$
for any $\alpha\in L_{\widetilde{\Phi}}.$ Thus $L_\Phi=L_{\widetilde{\Phi}}$
and by Proposition \ref{dim-compare-prop}, we have
$\D_\Phi(\alpha)=\D_{\widetilde{\Phi}}^\Psi(\alpha).$

\noindent{\it Suppose $J$ is the attractor of a  conformal IFS with
SOSC.} Let $\alpha\in L_{\widetilde \Phi}$. Let
$Z=\chi^{-1}(\partial V)$. The set $Z$ is closed and by Lemma \ref{boundary},  $\mu(Z)=0$ for any Gibbs
measure $\mu$. By Proposition \ref{lower-unify} we can construct a
Moran set $\Theta\subset E_{\widetilde{\Phi}}(\alpha)$ such that
$\dim_H^\Psi(\Theta)\geq {\mathcal E}_{\widetilde
\Phi}^\Psi(\alpha)=\D_{\widetilde{\Phi}}^\Psi(\alpha)$ and there
exists a sequence $g_j\nearrow\infty$ such that $T^{g_j}x\not\in Z$
for any $x\in\Theta$ and any $j\geq 1.$ The last property means that
$\Theta\subset \Sigma_A\setminus Z_\infty.$ Since $\chi$ is a
bijection between $\Sigma_A\setminus Z_\infty$ and
$J\setminus\widetilde Z_\infty$, we conclude that
$\chi^{-1}\circ\chi(\Theta)=\Theta,$ thus by Proposition
\ref{dim-compare-prop}, $ \dim_H\chi(\Theta)=\dim_H^\Psi
\Theta\geq\D_{\widetilde{\Phi}}^\Psi(\alpha).$ Since we also have
$\chi\circ T=\widetilde{g}\circ\chi$ on $\Sigma_A\setminus
Z_\infty,$ we get that $\chi(\Theta)\subset E_\Phi(\alpha).$ Thus
$\alpha\in L_\Phi$ and $\D_\Phi(\alpha)\geq
\D_{\widetilde{\Phi}}^\Psi(\alpha)$.

\noindent (2) Take $E=J$ in Proposition \ref{dim-compare-prop}, then
use (1) and Theorem \ref{main-one-sided}(2). \hfill $\Box$

\noindent{\bf Proof of Theorem~\ref{appli-fun-level-one-sided}.}\
Define $\widetilde{\xi}:=\xi\circ\chi.$

\noindent{\it Case 1:}  $J$ is a conformal repeller.  One checks easily
that $\chi^{-1}(E_\Phi(\xi))=E_{\widetilde{\Phi}}(\widetilde{\xi}).$ Then
the result is a consequence of Theorem \ref{appli-one-sided} and
Theorem \ref{main-fun-level-one-sided}.

\noindent{\it Case 2:}  $J$ is the attractor of a conformal IFS with
SOSC. By using Proposition \ref{lower-unify-func}, Theorem
\ref{main-fun-level-one-sided} and the same argument as in the proof of
Theorem \ref{appli-one-sided}, we get the result.
 \hfill$\Box$

\noindent{\bf Proof of Theorem \ref{carpet}.}\ Let $\xi=\Phi,$
then ${\mathcal
 F}(\widetilde{J},\widetilde{g})=E_{\Phi}(\xi)$. To show the result we need
 only to check the condition of Theorem
 \ref{appli-fun-level-one-sided} and the only condition we need to
 check is that
 $$
\sup\{\D_\Phi(\alpha): \alpha\in \xi(J)\cap{\rm
ri}(L_\Phi)\}=\sup\{\D_\Phi(\alpha): \alpha\in \xi(J)\cap L_\Phi\}.\
\ (\ast)
 $$
Notice that in this special case we have $\xi(J)=J$ and $L_\Phi={\rm
Co}(J),$ thus $\xi(J)\cap L_\Phi=J.$  Recall  that in this case
$L_\Phi$ is a convex polyhedron, thus by Proposition
\ref{regularity} and Theorem \ref{appli-one-sided}, $\D_\Phi$ is
continuous on $L_\Phi.$ Thus the supremum  in the right hand side of
$(\ast)$ can be reached. If the maximum is attained in  ${\rm
ri}(L_\Phi)$, then the result is obvious. Now suppose that there
exists  $\alpha_0\in
\partial L_\Phi\cap J$ such that $\D_\Phi(\alpha_0)=\sup\{\D_\Phi(\alpha):
\alpha\in J\}.$
 By the structure of $J$, it is ready to see that $B(\alpha_0,\delta)\cap
J\cap {\rm ri}(L_\Phi)\ne \emptyset$ for any $\delta>0$. By the
continuity of $\D_\Phi$, $(\ast)$ holds immediately.
 \hfill $\Box$

\section{Appendix}\label{appendix}
\noindent
{\bf Proof of Lemma~\ref{Psi-dependence}.} (1) \ Write $\Psi=(\psi_n)_{n=1}^{\infty}$ and $\Upsilon=(\upsilon_n)_{n=1}^{\infty}.$
  By the definition of $\|\cdot\|_{\mbox{\tiny \rm lim}}$, for any $\delta>\|\Psi-\Upsilon\|_{\mbox{\tiny \rm lim}},$
  there exist $N\in\N,$ such that for any $n\geq N$ we have $
 \psi_n(x)-n\delta\leq \upsilon_n(x)\leq \psi_n(x)+n\delta,
 $ hence for any $w\in \Sigma_{A,\ast}$ with $|w|$ large enough  %
 $\Psi[w]\exp (-|w|\delta)\leq\Upsilon[w]\leq\Psi[w] \exp (|w|\delta).$
 Given $w\in {\mathcal B}_n(\Psi)$, we have $e^{\Psi_{\min}-C(\Psi)-\|\Psi\|_{|w|}-n}\leq \Psi[w]\leq
 e^{-n}$ and $C_1(\Psi)n\leq |w|\leq C_2(\Psi)n,$ where $C_1(\Psi)=1/|\Psi_{\min}|$ and $C_2(\Psi)=1+1/|\Psi_{\max}|$. So we have
 $$
 e^{\Psi_{\min}-C(\Psi)-\|\Psi\|_{n}^\star-n(1+C_2(\Psi)\delta)}\leq\Upsilon[w]\leq e^{-n(1-C_2(\Psi)\delta)}.
 $$
 This implies that there exists $u\prec w$ such
that $u\in {\mathcal B}_{[n(1-C_2(\Psi)\delta)]}(\Upsilon)$. So we
conclude that
\begin{equation}\label{upper}
\#{\mathcal B}_{[n(1-C_2(\Psi)\delta)]}(\Upsilon)\leq \#{\mathcal
B}_{n}(\Psi).
\end{equation}

Let  $c_1(n)=-\Psi_{\min}+C(\Psi)+\|\Psi\|_{n}^\star$, then
$c_1(n)>0$ and $c_1(n)=o(n)$.  Write $w=uw^\prime$. The same
proof as that of
 the claim in Proposition \ref{dim-formular-1} yields $|w^\prime|\leq (c_1(n)+2nC_2(\Psi)\delta+C(\Upsilon))/
 |\Upsilon_{\max}|.$ Thus we can conclude that
\begin{equation}\label{lower}
\#{\mathcal B}_{[n(1-C_2(\Psi)\delta)]}(\Upsilon)\geq \#{\mathcal
B}_{n}(\Psi)m^{-(c_1(n)+2nC_2(\Psi)\delta+C(\Upsilon))/
 |\Upsilon_{\max}|}.
 \end{equation}
 Combining $(\ref{upper})$, $(\ref{lower})$ and $(\ref{full-dim})$ we get
 $$(1-C_2(\Psi)\delta)D(\Upsilon)\leq D(\Psi)\leq (1-C_2(\Psi)\delta)D(\Upsilon)-2C_2(\Psi)\delta\log m/|\Upsilon_{\max}|.$$
By using  $(\ref{ful-dim-1})$  we get $|D(\Psi)-D(\Upsilon)|\leq
a(m,\Psi,\Upsilon)\delta,$ where
$$a(m,\Psi,\Upsilon)=2C_2(\Psi)C_2(\Upsilon)\log
m=2\Big (1+\frac{1}{|\Psi_{\max}|}\Big )\Big (1+\frac{1}{|\Upsilon_{\max}|}\Big )\log
m.$$ Since $\delta>\|\Psi-\Upsilon\|_{\mbox{\tiny \rm lim}}$ is arbitrary, we get
$|D(\Psi)-D(\Upsilon)|\leq
a(m,\Psi,\Upsilon)\|\Psi-\Upsilon\|_{\mbox{\tiny \rm lim}}.$

(2) Now given  $\Phi,\Theta\in\C_{aa}(\Sigma_A,T,d)$. Assume
$0<\epsilon<\|\Phi\|$ and $\beta\in L_\Theta$.  Fix $\alpha\in
B(\beta,\eta)\cap L_\Phi.$ Fix $\delta>\|\Psi-\Upsilon\|_{\mbox{\tiny \rm lim}}$.

For any $\epsilon>0$, pick up  $w\in F(\alpha,n,\epsilon,\Phi,\Psi).$
Then $w\in {\mathcal B}_n(\Psi)$
 and there exists $x\in [w]$ such that $|\phi_{|w|}(x)-|w|\alpha|\leq |w|\epsilon.$
We have seen in proving (1) that  $w=uw^\prime$, where $u\in {\mathcal
B}_{[n(1-C_2(\Psi)\delta)]}(\Upsilon)$ and $|w^\prime|\leq
(c_1(n)+2nC_2(\Psi)\delta+C(\Upsilon))/
 |\Upsilon_{\max}|.$ Notice that $\mathrm{diam}(L_\Phi)\leq \|\Phi\|,$ thus $|\alpha|\leq \|\Phi\|.$ So we have
 \begin{eqnarray*}
 |\phi_{|u|}(x)-|u|\alpha|&\leq&
 |\phi_{|w|}(x)-|w|\alpha|+|w^\prime|(\|\Phi\|+|\alpha|)+|C(\Phi)|\\
 &\leq&
 |w|\epsilon+2|w^\prime|\|\Phi\|+|C(\Phi)|\\
 &\leq&|u|\Big (\epsilon+\frac{3\|\Phi\||w^\prime|+|C(\Phi)|}{|u|}\Big ).
 \end{eqnarray*}
Since $0<c_1(n)=o(n)$, for large $n$ we have
$$
3\|\Phi\|(c_1(n)+C(\Upsilon))+C(\Phi)|\Upsilon_{\max}|\leq
n\|\Phi\|C_2(\Psi)\delta.
$$
Combining this with $|u|\geq
C_1(\Upsilon)n(1-C_2(\Psi)\delta)$ we get that for $\delta<
1/(2C_2(\Psi))$,
$$(3\|\Phi\||w^\prime|+|C(\Phi)|)/|u|\leq 14\|\Phi\|C_2(\Psi)\frac{|\Upsilon_{\min}|}
{|\Upsilon_{\max}|}\delta=:\kappa(\Psi,\Upsilon,\Phi)\delta=\kappa\delta.$$

Fix any $a>\|\Theta-\Phi\|_{\mbox{\tiny \rm lim}}$. For $n$ large enough we have
\begin{eqnarray*}
|\theta_{|u|}(x)-|u|\beta|&=&|\theta_{|u|}(x)-\phi_{|u|}(x)|+|\phi_{|u|}(x)-|u|\alpha|+|u||\alpha-\beta|\\
&\leq&a|u|+(\epsilon+\kappa\delta)|u|+\eta|u|.
\end{eqnarray*}
As a result $u\in F(\beta,
n(1-C_2(\Psi)\delta),a+\epsilon+\kappa\delta+\eta,\Theta,\Upsilon).$
Thanks to our control of $|w^\prime|$, we have
$$f(\beta, n(1-C_2(\Psi)\delta),a+\epsilon+\kappa\delta+\eta,\Theta,\Upsilon)
\geq
f(\alpha,n,\epsilon,\Phi,\Psi)m^{-(c_1(n)+2nC_2(\Psi)\delta+C(\Upsilon))/
 |\Upsilon_{\max}|}.$$
This yields
 $
 \Lambda_{\Phi}^{\Psi}(\alpha,\epsilon)\leq \frac{2C_2(\Psi)\log m}{|\Upsilon_{\max}|}\delta+(1-C_2(\Psi)\delta)
 \Lambda_{\Theta}^{\Upsilon}(\beta, a+\epsilon+\kappa\delta+\eta).
 $
 Letting $\epsilon\downarrow 0$ and then
$a\downarrow a_0$ and $\delta\downarrow \delta_0$ we get
\begin{eqnarray*}
\Lambda_{\Phi}^{\Psi}(\alpha)\leq \frac{2C_2(\Psi)\log
m}{|\Upsilon_{\max}|}\delta_0+(1-C_2(\Psi)\delta_0)\Lambda_{\Theta}^{\Upsilon}(\beta,
(a_0+\kappa\delta_0+\eta)+2\eta).
%\\
%&\leq& \frac{2C_2(\Psi)\log
%m}{|\Upsilon_{\max}|}\delta_0+(1-C_2(\Psi)\delta_0)\Lambda_{\Theta}^{\Upsilon}(\beta,
%a_0+\kappa\delta_0+2\eta).
\end{eqnarray*}

\noindent
{\bf Proof of Lemma~\ref{appro-large-de}.} (1) For $\lambda\in \R$ define $f_n(\lambda):=P(\langle
z,\Phi^n-\alpha\rangle+\lambda\Psi^n)$ and $f(\lambda):=P(\langle
z,\Phi-\alpha\rangle+\lambda\Psi)$. Since
$
|\langle z,\phi^n_k(x)-k\alpha\rangle+\lambda\psi^n_k(x)-(\langle
z,\phi_k(x)-k\alpha\rangle+\lambda\psi_k(x))|\leq
|z|\|\phi^n_k-\phi_k\|+|\lambda|\|\psi^n_k-\psi_k\|
$, we have
$|f_n(\lambda)-f(\lambda)|\leq |z|\cdot\|\Phi^n-\Phi\|_{\mbox{\tiny \rm lim}}+|\lambda|\cdot\|\Psi^n-\Psi\|_{\mbox{\tiny \rm lim}}.$
Thus $f_n$ converges uniformly to $f$ over any bounded interval $I$.  By Lemma \ref{pre-sloope}, $f_n(\lambda)=0$ and
$f(\lambda)=0$ have unique solutions. Assume $f(\lambda_0)=0$ and
$f_n(\lambda_n)=0$. Then we have $\lambda_n\to \lambda_0$, i.e. $
{\tau_\Phi^\Psi}(z,\alpha,n)\to {\tau_\Phi^\Psi}(z,\alpha). $

Since ${\tau_\Phi^\Psi}(0,\alpha,n)=D(\Psi^n)$ and we have shown
that $D(\Psi^n)\to D(\Psi)$, we get
$
{\tau_\Phi^\Psi}(0,\alpha)=\lim_{n\to\infty}{\tau_\Phi^\Psi}(0,\alpha,n)=D(\Psi).
$, and then ${\tau_\Phi^\Psi}^\star(\alpha)\leq
{\tau_\Phi^\Psi}(0,\alpha)=D(\Psi).$

Now assume $\alpha\in \mathrm{int}(L_\Phi)$.  Since
$\|\Phi^n-\Phi\|_{\mbox{\tiny \rm lim}}\to 0$, it is ready to show
that $d_H(L_{\Phi^n},L_\Phi)\to 0.$ Since $L_{\Phi^n}$ and $L_\Phi$
are all compact convex sets and $L_\Phi$ has dimension $d$, it is
seen that $L_{\Phi^n}$ has nonempty interior for large $n.$ Moreover
we can find $N\in \N$ and $\delta_0>0$, such that
 $B(\alpha,\delta_0)\subset L_{\Phi^n}$ for any $n\geq N.$
By Lemma \ref{case-holder}(4) for any $z\in \R^d$ and any $n\geq N$
${\tau_\Phi^\Psi}(z,\alpha,n)\geq
\delta_0C_1(\Psi^n)|z|/2=\delta_0|z|/2|\Psi_{\min}^n|\geq
\delta_0|z|/2|\Psi_{\min}|.$
Letting $n$ tend to $\infty$ and then $|z|$ to $\infty$ we get
$\lim_{|z|\to\infty}{\tau_\Phi^\Psi}(z,\alpha)=+\infty.$
Thus we can find a $z_0\in \R^d$ such that
${\tau_\Phi^\Psi}^\star(\alpha)={\tau_\Phi^\Psi}(z_0,\alpha)$. By
Lemma \ref{case-holder}(3) ${\tau_\Phi^\Psi}(\cdot,\alpha,n)$ is
convex. By a well known theorem in convex analysis \cite{R} (p. 90),
${\tau_\Phi^\Psi}(\cdot,\alpha)$ is convex. Moreover the convergence
is uniform on any compact domain. Now by the uniform convergence of
${\tau_\Phi^\Psi}(z,\alpha,n)$ to ${\tau_\Phi^\Psi}(z,\alpha)$ over
the closed ball $B(z_0,R) $ with $R>0$ large enough, we can easily
show that ${\tau_\Phi^\Psi}^\star(\alpha,n)\to
{\tau_\Phi^\Psi}^\star(\alpha)$ as $n\to\infty$.

(2) Since $\Phi^n$ and $\Psi^n$ are H\"{o}lder continuous, by Lemma
\ref{case-holder}(2), if $\langle z,\alpha^\prime-\alpha\rangle\geq
0$, then
$$
C_1^n\langle z,\alpha^\prime-\alpha\rangle\leq
{\tau_\Phi^\Psi}(z,\alpha,n)-{\tau_\Phi^\Psi}(z,\alpha^\prime,n)\leq
C_2^n\langle z,\alpha^\prime-\alpha\rangle.
$$
where $C_1^n=1/|\Psi_{\min}^n|$ and $C_2^n=1+1/|\Psi_{\max}^n|$. By
Lemma \ref{control}, we know that $|\Psi_{\min}|\geq
|\Psi_{\min}^n|$ and $|\Psi_{\max}|\leq |\Psi_{\max}^n|$. Since
${\tau_\Phi^\Psi}(z,\cdot,n)\to{\tau_\Phi^\Psi}(z,\cdot)$, letting
$n\to\infty$ we get the result. \hfill $\Box$

Now we complete the proof of the Proposition \ref{gene-thermo}.

\noindent{\bf Proof of Proposition \ref{gene-thermo}(Continued)}\ We
prove it in two steps:

(1) For $\Phi$ H\"older and $\Psi\in\C_{aa}^-(\Sigma_A,T)$ the result
holds: Let
$$
c(\Psi,\Phi,n):=\sum_{w\in {\mathcal B}_n(\Psi)} \exp(\sup_{x\in
[w]}\langle z,\phi_{|w|}(x)\rangle).
$$
Fix $\Upsilon$ a H\"older potential and
$\delta>\|\Psi-\Upsilon\|_{\text{lim}}$. For any  $w\in{\mathcal
B}_n(\Psi)$  the proof of Lemma~\ref{Psi-dependence} yields $u\prec w$ such that
$u\in {\mathcal B}_{[n(1-C_2\delta)]}(\Upsilon)$ and $|w|-|u|\leq
C_2(o(n)+n\delta).$ Thus we get
\begin{eqnarray*}
c(\Psi,\Phi,n)
%&=&\sum_{w\in {\mathcal B}_n(\Psi)} \exp(\sup_{x\in
%[w]}\langle z,\phi_{|w|}(x)\rangle)\\
&=&\sum_{u\in {\mathcal B}_{[n(1-C_2\delta)]}(\Upsilon)}\ \
\sum_{v:uv\in{\mathcal B}_n(\Psi) } \exp(\sup_{x\in
[uv]}\langle z,\phi_{|uv|}(x)\rangle)\\
&\leq& (m^{C_2}e^{C(\Phi,z)})^{(o(n)+n\delta)}\sum_{u\in {\mathcal
B}_{[n(1-C_2\delta)]}(\Upsilon)}\exp(\sup_{x\in
[u]}\langle z,\phi_{|u|}(x)\rangle)\\
&=&C_3^{o(n)+n\delta}c(\Upsilon,\Phi,[n(1-C_2\delta])),
\end{eqnarray*}
where $C(\Phi,z)$ is a constant depending on $\Phi$ and $z$ only). Similarly we can get $c(\Psi,\Phi,n)\geq
C_4^{o(n)+n\delta}c(\Upsilon,\Phi,[n(1-C_2\delta])).$
%From this we get
%\begin{eqnarray*}
% &&\delta \log C_4+(1-C_2\delta)\tau_\Phi^\Upsilon(z,0)\leq
%\liminf_{n\to\infty}\frac{\log c(\Psi,\Phi,n)}{n}\\
%&\leq& \limsup_{n\to\infty}\frac{\log c(\Psi,\Phi,n)}{n}\leq\delta
%\log C_3+(1-C_2\delta)\tau_\Phi^\Upsilon(z,0).
%\end{eqnarray*}
Since $\|\Psi-\Upsilon\|_{\text{lim}}$ and hence $\delta$ can be taken arbitrarily small, this yields
$\tau_{\Phi}^\Psi(z,0)=\lim_{n\to\infty}\ln c(\Psi,\Phi,n)/n.$

(2) For general $\Phi$ and general $\Psi$ the result holds: Indeed,
once the previous step is established, by taking a sequence of
H\"older potentials $\Phi^j$ such that $\|\Phi^j-\Phi\|_{\rm lim}\to
0$ one can easily conclude. \hfill $\Box$

%%%%%%%%%%%%%%%%%%%%%%%%%%%%%%%%%%%%%%%%%%%%%%%%%%%%%%%%%%%%%%%%%%%%%%%%%%%%%%%%%%%%%%55


\begin{thebibliography}{9999}

\bibitem{BaFe09} J.~Barral and D.\,J. Feng, Weighted thermodynamic formalism and applications, arXiv:math/0909.4247v1.

%\bibitem{BQ} J. Barral, Y.-H. Qu,

\bibitem{Bar96} L. Barreira. A non-additive thermodynamic formalism and applications to dimension theory
of hyperbolic dynamical systems, Ergod. Theory $\&$ Dynam. Sys., {\bf 16} (1996), 871--927.

\bibitem{B} L. Barreira. Nonadditive thermodynamic formalism: equilibrium and
Gibbs measures. Discrete Contin. Dyn. Syst., {\bf 16} (2006), no. 2,
279--305.

\bibitem{BD} L. Barreira, P. Doutor. Almost additive
multifractal analysis. J. Math. Pures Appl., {\bf 92} (2009), 1--17.

\bibitem{BS} L. Barreira, B. Saussol, Multifractal analysis of hyperbolic flows, {\it Comm. Math. Phys.}, {\bf  214}  (2000),  339--371.

\bibitem{BSS} L. Barreira, B. Saussol, J. Schmeling. Higher-dimensional
multifractal analysis. J. Math. Pures Appl., {\bf 81} (2002), no. 1,
67--91.


\bibitem{Bow75} R. Bowen,
Equilibrium states and the ergodic theory of Anosov diffeomorphisms. Springer Lecture Notes No. 470 Springer-Verlag, Berlin. (1975)

\bibitem{Bow79}
R. Bowen, Hausdorff dimension of quasicircles. {\it Inst. Hautes
\'{E}tudes Sci. Publ. Math.} No. {\bf 50} (1979), 11--25.

\bibitem{Brooks} R.L. Brooks, C.A.B. Smith, A.H. Stone and W.T. Tutte.  The Dissection of Rectangles into Squares, Duke Math. J., {\bf 7} (1940), 312--340.


\bibitem{BMP} G. Brown, G. Michon, J. Peyri\`ere, On the multifractal analysis of measures, J. Stat. Phys., {\bf 66}, (1992) 775--790.


\bibitem{CFH} Y.-L. Cao, D.\,J. Feng, W. Huang. The thermodynamic formalism
for sub-additive potentials. Discrete Contin. Dyn. Syst., {\bf 20} (2008),
no. 3, 639--657.


\bibitem{Collet} P. Collet, J.L. Lebowitz, A. Porzio, The dimension spectrum of some dynamical systems,
 Proceedings of the symposium on statistical mechanics of phase transitions---mathematical and physical
  aspects (Trebon, 1986). J. Statist. Phys., {\bf 47} (1987), 609--644.

\bibitem{Fal88} K.J. Falconer, A subadditive thermodynamic formalism for
mixing repellers. {\it J. Phys. A} {\bf 21} (1988), no. 14,
L737--L742.


%\bibitem{Falc} K.\,J. Falconer,
%Fractal Geometry: Mathematical Foundations and Applications, 2nd Edition. Wiley, 2003

\bibitem{FF} A.\,H. Fan, D.\,J. Feng, On the distribution of long-term time
averages on symbolic space. J. Statist. Phys., {\bf 99} (2000), no. 3-4,
813--856.

\bibitem{FFW} A.\,H. Fan, D.\,J. Feng, J. Wu, Recurrence, dimension and
entropy. J. London Math. Soc. (2) 64 (2001), no. 1, 229--244.

%\bibitem{FanLau99} A.\,H. Fan, K.\, S. Lau, Iterated function system and Ruelle operator, J. Math. Anal. Appl., {\bf 231} (1999), 319--344.

\bibitem{Fen03}D.\,J. Feng, Lyapounov exponents for products of matrices and multifractal analysis.
 Part I: positive matrices. {\em Isra\"el J. Math.} {\bf 138} (2003), 353--376.

\bibitem{Fen04}D.\,J. Feng,
The variational principle for products of non-negative matrices.
{\it  Nonlinearity} {\bf 17} (2004) 447--457.

\bibitem{FH}  D.\,J. Feng, W. Huang, Lyapunov spectrum of asymptotically sub-additive potentials. arXiv:0905.2680v1, to appear in {\it Commun. Math. Phys.}


\bibitem{FeLa02}
D.\,J. Feng, K.\,S. Lau, The pressure function for products of
non-negative matrices, {\it Math. Res. Lett.} {\bf 9} (2002),
363-378.


\bibitem{FLW}  D.\,J. Feng, K.\, S. Lau, J. Wu, Ergodic limits on the conformal
repellers. Adv. Math., {\bf 169} (2002), no. 1, 58--91.

 \bibitem{GP97}
D. Gatzouras, Y.  Peres, Invariant measures of full dimension for some expanding maps. {\it Ergod. Th. $\&$ Dynam. Sys.} {\bf 17} (1997), no. 1, 147--167.

\bibitem{K} M. Kesseb\"ohmer,  Large deviation for weak Gibbs measures and multifractal
spectra, {\it Nonlinearity } {\bf 14} (2001), 395-409.


\bibitem{KS04} M. Kesseb\"ohmer, B. Stratman, A multifractal formalism for growth rates and
applications to geometrically finite Kleinian groups, {\it Ergod. Th. $\&$ Dynam. Sys.}, {\bf 24} (2004), 141--170.

\bibitem{Ma} N.G. Makarov, Fine structure of harmonic measure. {\it St. Peterburg Math.
J.,} 10(2) (1999), 217--268.


\bibitem{Mum06} A. Mummert, The thermodynamic formalism for almost-additive sequences.  {
\it Discrete Contin. Dyn. Syst.} {\bf 16} (2006), 435--454.

%\bibitem{Pat97} N. Patzschke, Self-conformal multifractal measures, Adv. Applied Math., {\bf 19} (1997), 486--513.

\bibitem{O} E. Olivier, Multifractal analysis in symbolic dynamics and distribution of pointwise dimension for $g$-measures , {\it Nonlinearity},  {\bf 12} (1999), 1571--1585.


\bibitem{Olsen} L. Olsen, Multifractal analysis of divergence points of deformed measure theoretical Birkhoff averages,  {\it J. Math. Pures Appl.},  {\bf 82}  (2003), 1591--1649. 

\bibitem{PRSS} Peres, Y.; Rams, M.; Simon, K.; Solomyak, B.
Equivalence of positive Hausdorff measure and the open set condition
for self-conformal sets. Proc. Amer. Math. Soc. 129 (2001), no. 9,
2689--2699.

\bibitem{P} Y. Pesin, Dimension theory in dynamical systems. Contemporary
views and applications. Chicago Lectures in Mathematics. University
of Chicago Press, Chicago, IL, 1997.

\bibitem{PW} Y. Pesin, H. Weiss, A multifractal analysis of equilibrium
measures for conformal expanding maps and Moran-like geometric
constructions. J. Statist. Phys., {\bf 86} (1997), no. 1-2, 233--275.

\bibitem{PeW} Y. Pesin, H. Weiss, The multifractal analysis of Gibbs measures: Motivation, Mathematical Foundation, and Examples, Chaos, {\bf 7} (1997) 89--106.


\bibitem{Rand} D.A. Rand, The singularity spectrum $f(\alpha)$ for cookie-cutters, Ergodic Theory Dynam. Systems, {\bf  9} (1989), 527--541.

\bibitem{R}  Rockafellar, R.
 Convex analysis. Princeton Mathematical Series, No. 28 Princeton University Press, Princeton, N.J. 1970.

\bibitem{Ruelle}D. Ruelle, {\it Thermodynamic formalism.
The mathematical structures of classical equilibrium statistical
mechanics}. Encyclopedia of Mathematics and its Applications, 5.
Addison-Wesley Publishing Co., Reading, Mass., 1978.


\end{thebibliography}
\end{document}